\documentclass[12pt]{amsart}
\usepackage{fullpage}

\usepackage{graphics}
\usepackage[all]{xy}
\xyoption{v2}
\xyoption{knot}

\let\url=\undefined
\usepackage
[bookmarks,colorlinks,linkcolor=blue,citecolor=blue,urlcolor=blue]
{hyperref}

\title{Operads in iterated monoidal categories}
\author {Stefan Forcey, Jacob Siehler and E. Seth Sowers}

\thanks{Thanks to {\Xy-pic} for the diagrams. }

\address{Department of Physics and Mathematics\\
Tennessee State University\\
Nashville, TN 37209 \\
USA\\
}


\swapnumbers
\theoremstyle{plain}
\newtheorem{theorem}{Theorem}[section]
\newtheorem{lemma}[theorem]{Lemma}
\newtheorem{corollary}[theorem]{Corollary}

\theoremstyle{definition}
\newtheorem{definition}[theorem]{Definition}
\newtheorem{example}[theorem]{Example}

\theoremstyle{remark}
\newtheorem{remark}[theorem]{Remark}

\def\cal#1{\mathcal{#1}}
\newcommand{\bcal}[1]{\mbox{\boldmath${\cal {#1}}$}}

\def\nat{\mathbb{N}}
\def\th{^\text{th}}
\def\dash{\text{---}}
\DeclareMathOperator{\concat}{Concat}
\DeclareMathOperator{\seq}{Seq}
\DeclareMathOperator{\modseq}{ModSeq}
\DeclareMathOperator{\col}{Col}
\DeclareMathOperator{\mon}{Mon}
\DeclareMathOperator{\oper}{Oper}

\begin{document}

\begin{abstract}
The structure of a $k$-fold monoidal category as introduced by
Balteanu, Fiedor\-owicz, Schw\"anzl and Vogt in \cite{Balt} can be
seen as a weaker structure than a symmetric or even braided monoidal
category.  In this paper we show that it is still sufficient to
permit a good definition of ($n$-fold) operads in a $k$-fold monoidal
category which generalizes the definition of operads in a braided
category. Furthermore, the inheritance of structure by the category
of operads is actually an inheritance of iterated monoidal structure,
decremented by at least two iterations.  We prove that the category of $n$-fold operads in a
$k$-fold monoidal category is itself a $(k-n)$-fold monoidal, strict
$2$-category, and  show
that $n$-fold operads are automatically $(n-1)$-fold operads.
We also introduce a family of simple examples of
$k$-fold monoidal categories and classify operads in the example
categories.
\end{abstract}

\maketitle
\tableofcontents

\section{Introduction}
In this introductory section we will give a brief, non-chronological
overview of the relationship between operads, higher category theory,
and topology. This will serve to motivate the study of iterated
monoidal categories and their operads that comprises the remaining
sections.  In the second section, in order to be self contained,
we repeat the definition of the iterated monoidal categories first
set down in \cite{Balt}.  In the fourth section we seek to fill a
gap in the literature which currently contains few good examples
of that definition.  Thus our first new contribution consists of a
series of simple and very geometric iterated monoidal categories
based on totally ordered monoids. By simple we mean that axioms are
largely fulfilled due to relationships between max, plus, concatenation,
sorting and lexicographic ordering as well as the fact that all
diagrams commute since the underlying directed graph of the category
is merely the total order. The most interesting examples of $n$-fold
monoidal categories are those whose objects can be represented by
Ferrer or Young diagrams  (the underlying shapes of Young tableaux.)
These  exhibit products with the geometrical interpretation ``combining
stacks of boxes.'' Managers of  warehouses or quarries perhaps may
already be well aquainted with the three dimensional version of the
main example of iterated monoidal categories we introduce here.
Imagine that floor space in the quarry or warehouse is at a premium
and that therefore you are charged with combining several stacks
of crates or stone blocks by restacking some together vertically
and shifting others together horizontally. It turns out that the
best result in terms of gained floor space is always to be achieved
most efficiently by doing the restacking and shifting in  a very
particular order--horizontally first, then vertically.

The main new contribution is the theory of operads within, or
enriched in, iterated monoidal categories.  This theory is based
upon the fact that the natural setting of operads turns out to be
in a category with lax interchange between multiple operations, as
opposed to the full strength of a braiding or symmetry as is
classically assumed.  Batanin's definition of $n$-operad also relies
on this insight \cite{bat}. In that paper he notes that an iterated
monoidal category ${\cal V}$ would be an example of a globular
monoidal category with a single object, and a single arrow in each
dimension up to $n$, in which last dimension the arrows would
actually be the objects of ${\cal V}.$ Of course the invertibility
of the interchange would also have to be dropped from his definition.
In that case the $n$-fold operads defined here would correspond to
Batanin's  $n$-operads. The advantages of seeing them in a single
categorical dimension are in the way that doing so generalizes the
fact that operads in a symmetric monoidal category inherit its
symmetric structure.  We  investigate the somewhat flexible structure
of the iterated monoidal 2-category that $n$-fold operads comprise.
Flexibility arises from the difference between $n$ and $k$, where
one is investigating $n$-fold operads in a $k$-fold monoidal category
${\cal V}$, where $n < k-1$.  It turns out that choosing $n$ much
smaller than $k$ allows multiple interchanging products to be defined
on the category of operads, whereas choosing $n$ closer to $k$
allows the operads to take on multiple operad structures at once
with respect to the products in ${\cal V}.$ Examples of combinatorial
operads living in the previously introduced combinatorically defined
categories are utilized to demonstrate the sharpness of several of
the resulting theorems, i.e. to provide counterexamples. The examples
start to take on a life of their own, however, as theorems and open
questions about the classification of operads   in combinatorial
$n$-fold monoidal categories arise.  The definition of operad in
the categories with morphisms given by ordering leads to descriptions
of interesting kinds of growth.  We give a complete description of
the simple example of 2-fold operads in the natural numbers.  We
then give the elementary results for operads in the category of
Young diagrams.  In the basic examples linear and logarithmic growth
characterize respective dimensions in a single sequence of Young
diagrams.  These phenomena hint towards a theory of operadic growth.
Full investigation and further classification must await a future
sequel to this paper. Applications might be found in scientific
fields such as the theory of small world networks, where the diameter
of a network is the logarithm of the number of nodes.

First, however, we look at  some of the  history and philosophy of
the two major players here, operads and iterated monoidal categories.
Operads in a category of topological spaces are the crystallization
of several approaches to the recognition problem for iterated loop
spaces. Beginning with Stasheff's associahedra and Boardman and
Vogt's little $n$-cubes, and continuing with more general $A_{\infty}$,
$E_n$ and $E_{\infty}$ operads described by May and others, that
problem has largely been solved \cite{Sta}, \cite{BV1}, \cite{May}.
Loop spaces are characterized by admitting an operad action of the
appropriate kind. More lately Batanin's approach to higher categories
through internal and higher operads  promises to shed further light
on the combinatorics of $E_n$ spaces \cite{bat2}, \cite{bat3}.

Recently there has also been growing interest in the application
of higher dimensional structured categories to the characterization
of loop spaces. The program being advanced by many categorical
homotopy theorists seeks to model the coherence laws governing
homotopy types with the coherence axioms of structured  $n$-categories.
By modeling we mean a connection that will  be in the form of a
functorial equivalence between categories of special categories and
categories of special spaces.  The largest challenges currently are
to find the most natural and efficient definition of (weak)
$n$-category, and to determine the nature of the functor from
categories to spaces.  The latter will almost certainly be analogous
to the nerve functor on 1-categories, which preserves homotopy
equivalence.  In \cite{StAlg} Street defines the nerve of a strict
$n$-category.  Recently Duskin in \cite{Dusk} has worked out the
description of the nerve of a bicategory.  A second part of the
latter paper promises the full description  of the functor including
how it takes morphisms of bicategories to continuous maps.

One major recent advance is the discovery of Balteanu, Fiedorowicz,
Schw\"anzl and Vogt in \cite{Balt} that the nerve functor
on categories gives a direct connection between iterated monoidal
categories and iterated loop spaces.  Stasheff \cite{Sta} and Mac
Lane \cite{Mac} showed that monoidal categories are precisely
analogous to 1-fold loop spaces. There is a similar connection
between symmetric monoidal categories and infinite loop spaces. The
first step in filling in the gap between 1 and infinity was made
in \cite{ZF} where it is shown that the group completion of the
nerve of a braided monoidal category is a 2-fold loop space.   In
\cite{Balt} the authors finish this process by, in their words,
``pursuing an analogy to the tautology that an $n$-fold loop space
is a loop space in the category of $(n-1)$-fold loop spaces.'' The
first thing they focus on is the fact that a braided category is a
special case of a carefully defined 2-fold monoidal category. Based
on their observation of the  correspondence between loop spaces and
monoidal categories, they iteratively define the notion of $n$-fold
monoidal category as a monoid in the category of $(n-1)$-fold
monoidal categories.  In \cite{Balt} a symmetric category is seen
as a category that is $n$-fold monoidal for all $n$.  The main
result in that paper is that the group completion of the nerve of
an $n$-fold monoidal category is an $n$-fold loop space.  It is
still open whether this is a complete characterization, that is,
whether every $n$-fold loop space arises as the nerve of an $n$-fold
monoidal category. Much progress towards the answer to this question
was made by the original authors in their sequel paper, but the
desired result was later shown to remain unproven. One of the future
goals of the program begun here is to use weakenings or deformations
of the examples of $n$-fold monoidal categories introduced here to
model specific loop spaces in a direct way.

The connection between the $n$-fold monoidal categories of Fiedorowicz
and the theory of higher categories is through the  periodic table
as laid out in \cite{Baez1}.  Here Baez organizes the $k$-tuply
monoidal $n$-categories, by which terminology he refers to
$(n+k)$-categories  that are trivial below dimension $k.$ The
triviality of lower cells allows the higher ones to compose freely,
and thus these special cases of $(n+k)$-categories are viewed as
$n$-categories with $k$ multiplications.  Of course a $k$-tuply
monoidal $n$-category is a special $k$-fold monoidal $n$-category.
The specialization results from the definition(s) of $n$-category,
all of which seem to include the axiom that the interchange
transformation between two ways of composing four higher morphisms
along two different lower dimensions is required to be an isomorphism.
As will be mentioned in the next section the  property of having
iterated loop space nerves held by the $k$-fold monoidal categories
relies on interchange transformations that are not isomorphisms.
If those transformations are indeed isomorphisms then the $k$-fold
monoidal 1-categories do reduce to the braided and symmetric
1-categories of the periodic table. Whether this continues for
higher dimensions, yielding for example the sylleptic monoidal
2-categories of the periodic table as 3-fold monoidal 2-categories
with interchange isomorphisms, is yet to be determined.

A further refinement of higher categories is to require all morphisms
to have inverses. These special cases are referred to as $n$-groupoids,
and since their nerves are simpler to describe it has long been
suggested that they model homotopy $n$-types through a construction
of a fundamental $n$-groupoid. This has in fact been shown to hold
in Tamsamani's definition of weak $n$-category \cite{tam}, and in
a recent paper by Cisinski to hold in the definition of Batanin as
found in \cite{bat}.  A homotopy $n$-type is a topological space
$X$ for which $\pi_k(X)$ is trivial for all $k>n.$ Thus the homotopy
$n$-types are classified by $\pi_k$ for $k\le n.$ It has been
suggested that a key requirement for any useful definition of
$n$-category is that a $k$-tuply monoidal $n$-groupoid be associated
functorially (by a nerve) to a topological space which is a  homotopy
$n$-type and a $k$-fold loop space \cite{Baez1}.  The loop degree
will be precise for $k<n+1,$ but for $k>n$ the associated homotopy
$n$-type will be an infinite loop space.  This last statement is a
consequence of the stabilization hypothesis , which states that
there should be a left adjoint to forgetting monoidal structure
that is an equivalence of $(n+k+2)$-categories between $k$-tuply
monoidal $n$-categories and $(k+1)$-tuply monoidal $n$-categories
for $k>n+1.$ This hypothesis has been shown by Simpson to hold in
the case of Tamsamani's definition \cite{sim}. For the case of $n=1$
if the interchange transformations are isomorphic then a $k$-fold
monoidal 1-category is equivalent to a symmetric category for $k>2.$
With these facts in mind it is possible that if we wish to precisely
model homotopy $n$-type $k$-fold loop spaces for $k>n$ then we may
need to consider $k$-fold as well as $k$-tuply monoidal $n$-categories.
This paper is part of an embryonic  effort in that direction.

Since a loop space can be efficiently described as an operad algebra,
it is not surprising that there are several existing definitions
of $n$-category that utilize operad actions.  These definitions
fall into two main classes: those that define an $n$-category as
an algebra of a higher order operad, and those that achieve an
inductive definition using classical operads in symmetric monoidal
categories to parameterize iterated enrichment. The first class of
definitions is typified by Batanin and Leinster \cite{bat},\cite{lst}.

The former author defines monoidal globular categories in which
interchange transformations are isomorphisms and which thus resemble
free strict $n$-categories.  Globular operads live in these, and
take all sorts of pasting diagrams as input types, as opposed to
just a string of objects as in the case of classical operads.  The
binary composition in an $n$-category derives from the action of a
certain one  of these globular operads.  Leinster expands this
concept to describe $n$-categories with unbiased composition of any
number of cells.  The second class of definitions is typified by
the works of Trimble and May \cite{trimble}, \cite{may2}.

The former parameterizes iterated enrichment with a series of operads
in $(n-1)$-Cat achieved by taking the fundamental $(n-1)$-groupoid
of the $k\th$ component of the topological path composition operad
$E.$ The latter begins with an $A_{\infty}$ operad in a symmetric
monoidal category ${\cal V}$ and requires his enriched categories
to be tensored over ${\cal V}$ so that the iterated enrichment
always refers to the same original operad.

Iterated enrichment over $n$-fold categories is described in
\cite{forcey1} and \cite{forcey2}.  It seems wothwhile to define
$n$-fold operads in $n$-fold monoidal categories in a way that is
consistent with the spirit of Batanin's globular operads.  Their
potential value may include using them to weaken enrichment over
$n$-fold monoidal categories in a way that is in the spirit of May
and Trimble.  This program carries with it the promise of characterizing
$k$-fold loop spaces with homotopy $n$-type for all $n,k$ by
describing  the categories with exactly those spaces as nerves.  As
a candidate for the type of category with such a nerve we suggest
a weak $n$-category with $k$ multiplications that interchange only
in the lax sense.

In this paper we follow May by  defining $n$-fold operads in terms
of monoids in a certain category of collections.  A more abstract
approach for future consideration would begin by finding an equivalent
definition in the language of Weber, where an operad lives within
a monoidal pseudo algebra of a 2-monad \cite{web}.  This latter is
a general notion of operad which includes as instances both the
classical operads and the higher operads of Batanin.

\section{\texorpdfstring{$k$}{K}-fold monoidal categories}
This sort of category was developed and defined in \cite{Balt}. The
authors describe its structure as arising recursively from its
description as a monoid in the category of $(k-1)$-fold monoidal
categories. Here we present that definition (in its expanded form)
altered only slightly to make visible the coherent associators as
in \cite{forcey1}.  That latter paper describes its structure in
terms of tensor objects in the category of $(k-1)$-fold monoidal
categories. Our variation has the effect of making visible the
associators $\alpha^{i}_{ABC}.$ It is desirable to do so for several
reasons. One is that this makes easier a direct comparison with
Batanin's definition of monoidal globular categories as in \cite{bat}.
A monoidal globular category can be seen as a quite special case
of a iterated monoidal category, with source and target maps that
take objects to those in a category with one less product, and with
interchanges that are isomorphisms.

The other reason is that in this paper we will consider a category
of collections in an iterated monoidal category which will be
(iterated) monoidal only up to natural associators.  That being
said, in much of the remainder of this paper we will consider
examples with strict associativity, where each $\alpha$ is the
identity, and in interest of clarity will often hide associators.
One actual simplification in the following definition is that all
the products are assumed to have the same unit. We note that this
is easily generalized, as in the case of collections which we will
consider.

\begin{definition}\label{iterated} An $n${\it -fold monoidal category}
is a category ${\cal V}$ with the following structure.
\begin{enumerate}
\item There are $n$ multiplications
\[\otimes_1,\otimes_2,\dots, \otimes_n:{\cal V}\times{\cal V}\to{\cal V}\]
each equipped with an associator $\alpha_{UVW}$, a natural isomorphism
which satisfies the pentagon equation:

\noindent
\begin{center}
\resizebox{4.5in}{!}{
$$
\xymatrix@R=45pt@C=-35pt{
&((U\otimes_i V)\otimes_i W)\otimes_i X \text{ }\text{ }
\ar[rr]^{ \alpha^{i}_{UVW}\otimes_i 1_{X}}
\ar[dl]^{ \alpha^{i}_{(U\otimes_i V)WX}}
&\text{ }\text{ }\text{ }\text{ }\text{ }&\text{ }\text{ }(U\otimes_i (V\otimes_i W))\otimes_i X
\ar[dr]^{ \alpha^{i}_{U(V\otimes_i W)X}}&\\
(U\otimes_i V)\otimes_i (W\otimes_i X)
\ar[drr]^{ \alpha^{i}_{UV(W\otimes_i X)}}
&&&&U\otimes_i ((V\otimes_i W)\otimes_i X)
\ar[dll]^{ 1_{U}\otimes_i \alpha^{i}_{VWX}}
\\&&U\otimes_i (V\otimes_i (W\otimes_i X))&&&
}
$$
}
\end{center}

\item ${\cal V}$ has an object $I$ which is a strict unit
for all the multiplications.

\item For each pair $(i,j)$ such that $1\le i<j\le n$ there is a
natural transformation
\[\eta^{ij}_{ABCD}: (A\otimes_j B)\otimes_i(C\otimes_j D)
    \to (A\otimes_i C)\otimes_j(B\otimes_i D).\]
These natural transformations $\eta^{ij}$ are subject to the following
conditions:
\begin{enumerate}
\item Internal unit condition:
    $\eta^{ij}_{ABII}=\eta^{ij}_{IIAB}=1_{A\otimes_j B}$
\item External unit condition:
    $\eta^{ij}_{AIBI}=\eta^{ij}_{IAIB}=1_{A\otimes_i B}$
\item Internal associativity condition: The following diagram commutes.

\noindent
\begin{center}
\resizebox{5in}{!}{
$$
\diagram
((U\otimes_j V)\otimes_i (W\otimes_j X))\otimes_i (Y\otimes_j Z)
\xto[rrr]^{\eta^{ij}_{UVWX}\otimes_i 1_{Y\otimes_j Z}}
\ar[d]^{\alpha^i}
&&&\bigl((U\otimes_i W)\otimes_j(V\otimes_i X)\bigr)\otimes_i (Y\otimes_j Z)
\dto^{\eta^{ij}_{(U\otimes_i W)(V\otimes_i X)YZ}}\\
(U\otimes_j V)\otimes_i ((W\otimes_j X)\otimes_i (Y\otimes_j Z))
\dto^{1_{U\otimes_j V}\otimes_i \eta^{ij}_{WXYZ}}
&&&((U\otimes_i W)\otimes_i Y)\otimes_j((V\otimes_i X)\otimes_i Z)
\ar[d]^{\alpha^i \otimes_j \alpha^i}
\\
(U\otimes_j V)\otimes_i \bigl((W\otimes_i Y)\otimes_j(X\otimes_i Z)\bigr)
\xto[rrr]^{\eta^{ij}_{UV(W\otimes_i Y)(X\otimes_i Z)}}
&&& (U\otimes_i (W\otimes_i Y))\otimes_j(V\otimes_i (X\otimes_i Z))
\enddiagram
$$
}
\end{center}

\item External associativity condition: The following diagram commutes.

\noindent
\begin{center}
\resizebox{5in}{!}{
$$
\diagram
((U\otimes_j V)\otimes_j W)\otimes_i ((X\otimes_j Y)\otimes_j Z)
\xto[rrr]^{\eta^{ij}_{(U\otimes_j V)W(X\otimes_j Y)Z}}
\ar[d]^{\alpha^j \otimes_i \alpha^j}
&&& \bigl((U\otimes_j V)\otimes_i (X\otimes_j Y)\bigr)\otimes_j(W\otimes_i Z)
\dto^{\eta^{ij}_{UVXY}\otimes_j 1_{W\otimes_i Z}}\\
(U\otimes_j (V\otimes_j W))\otimes_i (X\otimes_j (Y\otimes_j Z))
\dto^{\eta^{ij}_{U(V\otimes_j W)X(Y\otimes_j Z)}}
&&&((U\otimes_i X)\otimes_j(V\otimes_i Y))\otimes_j(W\otimes_i Z)
\ar[d]^{\alpha^j}
\\
(U\otimes_i X)\otimes_j\bigl((V\otimes_j W)\otimes_i (Y\otimes_j Z)\bigr)
\xto[rrr]^{1_{U\otimes_i X}\otimes_j\eta^{ij}_{VWYZ}}
&&& (U\otimes_i X)\otimes_j((V\otimes_i Y)\otimes_j(W\otimes_i Z))
\enddiagram
$$
}
\end{center}

\item Finally it is required  for each triple $(i,j,k)$ satisfying
$1\le i<j<k\le n$ that the giant hexagonal interchange diagram commutes.

\noindent
\begin{center}
\resizebox{5.5in}{!}{
$$
\xymatrix@C=-118pt{
&((A\otimes_k A')\otimes_j (B\otimes_k B'))\otimes_i((C\otimes_k C')\otimes_j (D\otimes_k D'))
\ar[ddl]|{\eta^{jk}_{AA'BB'}\otimes_i \eta^{jk}_{CC'DD'}}
\ar[ddr]|{\eta^{ij}_{(A\otimes_k A')(B\otimes_k B')(C\otimes_k C')(D\otimes_k D')}}
\\\\
((A\otimes_j B)\otimes_k (A'\otimes_j B'))\otimes_i((C\otimes_j D)\otimes_k (C'\otimes_j D'))
\ar[dd]|{\eta^{ik}_{(A\otimes_j B)(A'\otimes_j B')(C\otimes_j D)(C'\otimes_j D')}}
&&((A\otimes_k A')\otimes_i (C\otimes_k C'))\otimes_j((B\otimes_k B')\otimes_i (D\otimes_k D'))
\ar[dd]|{\eta^{ik}_{AA'CC'}\otimes_j \eta^{ik}_{BB'DD'}}
\\\\
((A\otimes_j B)\otimes_i (C\otimes_j D))\otimes_k((A'\otimes_j B')\otimes_i (C'\otimes_j D'))
\ar[ddr]|{\eta^{ij}_{ABCD}\otimes_k \eta^{ij}_{A'B'C'D'}}
&&((A\otimes_i C)\otimes_k (A'\otimes_i C'))\otimes_j((B\otimes_i D)\otimes_k (B'\otimes_i D'))
\ar[ddl]|{\eta^{jk}_{(A\otimes_i C)(A'\otimes_i C')(B\otimes_i D)(B'\otimes_i D')}}
\\\\
&((A\otimes_i C)\otimes_j (B\otimes_i D))\otimes_k((A'\otimes_i C')\otimes_j (B'\otimes_i D'))
}
$$
}
\end{center}
\end{enumerate}     
\end{enumerate}     
\end{definition}

Note that for $q>p$ we have natural transformations
\[
\eta^{pq}_{AIIB}:A\otimes_p B\to A\otimes_q B\qquad\text{ and }\qquad
\eta^{pq}_{IABI}:A\otimes_p B\to B\otimes_q A.
\]

If the authors of  \cite{Balt} had insisted a 2-fold monoidal
category be a tensor object in the category of monoidal categories
and {\it strictly monoidal\/} functors, this would be equivalent
to requiring that $\eta=1$.  In view of the above, they note that
this would imply $A\otimes_1 B = A\otimes_2 B = B\otimes_1 A$ and
similarly for morphisms.

Joyal and Street \cite{JS} considered a similar concept to Balteanu,
Fiedorowicz, Schw\"anzl and Vogt's idea of 2--fold
monoidal category.  The former pair required the natural transformation
$\eta_{ABCD}$ to be an isomorphism and showed that the resulting
category is naturally equivalent to a braided monoidal category.
As explained in \cite{Balt}, given such a category one obtains an
equivalent braided monoidal category by discarding one of the two
operations, say $\otimes_2$, and defining the commutativity isomorphism
for the remaining operation $\otimes_1$ to be the composite
\[
\diagram
A\otimes_1 B\rrto^{\eta_{IABI}}
&& B\otimes_2 A\rrto^{\eta_{BIIA}^{-1}}
&& B\otimes_1 A.
\enddiagram
\]

The authors of \cite{Balt} remark that a symmetric monoidal category
is $n$-fold monoidal for all $n$.  This they demonstrate by letting
\[\otimes_1=\otimes_2=\dots=\otimes_n=\otimes\]
and defining
\[\eta^{ij}_{ABCD}=\alpha^{-1}\circ (1_A\otimes \alpha)
    \circ(1_A\otimes (c_{BC}\otimes 1_D))\circ
    \circ(1_A\otimes \alpha^{-1})\circ \alpha \]
for all $i<j$. Here $c_{BC}: B\otimes C \to C\otimes B$ is the symmetry natural transformation.

Joyal and Street \cite{JS} require that the interchange natural
transformations $\eta^{ij}_{ABCD}$ be isomorphisms and observed
that for $n\ge3$ the resulting sort of category is equivalent to a
symmetric monoidal category. Thus as Balteanu, Fiedorowicz, Schw\"anzl
and Vogt point out, the nerves of such categories have group
completions which are infinite loop spaces rather than only $n$--fold
loop spaces.

Because of the recursive nature of the definition of iterated
monoidal category, there are multiple forgetful functors implied.
Specifically, letting $n<k$, from the category of $k$-fold monoidal
categories to the category of $n$-fold monoidal categories there
are $\binom k n$ forgetful functors which forget all but the chosen
set of products.

The coherence theorem for iterated monoidal categories states  that
any diagram composed solely of interchange transformations commutes;
i.e.  if two compositions of various interchange transformations
(legs of a diagram) have the same source and target then they
describe the same morphism. Furthermore we can easily determine
when a composition of interchanges exists between objects. Here are
the necessary definitions and Theorem as given in \cite{Balt}.

\begin{definition}
Let ${\cal F}_n(S)$ be the free $n$-fold monoidal category on the
finite set $S.$ Its objects are all finite expressions generated
by the elements of $S$ using the products $\otimes_i, i=1..n.$ By
${\cal M}_n(S)$ we denote the sub-category of ${\cal F}_n(S)$ whose
objects are expressions in which each element of $S$ occurs exactly
once.
\end{definition}

If $S \subset T$ then there is a restriction functor ${\cal M}_n(T)
\to {\cal M}_n(S)$, induced by the functor ${\cal F}_n(T) \to {\cal
F}_n(S)$, which sends $T - S$ to the empty expression 0.

\begin{definition}
Let $A$ be an object of ${\cal M}_n(S).$ For $a,b \in S$ we say
that $a\otimes_ib$ $in$ $A$ if the restriction functor ${\cal M}_n(S)
\to {\cal M}_n({a,b})$ sends $A$ to $a\otimes_ib.$
\end{definition}

\begin{theorem}\cite{Balt}
Let $A$ and $B$ be objects of ${\cal M}_n(S).$ Then
\begin{enumerate}
\item There is at most one morphism $A\to B.$
\item Moreover, there exists a morphism $A\to B$ if and only if, for
every $a,b \in S$, $a\otimes_i b\in A$ implies that either
$a\otimes_j b$ or $b\otimes_j a$ is in $B$ for some $j> i.$
\end{enumerate}
\end{theorem}

\section{\texorpdfstring{$n$}{N}-fold operads}
The two principle components of an operad are a collection,
historically a sequence, of objects in a monoidal category and a
family of composition maps. Operads are often described as
parameterizations  of $n$-ary operations.  Peter May's original
definition of operad in a symmetric (or braided) monoidal category
\cite{May} has a composition $\gamma$ that takes the tensor product
of the $n\th$ object ($n$-ary operation) and $n$ others (of various
arity) to a resultant that sums the arities of those others.  The
$n\th$ object or $n$-ary operation is often pictured as a tree with
$n$ leaves, and the composition appears like this:

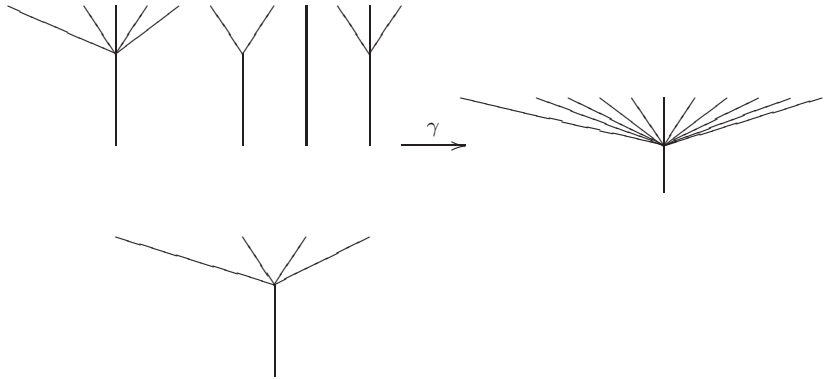
\begin{tabular}{ll}
$$
\xymatrix@W=2.2pc @H=2.2pc @R=1.5pc @C=1pc{
*=0{}\ar@{-}[drr] &*=0{}\ar@{-}[dr] &*=0{}\ar@{-}[d] &*=0{}\ar@{-}[dl] &*=0{}\ar@{-}[dll] &*=0{}\ar@{-}[dr] &*=0{} &*=0{}\ar@{-}[dl] &*=0{}\ar@{-}[d] &*=0{}\ar@{-}[dr] &*=0{}\ar@{-}[d] &*=0{}\ar@{-}[dl]\\
*=0{} &*=0{} &*=0{}\ar@{-}[d] &*=0{} &*=0{} &*=0{} &*=0{}\ar@{-}[d] &*=0{} &*=0{}\ar@{-}[d] &*=0{}  &*=0{}\ar@{-}[d] &*=0{} &*=0{}&*=0{}&*=0{}&*=0{}&*=0{}&*=0{}\\&*=0{}&*=0{}&*=0{}&*=0{}&*=0{}&*=0{}&*=0{}&*=0{}&*=0{}&*=0{}&*=0{}\ar[rr]^{\gamma}&*=0{}&*=0{}&*=0{}&*=0{}\\
*=0{} &*=0{} &*=0{}\ar@{-}[drrrrr] &*=0{} &*=0{} &*=0{}&*=0{}\ar@{-}[dr] &*=0{} &*=0{}\ar@{-}[dl] &*=0{} &*=0{}\ar@{-}[dlll] &*=0{} &*=0{}&*=0{}&*=0{}&*=0{}&*=0{}&*=0{}&*=0{}&*=0{}&*=0{}\\
*=0{} &*=0{} &*=0{} &*=0{}      &*=0{} &*=0{} &*=0{}&*=0{}\ar@{-}[d] &*=0{} &*=0{}\\ &*=0{} &*=0{} &*=0{} &*=0{} &*=0{} &*=0{} &*=0{} &*=0{}&*=0{}&*=0{}&*=0{}
}
$$ & \hspace{-10pc}
$$
\xymatrix@W=2.2pc @H=2.2pc @R=1.5pc @C=1pc{\\
*=0{}\ar@{-}[drrrrr]&*=0{}\ar@{-}[drrrr]&*=0{}\ar@{-}[drrr]&*=0{}\ar@{-}[drr]&*=0{}\ar@{-}[dr]&*=0{}\ar@{-}[d]&*=0{}\ar@{-}[dl]&*=0{}\ar@{-}[dll]&*=0{}\ar@{-}[dlll]&*=0{}\ar@{-}[dllll]&*=0{}\ar@{-}[dlllll]\\
*=0{}&*=0{}&*=0{}&*=0{}&*=0{}&*=0{}\ar@{-}[d]\\
*=0{}*=0{}*=0{}*=0{}*=0{}*=0{}*=0{}&*=0{}&*=0{}&*=0{}&*=0{}&*=0{}&*=0{}&*=0{}&*=0{}&*=0{}\\
}
$$
\end{tabular}

By requiring this composition to be associative we mean that it obeys this sort of pictured commuting diagram:

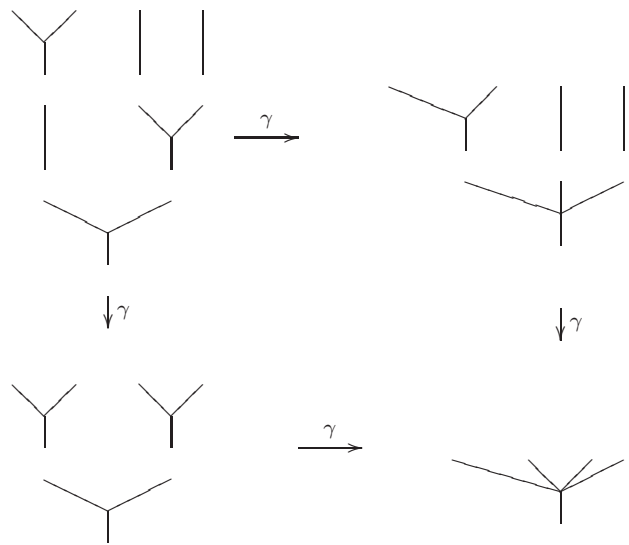
\begin{center}
\begin{tabular}{ll}
$$
\xymatrix@W=2.2pc @H=2.2pc @R=1pc @C=1pc{
*=0{}\ar@{-}[dr] &*=0{} &*=0{}\ar@{-}[dl] &*=0{} &*=0{}\ar@{-}[d] &*=0{} &*=0{}\ar@{-}[d] \\
*=0{} &*=0{}\ar@{-}[d] &*=0{} &*=0{} &*=0{}\ar@{-}[d] &*=0{} &*=0{}\ar@{-}[d] \\
*=0{}&*=0{}&*=0{}&*=0{}&*=0{}&*=0{}&*=0{}&*=0{}&*=0{}\\
*=0{} &*=0{}\ar@{-}[d] &*=0{} &*=0{} &*=0{}\ar@{-}[dr] &*=0{} &*=0{}\ar@{-}[dl] \\
*=0{} &*=0{}\ar@{-}[d] &*=0{} &*=0{} &*=0{} &*=0{}\ar@{-}[d] &*=0{}&*=0{}\ar[rr]^{\gamma}&*=0{}&*=0{}&*=0{}&*=0{} \\
*=0{}&*=0{}&*=0{}&*=0{}&*=0{}&*=0{}&*=0{}&*=0{}&*=0{}\\
*=0{} &*=0{}\ar@{-}[drr]  &*=0{} &*=0{} &*=0{} &*=0{}\ar@{-}[dll]\\
*=0{} &*=0{} &*=0{} &*=0{}\ar@{-}[d] \\
*=0{} &*=0{} &*=0{} &*=0{} \\
*=0{} &*=0{} &*=0{} &*=0{}\ar[d]^{\gamma}&*=0{}\\
*=0{} &*=0{} &*=0{} &*=0{} &*=0{}
}
$$ & \hspace{-4pc}
$$
\xymatrix@W=2.2pc @H=2.2pc @R=1pc @C=1pc{\\
*=0{}\ar@{-}[dr] &*=0{} &*=0{}\ar@{-}[dl] &*=0{} &*=0{}\ar@{-}[d] &*=0{} &*=0{}\ar@{-}[d] \\
*=0{} &*=0{}\ar@{-}[d] &*=0{} &*=0{} &*=0{}\ar@{-}[d] &*=0{} &*=0{}\ar@{-}[d] \\
*=0{}&*=0{}&*=0{}&*=0{}&*=0{}&*=0{}&*=0{}&*=0{}&*=0{}\\
*=0{} &*=0{}\ar@{-}[drrr] &*=0{} &*=0{} &*=0{}\ar@{-}[d] &*=0{} &*=0{}\ar@{-}[dll] \\
*=0{} &*=0{} &*=0{} &*=0{} &*=0{}\ar@{-}[d] &*=0{} &*=0{} \\
*=0{}&*=0{}&*=0{}&*=0{}&*=0{}&*=0{}&*=0{}\\
*=0{} \\
*=0{} &*=0{} &*=0{} &*=0{} &*=0{}\ar[d]^{\gamma}&*=0{}\\
*=0{} &*=0{} &*=0{} &*=0{} &*=0{} &*=0{}
}
$$ \\
$$
\xymatrix@W=2.2pc @H=2.2pc @R=1pc @C=1pc{
*=0{}\ar@{-}[dr] &*=0{} &*=0{}\ar@{-}[dl] &*=0{} &*=0{}\ar@{-}[dr] &*=0{} &*=0{}\ar@{-}[dl] \\
*=0{} &*=0{}\ar@{-}[d] &*=0{} &*=0{} &*=0{} &*=0{}\ar@{-}[d] &*=0{} \\
*=0{}&*=0{}&*=0{}&*=0{}&*=0{}&*=0{}&*=0{}&*=0{}&*=0{}&*=0{}\ar[rr]^{\gamma}&*=0{}&*=0{}&*=0{}&*=0{} \\
*=0{} &*=0{}\ar@{-}[drr]  &*=0{} &*=0{} &*=0{} &*=0{}\ar@{-}[dll]\\
*=0{} &*=0{} &*=0{} &*=0{}\ar@{-}[d] \\
*=0{} &*=0{} &*=0{} &*=0{} \\
}
$$
& \hspace{-2pc}
$$
\xymatrix@W=2.2pc @H=2.2pc @R=1pc @C=1pc{\\
*=0{}\ar@{-}[drr] &*=0{}\ar@{-}[dr] &*=0{} &*=0{}\ar@{-}[dl] &*=0{}\ar@{-}[dll]\\
*=0{} &*=0{} &*=0{}\ar@{-}[d] &*=0{} \\
*=0{} &*=0{} &*=0{} &*=0{}
}
$$
\end{tabular}
\end{center}

In the above pictures the tensor products are shown just by
juxtaposition, but now we would like to think about the products
more explicitly. If the monoidal category is not strict, then there
is actually required another leg of the associativity diagram, where
the tensoring is reconfigured  so that the composition can operate
in an alternate order.  Here is how that rearranging looks in a
symmetric (braided) category, where the shuffling is accomplished
by use of the symmetry (braiding):

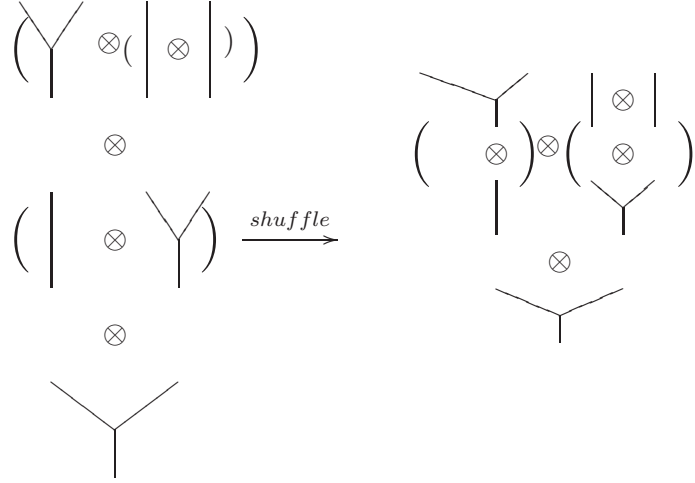
\begin{center}
\begin{tabular}{ll}
$$
\xymatrix@W=2.2pc @H=2.2pc @R=1.5pc @C=1pc{
*=0{}\ar@{-}[dr] &*=0{} &*=0{}\ar@{-}[dl] &*=0{} &*=0{}\ar@{-}[d] &*=0{} &*=0{}\ar@{-}[d] \\
*=0{\txt{\huge(}} &*=0{}\ar@{-}[d] &*=0{} &*=0{\otimes\txt{(}} &*=0{}\ar@{-}[d] &*=0{\otimes} &*=0{}\ar@{-}[d] &*=0{\left.\right)\txt{\huge)}}\\
*=0{}&*=0{}&*=0{}&*=0{}&*=0{}&*=0{}&*=0{}&*=0{}&*=0{}\\
*=0{}&*=0{}&*=0{}&*=0{\otimes}&*=0{}&*=0{}&*=0{}&*=0{}&*=0{}\\
*=0{} &*=0{}\ar@{-}[d] &*=0{} &*=0{} &*=0{}\ar@{-}[dr] &*=0{} &*=0{}\ar@{-}[dl] \\
*=0{\txt{\huge(}} &*=0{}\ar@{-}[d] &*=0{} &*=0{\otimes} &*=0{} &*=0{}\ar@{-}[d] &*=0{\txt{\huge)}}&*=0{}\ar[rrr]^{shuffle}&*=0{}&*=0{}&*=0{}&*=0{} \\
*=0{}&*=0{}&*=0{}&*=0{}&*=0{}&*=0{}&*=0{}&*=0{}&*=0{}\\
*=0{}&*=0{}&*=0{}&*=0{\otimes}&*=0{}&*=0{}&*=0{}&*=0{}&*=0{}\\
*=0{} &*=0{}\ar@{-}[drr]  &*=0{} &*=0{} &*=0{} &*=0{}\ar@{-}[dll]\\
*=0{} &*=0{} &*=0{} &*=0{}\ar@{-}[d] \\
*=0{} &*=0{} &*=0{} &*=0{} \\
}
$$ & \hspace{-1.25pc}
$$
\xymatrix@W=2.2pc @H=2.2pc @R=.85pc @C=1pc{\\
*=0{}\ar@{-}[dr] &*=0{} &*=0{}\ar@{-}[dl] &*=0{} &*=0{}\ar@{-}[d] &*=0{} &*=0{}\ar@{-}[d] \\
*=0{} &*=0{}\ar@{-}[d] &*=0{} &*=0{} &*=0{}\ar@{-}[d] &*=0{\otimes} &*=0{}\ar@{-}[d] &*=0{}\\
*=0{}&*=0{}&*=0{}&*=0{}&*=0{}&*=0{}&*=0{}&*=0{}&*=0{}\\
*=0{\txt{\Huge(}}&*=0{\otimes}&*=0{\txt{\Huge)}}&*=0{\otimes\txt{\Huge(}}&*=0{}&*=0{\otimes}&*=0{}&*=0{\txt{\Huge)}}&*=0{}\\
*=0{} &*=0{}\ar@{-}[d] &*=0{} &*=0{} &*=0{}\ar@{-}[dr] &*=0{} &*=0{}\ar@{-}[dl] \\
*=0{} &*=0{}\ar@{-}[d] &*=0{} &*=0{} &*=0{} &*=0{}\ar@{-}[d] &*=0{}&*=0{} \\
*=0{}&*=0{}&*=0{}&*=0{}&*=0{}&*=0{}&*=0{}&*=0{}&*=0{}\\
*=0{}&*=0{}&*=0{}&*=0{\otimes}&*=0{}&*=0{}&*=0{}&*=0{}&*=0{}\\
*=0{} &*=0{}\ar@{-}[drr]  &*=0{} &*=0{} &*=0{} &*=0{}\ar@{-}[dll]\\
*=0{} &*=0{} &*=0{} &*=0{}\ar@{-}[d] \\
*=0{} &*=0{} &*=0{} &*=0{} \\
}
$$
\end{tabular}
\end{center}

We now foreshadow our definition of operads in an iterated monoidal
category with the same picture as above but using two tensor products,
$\otimes_1$ and $\otimes_2.$ It becomes clear that the true nature
of the shuffle is in fact that of an interchange transformation.

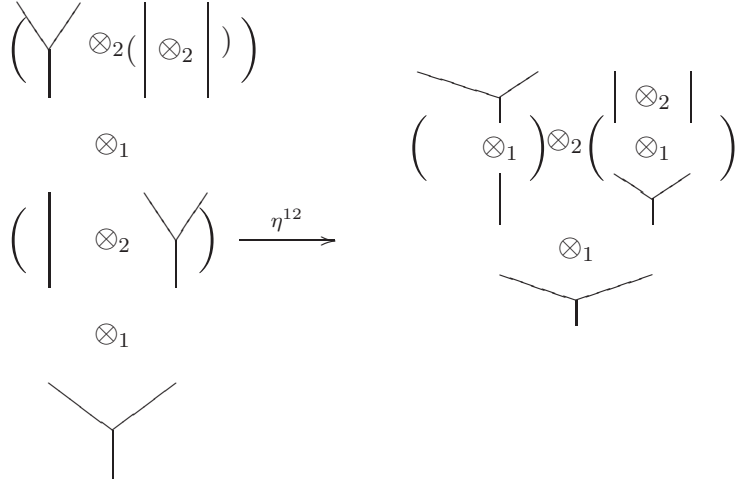
\begin{center}
\begin{tabular}{ll}
$$
\xymatrix@W=2.2pc @H=2.2pc @R=1.5pc @C=1pc{
*=0{}\ar@{-}[dr] &*=0{} &*=0{}\ar@{-}[dl] &*=0{} &*=0{}\ar@{-}[d] &*=0{} &*=0{}\ar@{-}[d] \\
*=0{\txt{\huge(}} &*=0{}\ar@{-}[d] &*=0{} &*=0{\otimes_2\txt{(}} &*=0{}\ar@{-}[d] &*=0{\otimes_2} &*=0{}\ar@{-}[d] &*=0{\left.\right)\txt{\huge)}}\\
*=0{}&*=0{}&*=0{}&*=0{}&*=0{}&*=0{}&*=0{}&*=0{}&*=0{}\\
*=0{}&*=0{}&*=0{}&*=0{\otimes_1}&*=0{}&*=0{}&*=0{}&*=0{}&*=0{}\\
*=0{} &*=0{}\ar@{-}[d] &*=0{} &*=0{} &*=0{}\ar@{-}[dr] &*=0{} &*=0{}\ar@{-}[dl] \\
*=0{\txt{\huge(}} &*=0{}\ar@{-}[d] &*=0{} &*=0{\otimes_2} &*=0{} &*=0{}\ar@{-}[d] &*=0{\txt{\huge)}}&*=0{}\ar[rrr]^{\eta^{12}}&*=0{}&*=0{}&*=0{}&*=0{} \\
*=0{}&*=0{}&*=0{}&*=0{}&*=0{}&*=0{}&*=0{}&*=0{}&*=0{}\\
*=0{}&*=0{}&*=0{}&*=0{\otimes_1}&*=0{}&*=0{}&*=0{}&*=0{}&*=0{}\\
*=0{} &*=0{}\ar@{-}[drr]  &*=0{} &*=0{} &*=0{} &*=0{}\ar@{-}[dll]\\
*=0{} &*=0{} &*=0{} &*=0{}\ar@{-}[d] \\
*=0{} &*=0{} &*=0{} &*=0{} \\
}
$$ & \hspace{-1.25pc}
$$
\xymatrix@W=2.2pc @H=2.2pc @R=.8pc @C=1.2pc{\\
*=0{}\ar@{-}[dr] &*=0{} &*=0{}\ar@{-}[dl] &*=0{} &*=0{}\ar@{-}[d] &*=0{} &*=0{}\ar@{-}[d] \\
*=0{} &*=0{}\ar@{-}[d] &*=0{} &*=0{} &*=0{}\ar@{-}[d] &*=0{\otimes_2} &*=0{}\ar@{-}[d] &*=0{}\\
*=0{}&*=0{}&*=0{}&*=0{}&*=0{}&*=0{}&*=0{}&*=0{}&*=0{}\\
*=0{\txt{\Huge(}}&*=0{\otimes_1}&*=0{\txt{\Huge)}}&*=0{\otimes_2\txt{\Huge(}}&*=0{}&*=0{\otimes_1}&*=0{}&*=0{\txt{\Huge)}}&*=0{}\\
*=0{} &*=0{}\ar@{-}[d] &*=0{} &*=0{} &*=0{}\ar@{-}[dr] &*=0{} &*=0{}\ar@{-}[dl] \\
*=0{} &*=0{}\ar@{-}[d] &*=0{} &*=0{} &*=0{} &*=0{}\ar@{-}[d] &*=0{}&*=0{} \\
*=0{}&*=0{}&*=0{}&*=0{}&*=0{}&*=0{}&*=0{}&*=0{}&*=0{}\\
*=0{}&*=0{}&*=0{}&*=0{\otimes_1}&*=0{}&*=0{}&*=0{}&*=0{}&*=0{}\\
*=0{} &*=0{}\ar@{-}[drr]  &*=0{} &*=0{} &*=0{} &*=0{}\ar@{-}[dll]\\
*=0{} &*=0{} &*=0{} &*=0{}\ar@{-}[d] \\
*=0{} &*=0{} &*=0{} &*=0{} \\
}
$$
\end{tabular}
\end{center}

To see this just focus on the actual domain and range of $\eta^{12}$
which are the upper two levels of trees in the pictures, with the
tensor product $\left({\mathbf|}\otimes_2{\mathbf|}\right)$ considered
as a single object.

Now we are ready to give the technical definitions. We begin with
the definition of 2-fold operad in an $n$-fold monoidal category,
as in the above picture, and then mention how it generalizes the
case of operad in a braided category. Because of this generalization
of the well known case, and since there are easily described examples
of 2-fold monoidal categories based on a braided category as in
\cite{forcey3}, it seems worthwhile to work out the theory for the
2-fold operads in its entirety before moving on to $n$-fold operads.

Let ${\cal V}$ be an $n$-fold monoidal category as defined in Section 2.

\begin{definition}
A 2-fold operad ${\cal C}$ in ${\cal V}$ consists of
objects ${\cal C}(j)$, $j\ge 0$,
a unit map ${\cal J}:I\to {\cal C}(1)$,
and composition maps in ${\cal V}$
\[
\gamma^{12}:{\cal C}(k) \otimes_1 ({\cal C}(j_1) \otimes_2 \dots \otimes_2 {\cal C}(j_k))\to {\cal C}(j)
\]
for $k\ge 1$, $j_s\ge0$ for $s=1\dots k$ and
$\smash{\sum\limits_{n=1}^k j_n = j}$. The composition maps obey
the following axioms:
\begin{enumerate}
\item Associativity: The following diagram is required to commute
for all  $k\ge 1$, $j_s\ge 0$ and $i_t\ge 0$, and
where $\sum\limits_{s=1}^k j_s = j$ and $\sum\limits_{t=1}^j i_t = i.$
Let $g_s= \sum\limits_{u=1}^s j_u$ and
let $h_s=\sum\limits_{u=1+g_{s-1}}^{g_s} i_u$.  The $\eta^{12}$
labelling the leftmost arrow actually stands for a variety of
equivalent maps which factor into instances of the $12$
interchange.
\[
\xymatrix{
{\cal C}(k)\otimes_1\left(\bigotimes\limits_{s=1}^k{}_2 {\cal C}(j_s)\right)\otimes_1
\left(\bigotimes\limits_{t=1}^j{}_2 {\cal C}(i_t)\right)
\ar[rr]^>>>>>>>>>>>>{\gamma^{12} \otimes_1 \text{id}}
\ar[dd]_{\text{id} \otimes_1 \eta^{12}}
&&{\cal C}(j)\otimes_1 \left(\bigotimes\limits_{t=1}^j{}_2 {\cal C}(i_t)\right)
\ar[d]^{\gamma^{12}}
\\
&&{\cal C}(i)
\\
{\cal C}(k)\otimes_1 \left(\bigotimes\limits_{s=1}^k{}_2 {\cal C}(j_s)\otimes_1
\left(\bigotimes\limits_{u=1}^{j_s}{}_2 {\cal C}(i_{u+g_{s-1}})\right)\right)
\ar[rr]_>>>>>>>>>>{\text{id} \otimes_1(\otimes^k_2\gamma^{12})}
&&{\cal C}(k)\otimes_1\left(\bigotimes\limits_{s=1}^k{}_2 {\cal C}(h_s)\right)
\ar[u]_{\gamma^{12}}
}
\]

\item Respect of units is required just as in the symmetric case.
The following unit diagrams commute.
\[
\xymatrix{
{\cal C}(k)\otimes_1 (\otimes_2^k I)
\ar[d]_{1\otimes_1(\otimes_2^k {\cal J})}
\ar@{=}[r]^{}
&{\cal C}(k)\\
{\cal C}(k)\otimes_1(\otimes_2^k {\cal C}(1))
\ar[ur]_{\gamma^{12}}
}
\xymatrix{
I\otimes_1 {\cal C}(k)
\ar[d]_{{\cal J}\otimes_1 1}
\ar@{=}[r]^{}
&{\cal C}(k)\\
{\cal C}(1)\otimes_1 {\cal C}(k)
\ar[ur]^{\gamma^{12}}
}
\]
\end{enumerate}
\end{definition}

Note that operads in a braided monoidal category are examples of
2-fold operads. This is true based on the arguments of Joyal and
Street \cite{JS}, who showed that braided categories arise as 2-fold
monoidal categories where the interchanges are isomorphisms. Also
note that given such a perspective on a braided category, the two
products are equivalent and the use of the braiding to shuffle in
the operad associativity requirement can be rewritten as the use
of the interchange.

It is immediately clear that we can define operads using more than
just the first two products in an $n$-fold monoidal category. The
best way of going about this is to use the theory of monoids, (and
more generally enriched categories), in iterated monoidal categories.
We continue by first describing this procedure for 2-fold operads.
Operads in a symmetric (braided) monoidal category are often
efficiently defined as the monoids of a category of collections.
For a braided category $({\cal V}, \otimes)$ with coproducts that
are preserved by both functors $ (\dash\otimes A)$ and $(A\otimes\dash)$
the objects of $\col({\cal V})$ are functors from the category of
natural numbers to ${\cal V}.$ In other words the data for a
collection ${\cal C}$ is a sequence of objects ${\cal C}(j).$
Morphisms in $\col({\cal V})$ are natural transformations.  The
tensor product in $\col({\cal V})$ is given by
\[
({\cal B}\otimes{\cal C})(j) =
    \coprod\limits_{\substack{k\ge 0\\ j_1+\dots +j_k =j}}
    {\cal B}(k) \otimes({\cal C}(j_1) \otimes \dots \otimes {\cal C}(j_k))
\]
where $j_i \ge 0.$ This product is associative by use of the symmetry
or braiding, and due to the hypothesis that the tensor product
preserves the coproduct.  The unit is the collection
$(\emptyset, I, \emptyset, \dots)$ where $\emptyset$ is an initial
object in ${\cal V}.$

Now recall how the interchange transformations generalize braiding.
For ${\cal V}$ a $2$-fold monoidal category with all coproducts in
which both $\otimes_1$ and $\otimes_2$ preserve the coproduct,
define the objects and morphisms of      $ \col_2(\cal V)$ in
precisely the same way as in the braided case, but define the product
to be
\[
({\cal B}{\otimes}^{12}{\cal C})(j) =
    \coprod\limits_{\substack{k\ge 0\\ j_1+\dots +j_k =j}}
    {\cal B}(k) \otimes_1({\cal C}(j_1) \otimes_2
    \dots \otimes_2 {\cal C}(j_k))
\]
In general the interchangers will not be isomorphisms, so this
product can not be that of a monoidal category with the usual strong
associativity.  However the interchangers can be used to make the
product in question obey lax associativity, where the associator
is a coherent natural transformation.  This lax associativity is
seen by inspection of the two 3-fold products
$({\cal B}{\otimes}^{12}{\cal C}){\otimes}^{12}{\cal D}$
and
${\cal B}{\otimes}^{12}({\cal C}{\otimes}^{12}{\cal D})$.
In the braided case mentioned above, the two large coproducts in
question are seen to be composed of the same terms up to a braiding
between them.  Here the terms of the two coproducts are related by
instances of the interchange transformation $\eta^{12}$ from the
term in $(({\cal B}{\otimes}^{12}{\cal C}){\otimes}^{12}{\cal D})(j)$
to the corresponding term in
$({\cal B}{\otimes}^{12}({\cal C}{\otimes}^{12}{\cal D}))(j).$
For example upon expansion of the two three-fold products we see
that in the coproduct which is
$(({\cal B}{\otimes}^{12}{\cal C}){\otimes}^{12}{\cal D})(2)$
we have the term
\[{\cal B}(2)\otimes_1({\cal C}(1)\otimes_2{\cal C}(1))\otimes_1({\cal D}(1)\otimes_2{\cal D}(1))\]
while in
$({\cal B}{\otimes}^{12}({\cal C}{\otimes}^{12}{\cal D}))(2)$
we have the term
\[{\cal B}(2)\otimes_1({\cal C}(1)\otimes_1{\cal D}(1))\otimes_2({\cal C}(1)\otimes_1{\cal D}(1)).\]
Note that the first of these terms appears courtesy of the fact
that when a tensor product preserves coproducts, there is implied
a distributive law
$(\coprod B_n) \otimes A \cong \coprod(B_n\otimes A)$
as shown in \cite{MacLane}.

The coherence theorem of iterated monoidal categories as it is
stated in \cite{Balt} guarantees the commutativity of the pentagon
equation for the associators, since they are defined as compositions
of interchangers $\eta^{12}$ in ${\cal V}$.  Some remarks about the
non-invertibility of $\alpha$ are in order. Note that  Mac Lane
proves his coherence theorem in two steps \cite{MacLane}. First it
is shown that every diagram involving only $\alpha$ (no $\alpha^{-1}$)
commutes. Then it is noted that this suffices to make every diagram
of both $\alpha$ and $\alpha^{-1}$ commute since for  every binary
word there exists a path of just instances of $\alpha$ from that
word to the word parenthesized all to the right. (Here we are taking
the domain of $\alpha$ to be $(A\otimes B)\otimes C.$) Thus when
$\alpha$ is not invertible we still have that every diagram commutes.
There are still canonical maps from every binary word to the word
parenthesized all to the right.  However there are necessarily fewer
diagrams. For instance if $({\cal V}, \otimes)$ is lax monoidal
there is no canonical map between the two objects $(B\otimes
B)\otimes(B\otimes B)$ and $(B\otimes (B\otimes B))\otimes B$. This
affects the statement of the general associativity theorem for
monoids in a lax monoidal category.  Only the specific case of the
general associativity theorem as stated by Mac Lane  holds, as
follows.

\begin{theorem}
Let $(A,\mu)$ be a monoid in a (lax) monoidal category. Let $A^{n}$
be the product given by $A^2 = A \otimes A, A^{n+1} = A\otimes A^n$,
i.e. parenthesized to the right. Define the composition $\mu^{(n)}$
by $\mu^{(2)} = \mu, \mu^{(n+1)} = \mu \circ (1 \otimes \mu^{(n)})$.
Then
\[ \mu^{(n)} \circ (\mu^{(k_1)}\otimes \dots \otimes \mu^{(k_n)}) = \mu^{(k_1 + \dots +k_n)}\circ \alpha' \]
for all $n, k_i \ge 2$ where $\alpha'$ stands for the canonical map
to $A^{k_1 + \dots +k_n}.$
\end{theorem}

\begin{proof}
This is just the special case of the general associative law for
monoids shown by Mac Lane, which only depends on the existence of
the canonical map $\alpha'$ \cite{MacLane}.
\end{proof}
Now we have a condensed way of defining 2-fold operads.
\begin{theorem}
2-fold operads in $2$-fold monoidal ${\cal V}$ are monoids in
$\col_2(\cal V).$
\end{theorem}
\begin{proof}
A monoid in $ \col_2(\cal V)$ is an object ${\cal C}$ in $\col_2(\cal V)$
with multiplication and unit morphisms. Since morphisms of
$\col_2(\cal V)$ are natural transformations the multiplication and
unit consist of families of maps in ${\cal V}$ indexed by the natural
numbers, with source and target exactly as required for operad
composition and unit. The operad axioms are equivalent to the
associativity and unit requirements of monoids.
\end{proof}

This brings us back to the question of defining operads in $n$-fold
monoidal ${\cal V}$ using the higher products and interchanges.
This idea will correspond to a series of higher products, denoted
by $\otimes^{pq},$ in the category of collections. These are defined
just as for the first case $\otimes^{12}$ above.  Associators are
as described above for the first product, using $\eta^{pq}$ for the
associator
$\alpha\colon{\cal A}\otimes^{pq}({\cal B}\otimes^{pq}{\cal C})
    \to ({\cal A}\otimes^{pq}{\cal B})\otimes^{pq}{\cal C}$.
The unit for each is the collection  $(\emptyset, I, \emptyset, \dots)$
where $\emptyset$ is an initial object in ${\cal V}.$ Notice that
these products do not interchange; i.e they are not functorial with
respect to each other.  Notice also that the associators in these
categories of collections are not isomorphisms unless we are
considering the special cases of braiding or symmetry. Instead they
are lax monoidal, by which we will mean that the associator is
merely a natural transformation which obeys the pentagon coherence
condition.

Now we will focus on the products $\otimes^{(m-1)m}$ in the category
of collections in  $n$-fold monoidal ${\cal V}$, for $m\le n,$ since
these will be seen to suffice for defining all operad compositions.
Before defining $m$-fold operads as monoids with respect to
$\otimes^{(m-1)m},$ we note that there is also fibrewise monoidal
structure.  This will be important in the description of the monoidal
structure of the category of operads. In fact, we have the following

\begin{theorem}\label{foo}
If an $n$-fold monoidal category $ {\cal V}$ has coproducts and
$({\cal V},\coprod,\otimes_3,\dots, \otimes_n)$ is an $(n-1)$-fold
monoidal category for which each of the functors $(\dash\otimes_i A)$
and $(A\otimes_i\dash)$ preserves coproducts, then for $n\ge m \ge 2,$
the category of collections in ${\cal V}$ can be given the structure
of an $(n-m+1)$-fold lax monoidal category, denoted $\col_m(\cal V)$.
\end{theorem}

\begin{proof}
The first tensor product is $\hat{\otimes}_1 = \otimes^{(m-1)m}$
and the others are the higher fibrewise products starting with
fibrewise $\otimes_{m+1}.$ Thus the products of $\col_m(\cal V)$
are as follows:

\begin{align*}
({\cal B}\hat{\otimes}_1{\cal C})(j) &=
    \coprod\limits_{\substack {k\ge 0\\ j_1+\dots +j_k =j}}
    {\cal B}(k) \otimes_{m-1}
    ({\cal C}(j_1) \otimes_m \dots \otimes_m {\cal C}(j_k))
\intertext{and}
({\cal B}\hat{\otimes}_2{\cal C})(j) &= {\cal B}(j){\otimes}_{m+1}{\cal C}(j) \\
&\hspace{0.5em}\vdots\\
({\cal B}\hat{\otimes}_{n-m+1}{\cal C})(j) &= {\cal B}(j){\otimes}_{n}{\cal C}(j)
\end{align*}

The unit for $\hat{\otimes}_1$ is
${\cal I} = (\emptyset, I, \emptyset, \dots)$
and the unit for all the other products is
${\cal M} = (I, I, \dots)$.
First we must check that there are natural transformations
\[\xi^{1j}: (A\hat{\otimes}_jB)\hat{\otimes}_1(C\hat{\otimes}_jD)
    \to (A\hat{\otimes}_1C)\hat{\otimes}_j(B\hat{\otimes}_1D)\]
These utilize the $\eta^{ij}$ of ${\cal V}$ and thus exist by
inspection of the terms of the compound products.  For example, in
\[ ((A\hat{\otimes}_2 B)\hat{\otimes}_1(C\hat{\otimes}_2D))(2) \]
we find the term
\[(A(2)\otimes_3 B(2))\otimes_1((C(1)\otimes_3 D(1))
    \otimes_2(C(1)\otimes_3 D(1))), \]
while in
\[ (A\hat{\otimes}_1C)\hat{\otimes}_2(B\hat{\otimes}_1D) \]
we find the two terms
\[ (A(2)\otimes_1 (C(1)\otimes_2 C()1))\quad\text{and}\quad
    (B(2))\otimes_1((D(1)\otimes_2 D(1))) \]
in two separate coproducts which are joined by $\otimes_3.$

The map $\xi^{12}$ thus uses first $\eta^{23}$, then $\eta^{13}$
and finally the hypothesis that
$( {\cal V},\coprod,\otimes_3,\dots, \otimes_n)$
is an $(n-1)$-fold monoidal category; specifically instances of the
map $(X \otimes_3 Y)\coprod(Z\otimes_3 W)$
to $(X \coprod Z)\otimes_3(Y\coprod W)$.

It is also not hard to check the unit conditions which are required
for the fibrewise products to be the multiplication for a monoid
in the category of monoidal categories.  The extra requirement of
the two sorts of unit is that ${\cal M} \hat{\otimes}_1 {\cal M} =
{\cal M} $ and that ${\cal I} \hat{\otimes}_k {\cal I} = {\cal I}$
for $k >1.$ These equations do indeed hold. Thus the first product
together with any of the fibrewise products are those of a 2-fold
monoidal category.

For the products $\hat{\otimes}_{2}$ and higher the associators and
interchange transformations are fibrewise and the axioms hold since
they hold for each fiber.  Lastly we need to mention that the giant
hexagon diagram for $i,j,k= 1,j,k$ commutes.  This can be seen by
splitting the hexagon into two commuting diagrams, one made up of
the fibrewise applied interchangers of  ${\cal V}$ and another that
is a giant hexagon in the $(n-1)$-fold monoidal category $( {\cal
V},\coprod,\otimes_{m+1},\dots, \otimes_n).$
\end{proof}

\begin{remark}
In the context of \cite{batainf} the  lax functoriality of the
tensor product with respect to the coproduct is due to the hypothesis
that symmetric ${\cal V}$ is closed (from the right) with respect
to the tensor product.  This guarantees that that product preserves
colimits on the first operand, since the functor $(\dash\otimes B)$
has as a right adjoint the internal hom, denoted by  $[B,\dash].$
Applied to the coproduct this fact in turn implies that there is a
canonical map in ${\cal V}$ from $(A \otimes B)\coprod(C\otimes D)$
to $(A \coprod C)\otimes (B\coprod D).$ From the universal
properties of the coproduct it can be  checked that this map satisfies
the the middle interchange law that is required of a monoidal
functor.  Also in \cite{batainf} Batanin points out that a fibrewise
product is a monoidal functor with respect to the collection product.
In that paper the existence of the transformation $\xi$ depends on
the symmetry (braiding) and the lax functoriality of the tensor
product with respect to the coproduct. In this paper we chose to
simply include the necessary iterated monoidal structure as a
hypothesis, rather than the hypothesis of closedness, in the interest
of generality.

Theorem~\ref{foo} is quite useful for describing $n$-fold operads
and their higher-categorical structure, especially when coupled
with two other facts. The first is that monoids are equivalently
defined as single object enriched categories, and the second is the
following result from \cite{forcey1} and \cite{forcey2}.  In those
sources the quantifier lax is sometimes left off, but the proofs
in question nowhere require the associator to be an isomorphism.
\end{remark}

\begin{theorem}\label{enrich}
For ${\cal V}$ $n$-fold (lax) monoidal the category of
enriched categories over $({\cal V}, \otimes_1)$
is an $(n-1)$-fold monoidal 2-category.
\end{theorem}

For our purposes we translate the theorem about enriched categories
into its single object corollary about the category $\mon(\cal V)$
of monoids in ${\cal V}$.

\begin{corollary}\label{mon}
For ${\cal V}$ $n$-fold (lax) monoidal, the category $\mon(\cal V)$
is an $(n-1)$-fold monoidal 2-category.
\end{corollary}
\begin{proof}
The product of enriched categories always has as its object set the
cartesian product of the object sets of its components. Thus one
object enriched categories have products with one object as well.
\end{proof}

\begin{definition}\label{itop}
If an $n$-fold monoidal category $ {\cal V}$ has coproducts and
$({\cal V},\coprod,\otimes_3,\dots, \otimes_n)$ is an $(n-1)$-fold
monoidal category in which each of the functors $(\dash\otimes_i A)$
and $(A\otimes_i\dash)$ preserves coproducts we define the
category of $m$-fold operads $\oper_m(\cal V)$ to be the category
of monoids in the category of collections
$( \col_m(\cal V), \hat{\otimes}_1)$ for $n\ge m \ge 2.$
\end{definition}

\begin{corollary}\label{op}
$\oper_m(\cal V)$ is an $(n-m)$-fold monoidal 2-category.
\end{corollary}

\begin{proof}
Rather than starting with monoids in $n$-fold monoidal ${\cal V}$
as in the previous corollary we are actually beginning with monoids
in $(n-m+1)$-fold monoidal $\col_m(\cal V)$.  Note that in
\cite{forcey1} the products in ${\cal V}$ are assumed to have a
common unit. To generalize to our situation here, where the unit
for the first product in the category of collections is distinct,
we need to add slightly to the definitions in \cite{forcey1}. When
enriching (or more specifically taking monoids) we are doing so
with respect to the first available product. Thus the unit morphism
for enriched categories has its domain the unit for that first
product, ${\cal I}.$ However the unit enriched category ${\bcal I}$
has one object, denoted $0$, and ${\bcal I}(0,0) = {\cal M}$.
\end{proof}

\begin{remark}
This theorem justifies our focus on the first $m$ products of ${\cal
V}$ as opposed to any subset of the $n$ products. It is due to the
way in which this focus allows us to describe the resulting structure
on the category of $m$-fold operads. Of course, we can use the
forgetful functors mentioned in Section 2 to pass from $n$-fold
monoidal ${\cal V}$ to ${\cal V}$  with any of the subsets of
products.  The $m$-fold operads do behave as expected under this
forgetting, retaining all but the structure that depends on the
forgotten products. This will be seen more clearly upon inspection
of the unpacked definition to follow.  In short, we will see that
an $m$-fold operad is also an $(m-1)$-fold operad.
\end{remark}

\begin{remark}
We note that since a symmetric monoidal category is $n$-fold monoidal
for all $n$, then operads in a symmetric monoidal category are
$n$-fold monoidal for all $n$ as well.  More generally, if $n\ge
3$ and the interchanges are isomorphisms, then by the Eckmann-Hilton
argument the products collapse into one and the result is a a
symmetric monoidal category, and so operads in it are again $n$-fold
monoidal for all $n.$ Here we are always discussing ordinary
``non-symmetric,'' (``non-braided'') operads. The possible faithful
actions of symmetry or braid groups can be considered after the
definition, which we leave for a later paper. We do point out that
the proper direction in which to expand this work is seen in Weber's
paper \cite{web}. He generalizes by making a distinction between
the binary and $k$-ary products in the domain of the composition
map
$\gamma\colon{\cal C}(k)\otimes({\cal C}(j_1)\otimes
    \dots \otimes{\cal C}(j_k))
    \to {\cal C}(j)$.
The binary tensor product is seen formally as a pseudo-monoid
structure and the $k$-ary product as a pseudo-algebra structure for
a 2-monad which can contain the information needed to describe
actions of braid or symmetry groups. The two structures are defined
using strong monoidal morphisms, and so the products coincide and
give rise to the braiding which is used to describe the associativity
of composition. To encompass the definitions in this paper we would
move to operads in lax-monoidal pseudo algebras, where instead of
pseudo monoids and strong monoidal morphisms in a pseudo algebra
we would consider the same picture but with lax monoidal morphisms.
\end{remark}

The fact that monoids are single object enriched categories also
leads to an efficient expanded definition of $m$-fold operads in
an $n$-fold monoidal category.  Let ${\cal V}$ be an $n$-fold
monoidal category as defined in Section 2.

\begin{definition}\label{expandop}
For $2\le m\le n$ an $m$-fold operad ${\cal C}$ in ${\cal V}$ consists of
objects ${\cal C}(j)$, $j\ge 0$,
a unit map ${\cal J}:I\to {\cal C}(1)$,
and composition maps in ${\cal V}$
\[
\gamma^{pq}:{\cal C}(k) \otimes_p ({\cal C}(j_1) \otimes_q \dots \otimes_q {\cal C}(j_k))\to {\cal C}(j)
\]
for $m\ge q>p \ge 1$, $k\ge 1$, $j_s\ge0$ for $s=1\dots k$ and $\sum\limits_{n=1}^k j_n = j$. The composition maps obey the following axioms:
\begin{enumerate}
\item Associativity: The following diagram is required to commute for all $m\ge q>p \ge 1$, $k\ge 1$, $j_s\ge 0$ and $i_t\ge 0$, and
where $\sum\limits_{s=1}^k j_s = j$ and $\sum\limits_{t=1}^j i_t = i.$ Let $g_s= \sum\limits_{u=1}^s j_u$ and
let $h_s=\sum\limits_{u=1+g_{s-1}}^{g_s} i_u$.

The $\eta^{pq}$ labelling the leftmost arrow actually stands for a variety of equivalent maps which factor into instances of the $pq$ interchange.

\[
\xymatrix{
{\cal C}(k)\otimes_p\left(\bigotimes\limits_{s=1}^k{}_q {\cal C}(j_s)\right)\otimes_p
\left(\bigotimes\limits_{t=1}^j{}_q {\cal C}(i_t)\right)
\ar[rr]^>>>>>>>>>>>>{\gamma^{pq} \otimes_p \text{id}}
\ar[dd]_{\text{id} \otimes_p \eta^{pq}}
&&{\cal C}(j)\otimes_p \left(\bigotimes\limits_{t=1}^j{}_q {\cal C}(i_t)\right)
\ar[d]^{\gamma^{pq}}
\\
&&{\cal C}(i)
\\
{\cal C}(k)\otimes_p \left(\bigotimes\limits_{s=1}^k{}_q {\cal C}(j_s)\otimes_p
\left(\bigotimes\limits_{u=1}^{j_s}{}_q {\cal C}(i_{u+g_{s-1}})\right)\right)
\ar[rr]_>>>>>>>>>>{\text{id} \otimes_p(\otimes^k_q\gamma^{pq})}
&&{\cal C}(k)\otimes_p\left(\bigotimes\limits_{s=1}^k{}_q {\cal C}(h_s)\right)
\ar[u]_{\gamma^{pq}}
}
\]

\item Respect of units is required just as in the symmetric case.
The following unit diagrams commute for all $m\ge q>p \ge 1$.

\[
\xymatrix{
{\cal C}(k)\otimes_p (\otimes_q^k I)
\ar[d]_{1\otimes_p(\otimes_q^k {\cal J})}
\ar@{=}[r]^{}
&{\cal C}(k)\\
{\cal C}(k)\otimes_p(\otimes_q^k {\cal C}(1))
\ar[ur]_{\gamma^{pq}}
}
\xymatrix{
I\otimes_p {\cal C}(k)
\ar[d]_{{\cal J}\otimes_p 1}
\ar@{=}[r]^{}
&{\cal C}(k)\\
{\cal C}(1)\otimes_p {\cal C}(k)
\ar[ur]^{\gamma^{pq}}
}
\]

\end{enumerate}
\end{definition}

\begin{theorem}
The description of $m$-fold operad in Definition~\ref{expandop} is
equivalent to that given in Definition~\ref{itop}.
\end{theorem}

\begin{proof}
If a collection has an operad composition $\gamma^{q,q+1}$ using
$\otimes_q$ and $\otimes_{q+1}$ then it automatically has an operad
composition for any pair of products $\otimes_p$ and $\otimes_{s}$
for $p< s \le q+1$. This  follows from the fact that for $p < q$
we have natural transformations
$ \eta^{pq}_{AIIB}:A\otimes_p B\to A\otimes_q B$,
as described at the end of Definition~\ref{iterated}.
Thus if we have $\gamma^{q,q+1}$ then we can form
$\gamma^{ps} = \gamma^{q,q+1} \circ (\eta^{pq}
    \circ (1 \otimes_q\eta^{s,q+1})).$
The new $\gamma^{ps} $ is associative based on the old $\gamma$'s
associativity, the naturality of $\eta,$ and the coherence of $\eta.$
Thus follows our claim that generally operads are preserved as such
by the forgetful functors mentioned in Section 2 and specifically
that an $m$-fold operad is also an $(m-1)$-fold operad. The converse
of this latter statement is not true, as we will see by counterexample
in the last section. It will demonstrate the existence of $m$-fold
operads which are not $(m+1)$-fold operads.
\end{proof}

It is also worth while to expand the definition of the tensor
products of $m$-fold operads that is implicit in their depiction
as monoids in the category of collections in an $n$-fold monoidal
category. Here is the expanded version of the definition:

\begin{definition}\label{tensor}
Let ${\cal C}$,${\cal D}$ be  $m$-fold operads.  For $1 \le i \le (n-m)$
and using a $\otimes'_k$ to denote the product of two $m$-fold
operads, we define that product to be given by:
\[({\cal C}\otimes'_i {\cal D})(j) = {\cal C}(j) \otimes_{i+m} {\cal D}(j).\]
\end{definition}

We note that the new $\gamma$ is in terms of the two old ones, for
$m\ge q>p \ge 1$:
\[
\gamma^{pq}_{{\cal C}\otimes'_i {\cal D}} =
(\gamma^{pq}_{{\cal C}}\otimes_{i+m}\gamma^{pq}_{{\cal D}})\circ \eta^{p,i+m} \circ (1 \otimes_p \eta^{q,i+m})
\]
where the subscripts denote the $n$-fold operad the $\gamma$ belongs
to and the $\eta$'s actually stand for any of the equivalent maps
which factor into them. Note that this expansion also helps make
clear why it is that the monoidalness, or number of products, of
$m$-fold operads must decrease by the same number $m.$ From the
condensed version this is expected due to the iterated enrichment.
From the expanded view this allows us to define the new composition
since in order for the products of operads to be closed, $\gamma$
for the $i^{th}$ product  utilizes an interchange with superscript
$i+m.$ Defined this way $i$ can only be allowed to be as large as
$n-m.$ We demonstrate in the last section in fact a counterexample
which shows that the degree of monoidalness for the category of
$m$-fold operads in an $n$-fold  monoidal category is in general
no greater than $n-m.$

\section{Examples of iterated monoidal categories}
\label{section:categoryexamples}
\begin{lemma}\label{one}
Given a totally ordered set $S$ with a least element $e \in S$,
then the elements of $S$ make up the objects of a strict monoidal
category.
\end{lemma}
The category will also be denoted $S.$  Its morphisms are given by
the ordering; there is only an arrow $a \to b$ if $a\le b.$ The
product is max and the 2-sided unit is the least element $e$. We
must check that the product is functorial since this defines monoidal
structure on morphisms.  Here it is so since if $a\le b$ and $a'\le
b'$ then $\max(a,a')\le\max(b,b').$ Also the identity is clearly
preserved.

\begin{example}\label{uno}
The basic example is the nonnegative integers $\nat$ with their
ordering~$\le.$
\end{example}

\begin{lemma}\label{semi}
Any ordered monoid with its identity element $e$ also its least
element forms the object set of a 2-fold monoidal category.
\end{lemma}

\begin{proof} Morphisms are again given by the ordering. The products
are given by $\max$ and the monoid operation: $a\otimes_1b = \max(a,b)$
and $a\otimes_2b = ab$.  The shared two-sided unit for
these products is the identity element $0.$ The products are both
strictly associative and functorial since if $a\le b$ and
$a'\le b'$ then $aa'\le bb'$ and $\max(a,a')\le\max(b,b').$ The
interchange natural transformations  exist since
$\max(ab,cd)\le\max(a,c)\max(b,d).$ That is because
\begin{align*}
a\le\max(a,c)\quad&\text{and}\quad b\le\max(b,d) \\
\intertext{so}
ab\le\max(a,c)\max(b,d)\quad&\text{and}\quad cd\le\max(a,c)\max(b,d)
\end{align*}

The internal and external unit and associativity conditions of
Definition~\ref{iterated} are all satisfied due to the fact that
there is only one morphism between two objects. More generally,
given any ordered $n$-fold monoidal category with $I$ the least
object we can potentially form an $(n+1)$-fold monoidal category
with  morphisms ordering, and the new $\otimes_1 = \max.$
\end{proof}

\begin{example}\label{nat}
Again we have in mind $\nat$ with its ordering and addition.
\end{example}
Other examples of such  monoids as in Lemma~\ref{semi} are the pure
braids on $n$ strands with only right-handed crossings \cite{rolf}.
Notice that braid composition is a non-symmetric example.  Further
examples are found in the papers on semirings and idempotent
mathematics, such as \cite{LitSob} and its references as well as
on the related concept of tropical geometry, such as \cite{Sturm}
and its references. Semirings that arise in these two areas of study
usually require some  translation of the lemmas we have stated thus
far, since the idempotent operation is usually $\min$ and its unit
$\infty.$ Also, since the operation given by addition has unit 0,
we have to broaden our definition of 2-fold monoidal category.
Working from the principle that the second operation is the
multiplication  of a categorical monoid with respect to the first,
the additional requirement is that the two distinct units obey each
other's operations: i.e $I_1\otimes_2I_1 = I_1$ and $I_2\otimes_1I_2 = I_2$.
For example, $\min(0,0)= 0$ and $\infty+\infty=\infty.$

\begin{example}\label{seq}
If $S$ is an ordered set then by $\seq(S)$ we denote the infinite
sequences $X_n$ of elements of $S$ for which there exists a natural
number  $l(X)$ called the length such that $k> l(X)$ implies $X_k
= e$ and $X_{l(X)}\ne e.$ Under lexicographic ordering $\seq(S)$
is in turn a totally ordered set with a least element. The latter
is the sequence 0 where $0_n = e$ for all $n.$ We let $l(0) = 0.$
The lexicographic order means that $A \le B$ if either $A_k = B_k$
for all $k$ or there is a natural number $n=n_{AB}$ such that
$A_k = B_k$ for all $k < n$, and such that $A_{n} < B_{n}.$

The ordering is easily shown to be reflexive, transitive, and
antisymmetric. See for instance \cite{Schrod} where the case of
lexicographic ordering of $n$-tuples of natural numbers is discussed.
In our case we will need to modify the proof given in that source
by always making comparisons of $\max(l(A),l(B))$-tuples.

As a category $\seq(S)$ is 2-fold monoidal since we can demonstrate
two interchanging products. They are max using the lexicographic
order: $A\otimes_1 B = \max(A,B);$ and concatenation of sequences:
\[ (A\otimes_2 B)_n =
\begin{cases}
A_n,    &n \le l(A) \\
B_n,    &n >   l(A)
\end{cases} \]
Concatenation clearly preserves the ordering.
\end{example}

\begin{example}
Letting $S$ be the set with a single element recovers Example~\ref{nat}
as $\seq(S)$.
\end{example}

\begin{lemma}
If we have an ordered monoid $(M,+)$ as in Lemma~\ref{semi} and
reconsider $\seq(M)$ as in Example~\ref{seq} then we can describe
a 3-fold monoidal category $\seq(M,+)$ (with $\seq(M)$ the image
of forgetting the third product of pointwise addition) iff the
monoid operation $+$ is such that $0<a<b$ and $c\le d$ imply both
$a+c < b+d$ and $c+a<d+b$ strictly.
\end{lemma}
\begin{proof}
The first two products are again lexicographic max and concatenation
of sequences. The third product $\otimes_3$ is pointwise application
of $+$, $(A\otimes_3B)_n = A_n+B_n$. The last condition that the
monoid operation $+$ strictly respect strict ordering is necessary
to guarantee that the third product both respect the lexicographic
ordering; and interchange correctly  with concatenation. To see the
former let sequences $A\le B, C\le D.$ Note that if  $A=B, C=D$
then $A \otimes_3C = B\otimes_3D$ and if instead (without loss of
generality) $A_j < B_j$ for $j$ such that
$A_i=B_i \text{ and } C_i = D_i \text{ for } i\le j$,
then $A\otimes_3C < B\otimes_3D,$ since $C_j\le D_j.$ To see the
converse, consider a case where $0<a<b$ and $c\le d$ but $a+c=b+d$.
Then the sequences $A=(a,a)$, $B=(b,0)$, $C=(c,0)$, $D=(d,0)$ are
such that lexicographically $A<B$ and $C\le D$ but
$A\otimes_3C = (a+c,a) > B\otimes_3D = (b+d,0).$
To see the interchange
$(A\otimes_3B)\otimes_2(C\otimes_3D) \le (A\otimes_2C)\otimes_3(B\otimes_2D) $
notice that we can assume that $l(A) > l(B).$ Then

\[ \concat(A+B,C+D)\le \concat(A,C)+ \concat(B,D) \]

due to the fact that if $D$ has a first non-zero term, it will be
added to an earlier term of the concatenation of $A$ and $C$ in the
second four-fold product.
\end{proof}

\begin{remark}
A non-example is seen if we begin with the monoid of Lemma~\ref{one},
an ordered set with a least element where the product is max. Here
max does not strictly preserve strict ordering, and so pointwise
max does not respect lexicographic ordering. Neither do concatenation
and pointwise max interchange.
\end{remark}

\begin{corollary}\label{duh}
Given any ordered $n$-fold monoidal category $C$ with $I$ the least
object and $\otimes_1$ the max, and whose higher products strictly
respect strict ordering, we can form an $(n+1)$-fold monoidal
category $\seq(C)$.
\end{corollary}
\begin{proof}
The new products of $\seq(C)$ are the lexicographic max, the
concatenation, and the pointwise application of $\otimes_i$ for $i
= 2\dots n.$ The pointwise application of the original products to
the sequences directly inherits the interchange properties. For
instance, if $A, B, C, D \in $ $\seq(C)$ then
$(A_n\otimes_2B_n)\otimes_1(C_n\otimes_2D_n) \le
(A_n\otimes_1C_n)\otimes_2(B_n\otimes_1D_n)$ for all $n$, which
certainly implies that the pointwise 4-fold products are ordered
lexicographically.
\end{proof}

\begin{example}\label{2castles}
Even more symmetrical structure is found in examples with a natural
geometric representation which allows the use of addition in each
product. One such category is that whose objects are Young diagrams,
by which we mean the underlying shapes or diagrams of Young tableaux.
These can be presented by a decreasing sequence of nonnegative
integers in two ways: the sequence that gives the heights of the
columns or the sequence that gives the lengths of the rows. We let
$\otimes_2$ be the product which adds the heights of columns of two
diagrams, $\otimes_1$ adds the length of rows. We often refer to
these as vertical and horizontal stacking respectively.  If
\[
A = \xymatrix@W=.75pc @H=.75pc @R=0pc @C=0pc @*[F-]{~&~&~&~\\~}
\text{ and } B = \xymatrix@W=.75pc @H=.75pc @R=0pc @C=0pc @*[F-]{~&~\\~\\~}
\]
\begin{align*}
\text{then } A\otimes_1 B &= \xymatrix@W=.75pc @H=.75pc @R=0pc @C=0pc @*[F-]{~&~&~&~&~&~\\~&~\\~} \\
\text{ and } A\otimes_2 B &= \xymatrix@W=.75pc @H=.75pc @R=0pc @C=0pc @*[F-]{~&~&~&~\\~&~\\~\\~\\~}
\end{align*}

There are several possibilities for morphisms.  We can take as
morphisms the totally ordered structure of the Young diagrams given
by lexicographic ordering. In interest of focusing on the stacking
products though we may choose to restrict these morphisms further,
and say an arrow given by ordering can only exist between similar
mass objects, i.e. the two objects in question have equal sums of
their respective sequences or, in reference to the pictures, an
equal total number of blocks.  This restriction eliminates the
product described by lexicographic max.  By the category of restricted
Young diagrams, we will refer to morphisms as restricted lexicographic
ordering, and the two stacking products demonstrated above. We will
often find occasion to relax the morphisms to include all ordering
and reintroduce the lexicographic max as $\otimes_1$, and will refer
to that category simply as the category of Young diagrams.

By previous discussion of sequences the Young diagrams with $\otimes_1$
the lexicographic max and  $\otimes_3$ the piecewise addition
(thought of here as vertical stacking) form a subcategory of the
3-fold monoidal category called $\seq(\nat,+)$.  To see that with
the additional $\otimes_2$ of horizontal stacking that this becomes
a valid 3-fold monoidal category  we look at that operation from
another point of view. Note that the horizontal product of Young
diagrams $A$ and $C$ can be described as a reorganization of all
the columns of both $A$ and $C$ into a new Young diagram made up
of those columns in descending order of height. Rather than (but
equivalent to) the addition of rows, we see horizontal stacking as
the concatenation of monotone decreasing sequences (of columns)
followed by sorting greatest to least.  We call this operation
merging.
\end{example}

\begin{lemma}\label{sort}
Let $(S, \le, +)$ be an ordered monoid and consider the sequences
$\seq(S,+)$ with lexicographic ordering, piecewise addition $+$ and
the function of sorting denoted by
\[ s\colon\seq(S,+)\to \seq(S,+) \]
Then the triangle inequality holds for two sequences:
$s(A+B)\le s(A)+s(B).$
\end{lemma}

\begin{proof}
Consider $s(A+B)$, where we start with the two sequences and add
them piecewise before sorting. We can metamorphose this into
$s(A)+s(B)$ in stages by using an algorithm to sort $A$ and $B$.
Note that if $A$ and $B$ are already sorted, the inequality becomes
an equality.  For our algorithm we choose parallel bubble sorting.
This consists of a series of passes through the sequences comparing
$A_n$ and $A_{n+1}$ and  comparing $B_n$ and $B_{n+1}$ simultaneously.
If the two elements of a given sequence are not already in strictly
decreasing order we switch their places.  We claim that switching
consecutive sequence elements into order always results in a
lexicographically larger sequence after adding piecewise and sorting.
If both the elements of $A$ and of $B$ are switched, or if neither,
then the result is unaltered. Therefore without loss of generality
we assume that $A_n < A_{n+1}$ and that $B_{n+1} < B_n.$ Then we
compare the original result of sorting after adding and the same
but after the switching of $A_n$ and $A_{n+1}.$ It is simplest to
note that the new result includes $A_{n+1} + B_n,$ which is larger
than both $A_n + B_n$ and $A_{n+1} + B_{n+1}.$ So after adding and
sorting the new result is indeed larger lexicographically.  Thus
since each move of the parallel bubble sort results in a larger
expression after first adding and then sorting, and after all the
moves the result of adding and then sorting the now presorted
sequences is the same as first sorting then adding, the triangle
inequality follows.
\end{proof}

\begin{theorem}
The category of restricted Young diagrams forms a 2-fold monoidal
category, and the category of Young diagrams forms a 3-fold monoidal
category.
\end{theorem}
\begin{proof}
We show the latter statement is true, and then note that the the
former statement follows since the restricted Young diagrams are
just the image of forgetting the first product on Young diagrams
and then passing to a subcategory by restricting morphisms.  The
products on Young diagrams are $\otimes_1 = $ lexicographic max,
$\otimes_2 = $ horizontal stacking and $\otimes_3 =$ vertical
stacking. We need to check first  that horizontal stacking, or
merging, is functorial with respect to morphisms (defined as the
$\le$ relations of the lexicographic ordering.)  The cases where
$A=B$ or $C=D$ are easy. For example let $A_k = B_k$ for all $k$
and $C_k = D_k$ for all $k < n_{CD},$ where $n_{CD}$ is as defined
in Example~\ref{seq}. Thus the columns from the copies of, for
instance $A$ in $A\otimes_1 C$ and $A\otimes_1 D$ fall into the
same final spot under the sortings right up to the critical location,
so if $C \le D$, then $A\otimes_1 C \le A\otimes_1 D.$ Similarly,
it is clear that $A \le B$ implies $(A \otimes_1 D) \le (B \otimes_1 D).$
Hence if $A\le B$ and $C\le D$, then
$A\otimes_1 C \le A\otimes_1 D \le B\otimes_1 D$
which by transitivity gives us our desired property.

Next we check that our interchange transformations will always
exist. $\eta^{1j}$ exists by the proof of Lemma~\ref{semi} for
$j=2,3$ since the higher products both respect morphisms(ordering)
and are thus ordered monoid operations. We need to check for existence
of $\eta^{23},$ i.e. we need to show that
$(A \otimes_3 B)\otimes_2(C \otimes_3 D)
    \le (A \otimes_2 C)\otimes_3 (B \otimes_2 D).$
This inequality follows from Lemma~\ref{sort} on the triangle
inequality for sorting.  To prove the new inequality we consider
the special case of two sequences formed by letting $A'$ be $A$
followed by $C$ and letting $B'$ be $B$ followed by $D$. By ``followed
by'' we mean padded by zeroes so that $l(A') = \max(l(A),l(B)) + l(C) $
and $l(B') = \max(l(A),l(B)) + l(D).$ Thus piecewise addition
of $A'$ and $B'$ results in piecewise addition of $A$ and $B$, and
respectively $C$ and $D.$ Then to our new sequences $A'$ and $B'$
we apply the result of Lemma~\ref{sort} and have our desired result.
\end{proof}

Here is an example of the inequality we have just shown to always hold.
Let four Young diagrams be as follow:
\[
A = \xymatrix@W=.75pc @H=.75pc @R=0pc @C=0pc @*[F-]{~&~\\~\\~}~~
B = \xymatrix@W=.75pc @H=.75pc @R=0pc @C=0pc @*[F-]{~&~&~\\~&~}~~
C = \xymatrix@W=.75pc @H=.75pc @R=0pc @C=0pc @*[F-]{~\\~}~~
D = \xymatrix@W=.75pc @H=.75pc @R=0pc @C=0pc @*[F-]{~&~}
\]
Then the fact that
$(A \otimes_3 B)\otimes_2(C \otimes_3 D)
    \le (A \otimes_2 C)\otimes_3 (B \otimes_2 D)$
appears as follows:
\[
\xymatrix@W=.75pc @H=.75pc @R=0pc @C=0pc @*[F-]{~&~&~&~&~\\~&~&~\\~&~&~\\~\\~} \le \xymatrix@W=.75pc @H=.75pc @R=0pc @C=0pc @*[F-]{~&~&~&~&~\\~&~&~\\~&~\\~&~\\~}
\]

\begin{remark}\label{hpre}
Alternatively we can create a category equivalent to the non-negative
integers in Example~\ref{uno} by pre-ordering the Young diagrams
by height.  Here the height $h(A)$ of the Young diagram is the
number of boxes in its leftmost column, and  we say $A \le B$ if
$h(A)\le h(B)$.  Two Young diagrams with the same height are
isomorphic objects, and the one-column stacks form both a full
subcategory and a  skeleton of the height preordered category.
Everything works as for the previous example of natural numbers
since $h(A\otimes_2 B) = h(A) +h(B)$ and $h(A\otimes_1 B) = \max(h(A),
h(B)).$ There is also a max product; the new max with respect to
the height preordering is defined as
\[ \max(A,B) = \begin{cases}
    A,& \text{if $B\le A$} \\
    B,& \text{otherwise.}
\end{cases} \]
In the height preordered category this latter product is equivalent
to horizontal stacking,~$\otimes_1$.
\end{remark}

\begin{remark}
Notice that we can start with any totally ordered monoids $\{M,\le,
+\}$ such that the identity 0 is less than any other element and
such that $0<a<b$ and $c\le d$ implies both $a+c < b+d$ and $c+a<d+b$
for all $a,b,c \in G.$ We create a 3-fold monoidal category
$\modseq(M,+)$ with objects monotone decreasing finitely non-zero
sequences of elements of $M$ and morphisms given by the lexicographic
ordering.  The products are as described for the category of Young
diagrams $\modseq(\nat,+)$ in the previous example.  The common
unit is the zero sequence.   The proofs we have given in the previous
example for $M=\nat$  are all still valid.
\end{remark}

By Corollary~\ref{duh} we can also consider 4-fold monoidal categories
such as $\seq(\modseq(M))$ and other combinations of Seq and ModSeq.
For instance if $\modseq(\nat,+)$ is our category of Young diagrams
then $\modseq(\modseq(\nat,+))$ has objects monotone decreasing
sequences of Young diagrams, which we can visualize along the
$z$-axis. Here the lexicographic-lexicographic max is  $\otimes_1$,
lexicographic merging is $\otimes_2$, pointwise merging (pointwise
horizontal or $y$-axis stacking) is $\otimes_3$ and pointwise-pointwise
addition (pointwise $x$-axis stacking) is $\otimes_4$. For example,
if:
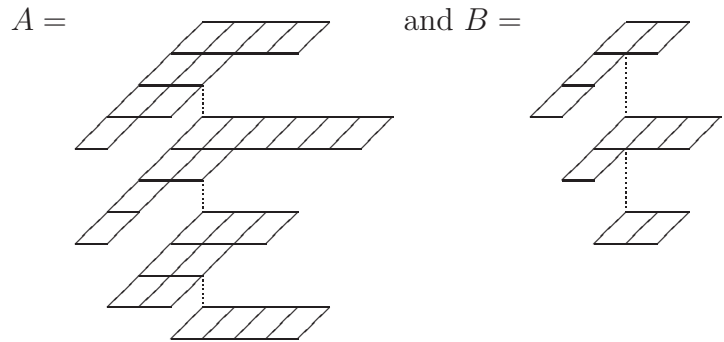
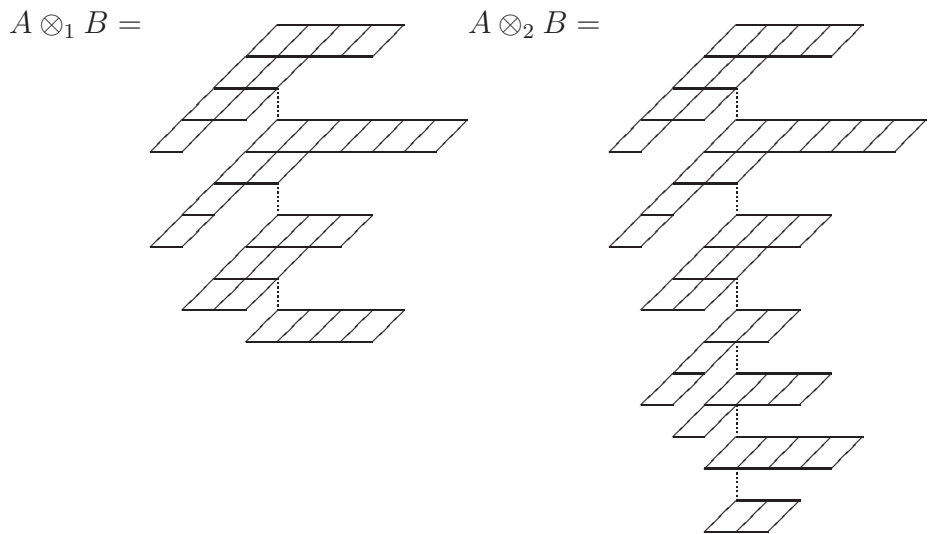
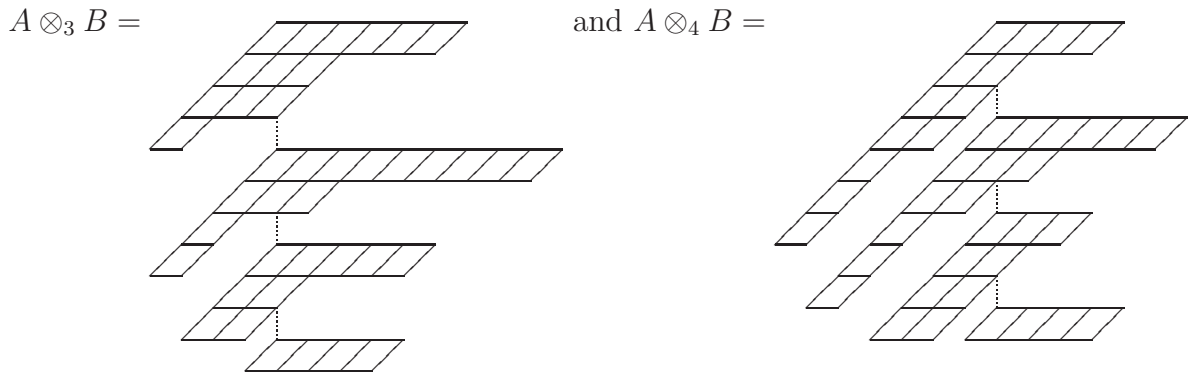
\[
A=\xymatrix @W=1.0pc @H=1.0pc @R=1.0pc @C=1.0pc
{*=0{}&*=0{}&*=0{}&*=0{}&*=0{}\ar@{-}[r]\ar@{-}[dl]&*=0{}\ar@{-}[r]\ar@{-}[dl]&*=0{}\ar@{-}[r]\ar@{-}[dl]&*=0{}\ar@{-}[r]\ar@{-}[dl]&*=0{}\ar@{-}[dl]
\\*=0{}&*=0{}&*=0{}&*=0{}\ar@{-}[r]\ar@{-}[dl]&*=0{}\ar@{-}[r]\ar@{-}[dl]&*=0{}\ar@{-}[r]\ar@{-}[dl]&*=0{}\ar@{-}[r]&*=0{}
\\*=0{}&*=0{}&*=0{}\ar@{-}[r]\ar@{-}[dl]&*=0{}\ar@{-}[r]\ar@{-}[dl]&*=0{}\ar@{-}[dl]\ar@{..}[d]
\\*=0{}&*=0{}\ar@{-}[r]\ar@{-}[dl]&*=0{}\ar@{-}[r]\ar@{-}[dl]&*=0{}&*=0{}\ar@{-}[r]\ar@{-}[dl]&*=0{}\ar@{-}[r]\ar@{-}[dl]&*=0{}\ar@{-}[r]\ar@{-}[dl]&*=0{}\ar@{-}[r]\ar@{-}[dl]&*=0{}\ar@{-}[r]\ar@{-}[dl]&*=0{}\ar@{-}[r]\ar@{-}[dl]&*=0{}\ar@{-}[dl]
\\*=0{}\ar@{-}[r]&*=0{}&*=0{}&*=0{}\ar@{-}[r]\ar@{-}[dl]&*=0{}\ar@{-}[r]\ar@{-}[dl]&*=0{}\ar@{-}[r]\ar@{-}[dl]&*=0{}\ar@{-}[r]&*=0{}\ar@{-}[r]&*=0{}\ar@{-}[r]&*=0{}
\\*=0{}&*=0{}&*=0{}\ar@{-}[r]\ar@{-}[dl]&*=0{}\ar@{-}[r]\ar@{-}[dl]&*=0{}\ar@{..}[d]
\\*=0{}&*=0{}\ar@{-}[r]\ar@{-}[dl]&*=0{}\ar@{-}[dl]&*=0{}&*=0{}\ar@{-}[r]\ar@{-}[dl]&*=0{}\ar@{-}[r]\ar@{-}[dl]&*=0{}\ar@{-}[r]\ar@{-}[dl]&*=0{}\ar@{-}[dl]
\\*=0{}\ar@{-}[r]&*=0{}&*=0{}&*=0{}\ar@{-}[r]\ar@{-}[dl]&*=0{}\ar@{-}[r]\ar@{-}[dl]&*=0{}\ar@{-}[r]\ar@{-}[dl]&*=0{}
\\*=0{}&*=0{}&*=0{}\ar@{-}[r]\ar@{-}[dl]&*=0{}\ar@{-}[r]\ar@{-}[dl]&*=0{}\ar@{-}[dl]\ar@{..}[d]
\\*=0{}&*=0{}\ar@{-}[r]&*=0{}\ar@{-}[r]&*=0{}&*=0{}\ar@{-}[r]\ar@{-}[dl]&*=0{}\ar@{-}[r]\ar@{-}[dl]&*=0{}\ar@{-}[r]\ar@{-}[dl]&*=0{}\ar@{-}[r]\ar@{-}[dl]&*=0{}\ar@{-}[dl]
\\*=0{}&*=0{}&*=0{}&*=0{}\ar@{-}[r]&*=0{}\ar@{-}[r]&*=0{}\ar@{-}[r]&*=0{}\ar@{-}[r]&*=0{}
\\*=0{}}
\text{ and }B=\xymatrix @W=1.0pc @H=1.0pc @R=1.0pc @C=1.0pc
{*=0{}&*=0{}&*=0{}&*=0{}\ar@{-}[r]\ar@{-}[dl]&*=0{}\ar@{-}[r]\ar@{-}[dl]&*=0{}\ar@{-}[dl]
\\*=0{}&*=0{}&*=0{}\ar@{-}[r]\ar@{-}[dl]&*=0{}\ar@{-}[r]\ar@{-}[dl]\ar@{..}[dd]&*=0{}
\\*=0{}&*=0{}\ar@{-}[r]\ar@{-}[dl]&*=0{}\ar@{-}[dl]
\\*=0{}\ar@{-}[r]&*=0{}&*=0{}&*=0{}\ar@{-}[r]\ar@{-}[dl]&*=0{}\ar@{-}[r]\ar@{-}[dl]&*=0{}\ar@{-}[r]\ar@{-}[dl]&*=0{}\ar@{-}[dl]
\\*=0{}&*=0{}&*=0{}\ar@{-}[r]\ar@{-}[dl]&*=0{}\ar@{-}[r]\ar@{-}[dl]\ar@{..}[dd]&*=0{}\ar@{-}[r]&*=0{}
\\*=0{}&*=0{}\ar@{-}[r]&*=0{}
\\
*=0{}&*=0{}&*=0{}&*=0{}\ar@{-}[r]\ar@{-}[dl]&*=0{}\ar@{-}[r]\ar@{-}[dl]&*=0{}\ar@{-}[dl]
\\*=0{}&*=0{}&*=0{}\ar@{-}[r]&*=0{}\ar@{-}[r]&*=0{}
\\*=0{}}
\]
then
\[
A\otimes_1 B =\xymatrix @W=1.0pc @H=1.0pc @R=1.0pc @C=1.0pc
{*=0{}&*=0{}&*=0{}&*=0{}&*=0{}\ar@{-}[r]\ar@{-}[dl]&*=0{}\ar@{-}[r]\ar@{-}[dl]&*=0{}\ar@{-}[r]\ar@{-}[dl]&*=0{}\ar@{-}[r]\ar@{-}[dl]&*=0{}\ar@{-}[dl]
\\*=0{}&*=0{}&*=0{}&*=0{}\ar@{-}[r]\ar@{-}[dl]&*=0{}\ar@{-}[r]\ar@{-}[dl]&*=0{}\ar@{-}[r]\ar@{-}[dl]&*=0{}\ar@{-}[r]&*=0{}
\\*=0{}&*=0{}&*=0{}\ar@{-}[r]\ar@{-}[dl]&*=0{}\ar@{-}[r]\ar@{-}[dl]&*=0{}\ar@{-}[dl]\ar@{..}[d]
\\*=0{}&*=0{}\ar@{-}[r]\ar@{-}[dl]&*=0{}\ar@{-}[r]\ar@{-}[dl]&*=0{}&*=0{}\ar@{-}[r]\ar@{-}[dl]&*=0{}\ar@{-}[r]\ar@{-}[dl]&*=0{}\ar@{-}[r]\ar@{-}[dl]&*=0{}\ar@{-}[r]\ar@{-}[dl]&*=0{}\ar@{-}[r]\ar@{-}[dl]&*=0{}\ar@{-}[r]\ar@{-}[dl]&*=0{}\ar@{-}[dl]
\\*=0{}\ar@{-}[r]&*=0{}&*=0{}&*=0{}\ar@{-}[r]\ar@{-}[dl]&*=0{}\ar@{-}[r]\ar@{-}[dl]&*=0{}\ar@{-}[r]\ar@{-}[dl]&*=0{}\ar@{-}[r]&*=0{}\ar@{-}[r]&*=0{}\ar@{-}[r]&*=0{}
\\*=0{}&*=0{}&*=0{}\ar@{-}[r]\ar@{-}[dl]&*=0{}\ar@{-}[r]\ar@{-}[dl]&*=0{}\ar@{..}[d]
\\*=0{}&*=0{}\ar@{-}[r]\ar@{-}[dl]&*=0{}\ar@{-}[dl]&*=0{}&*=0{}\ar@{-}[r]\ar@{-}[dl]&*=0{}\ar@{-}[r]\ar@{-}[dl]&*=0{}\ar@{-}[r]\ar@{-}[dl]&*=0{}\ar@{-}[dl]
\\*=0{}\ar@{-}[r]&*=0{}&*=0{}&*=0{}\ar@{-}[r]\ar@{-}[dl]&*=0{}\ar@{-}[r]\ar@{-}[dl]&*=0{}\ar@{-}[r]\ar@{-}[dl]&*=0{}
\\*=0{}&*=0{}&*=0{}\ar@{-}[r]\ar@{-}[dl]&*=0{}\ar@{-}[r]\ar@{-}[dl]&*=0{}\ar@{-}[dl]\ar@{..}[d]
\\*=0{}&*=0{}\ar@{-}[r]&*=0{}\ar@{-}[r]&*=0{}&*=0{}\ar@{-}[r]\ar@{-}[dl]&*=0{}\ar@{-}[r]\ar@{-}[dl]&*=0{}\ar@{-}[r]\ar@{-}[dl]&*=0{}\ar@{-}[r]\ar@{-}[dl]&*=0{}\ar@{-}[dl]
\\*=0{}&*=0{}&*=0{}&*=0{}\ar@{-}[r]&*=0{}\ar@{-}[r]&*=0{}\ar@{-}[r]&*=0{}\ar@{-}[r]&*=0{}
\\*=0{}}
A \otimes_2 B=\xymatrix @W=1.0pc @H=1.0pc @R=1.0pc @C=1.0pc
{*=0{}&*=0{}&*=0{}&*=0{}&*=0{}\ar@{-}[r]\ar@{-}[dl]&*=0{}\ar@{-}[r]\ar@{-}[dl]&*=0{}\ar@{-}[r]\ar@{-}[dl]&*=0{}\ar@{-}[r]\ar@{-}[dl]&*=0{}\ar@{-}[dl]
\\*=0{}&*=0{}&*=0{}&*=0{}\ar@{-}[r]\ar@{-}[dl]&*=0{}\ar@{-}[r]\ar@{-}[dl]&*=0{}\ar@{-}[r]\ar@{-}[dl]&*=0{}\ar@{-}[r]&*=0{}
\\*=0{}&*=0{}&*=0{}\ar@{-}[r]\ar@{-}[dl]&*=0{}\ar@{-}[r]\ar@{-}[dl]&*=0{}\ar@{-}[dl]\ar@{..}[d]
\\*=0{}&*=0{}\ar@{-}[r]\ar@{-}[dl]&*=0{}\ar@{-}[r]\ar@{-}[dl]&*=0{}&*=0{}\ar@{-}[r]\ar@{-}[dl]&*=0{}\ar@{-}[r]\ar@{-}[dl]&*=0{}\ar@{-}[r]\ar@{-}[dl]&*=0{}\ar@{-}[r]\ar@{-}[dl]&*=0{}\ar@{-}[r]\ar@{-}[dl]&*=0{}\ar@{-}[r]\ar@{-}[dl]&*=0{}\ar@{-}[dl]
\\*=0{}\ar@{-}[r]&*=0{}&*=0{}&*=0{}\ar@{-}[r]\ar@{-}[dl]&*=0{}\ar@{-}[r]\ar@{-}[dl]&*=0{}\ar@{-}[r]\ar@{-}[dl]&*=0{}\ar@{-}[r]&*=0{}\ar@{-}[r]&*=0{}\ar@{-}[r]&*=0{}
\\*=0{}&*=0{}&*=0{}\ar@{-}[r]\ar@{-}[dl]&*=0{}\ar@{-}[r]\ar@{-}[dl]&*=0{}\ar@{..}[d]
\\*=0{}&*=0{}\ar@{-}[r]\ar@{-}[dl]&*=0{}\ar@{-}[dl]&*=0{}&*=0{}\ar@{-}[r]\ar@{-}[dl]&*=0{}\ar@{-}[r]\ar@{-}[dl]&*=0{}\ar@{-}[r]\ar@{-}[dl]&*=0{}\ar@{-}[dl]
\\*=0{}\ar@{-}[r]&*=0{}&*=0{}&*=0{}\ar@{-}[r]\ar@{-}[dl]&*=0{}\ar@{-}[r]\ar@{-}[dl]&*=0{}\ar@{-}[r]\ar@{-}[dl]&*=0{}
\\*=0{}&*=0{}&*=0{}\ar@{-}[r]\ar@{-}[dl]&*=0{}\ar@{-}[r]\ar@{-}[dl]&*=0{}\ar@{-}[dl]\ar@{..}[d]
\\*=0{}&*=0{}\ar@{-}[r]&*=0{}\ar@{-}[r]&*=0{}&*=0{}\ar@{-}[r]\ar@{-}[dl]&*=0{}\ar@{-}[r]\ar@{-}[dl]&*=0{}\ar@{-}[dl]
\\*=0{}&*=0{}&*=0{}&*=0{}\ar@{-}[r]\ar@{-}[dl]&*=0{}\ar@{-}[r]\ar@{-}[dl]\ar@{..}[d]&*=0{}
\\*=0{}&*=0{}&*=0{}\ar@{-}[r]\ar@{-}[dl]&*=0{}\ar@{-}[dl]&*=0{}\ar@{-}[r]\ar@{-}[dl]&*=0{}\ar@{-}[r]\ar@{-}[dl]&*=0{}\ar@{-}[r]\ar@{-}[dl]&*=0{}\ar@{-}[dl]
\\*=0{}&*=0{}\ar@{-}[r]&*=0{}&*=0{}\ar@{-}[r]\ar@{-}[dl]&*=0{}\ar@{-}[r]\ar@{-}[dl]\ar@{..}[d]&*=0{}\ar@{-}[r]&*=0{}
\\*=0{}&*=0{}&*=0{}\ar@{-}[r]&*=0{}&*=0{}\ar@{-}[r]\ar@{-}[dl]&*=0{}\ar@{-}[r]\ar@{-}[dl]&*=0{}\ar@{-}[r]\ar@{-}[dl]&*=0{}\ar@{-}[r]\ar@{-}[dl]&*=0{}\ar@{-}[dl]
\\*=0{}&*=0{}&*=0{}&*=0{}\ar@{-}[r]&*=0{}\ar@{-}[r]\ar@{..}[d]&*=0{}\ar@{-}[r]&*=0{}\ar@{-}[r]&*=0{}
\\
*=0{}&*=0{}&*=0{}&*=0{}&*=0{}\ar@{-}[r]\ar@{-}[dl]&*=0{}\ar@{-}[r]\ar@{-}[dl]&*=0{}\ar@{-}[dl]
\\*=0{}&*=0{}&*=0{}&*=0{}\ar@{-}[r]&*=0{}\ar@{-}[r]&*=0{}
\\*=0{}}
\]
\[
A \otimes_3 B=\xymatrix @W=1.0pc @H=1.0pc @R=1.0pc @C=1.0pc
{*=0{}&*=0{}&*=0{}&*=0{}&*=0{}\ar@{-}[r]\ar@{-}[dl]&*=0{}\ar@{-}[r]\ar@{-}[dl]&*=0{}\ar@{-}[r]\ar@{-}[dl]&*=0{}\ar@{-}[r]\ar@{-}[dl]&*=0{}\ar@{-}[r]\ar@{-}[dl]&*=0{}\ar@{-}[r]\ar@{-}[dl]&*=0{}\ar@{-}[dl]
\\*=0{}&*=0{}&*=0{}&*=0{}\ar@{-}[r]\ar@{-}[dl]&*=0{}\ar@{-}[r]\ar@{-}[dl]&*=0{}\ar@{-}[r]\ar@{-}[dl]&*=0{}\ar@{-}[r]\ar@{-}[dl]&*=0{}\ar@{-}[r]&*=0{}\ar@{-}[r]&*=0{}
\\*=0{}&*=0{}&*=0{}\ar@{-}[r]\ar@{-}[dl]&*=0{}\ar@{-}[r]\ar@{-}[dl]&*=0{}\ar@{-}[r]\ar@{-}[dl]&*=0{}\ar@{-}[dl]
\\*=0{}&*=0{}\ar@{-}[r]\ar@{-}[dl]&*=0{}\ar@{-}[r]\ar@{-}[dl]&*=0{}\ar@{-}[r]&*=0{}\ar@{..}[d]
\\*=0{}\ar@{-}[r]&*=0{}&*=0{}&*=0{}&*=0{}\ar@{-}[r]\ar@{-}[dl]&*=0{}\ar@{-}[r]\ar@{-}[dl]&*=0{}\ar@{-}[r]\ar@{-}[dl]&*=0{}\ar@{-}[r]\ar@{-}[dl]&*=0{}\ar@{-}[r]\ar@{-}[dl]&*=0{}\ar@{-}[r]\ar@{-}[dl]&*=0{}\ar@{-}[r]\ar@{-}[dl]&*=0{}\ar@{-}[r]\ar@{-}[dl]&*=0{}\ar@{-}[r]\ar@{-}[dl]&*=0{}\ar@{-}[dl]
\\*=0{}&*=0{}&*=0{}&*=0{}\ar@{-}[r]\ar@{-}[dl]&*=0{}\ar@{-}[r]\ar@{-}[dl]&*=0{}\ar@{-}[r]\ar@{-}[dl]&*=0{}\ar@{-}[r]\ar@{-}[dl]&*=0{}\ar@{-}[r]&*=0{}\ar@{-}[r]&*=0{}\ar@{-}[r]&*=0{}\ar@{-}[r]&*=0{}\ar@{-}[r]&*=0{}
\\*=0{}&*=0{}&*=0{}\ar@{-}[r]\ar@{-}[dl]&*=0{}\ar@{-}[r]\ar@{-}[dl]&*=0{}\ar@{-}[r]\ar@{..}[d]&*=0{}
\\*=0{}&*=0{}\ar@{-}[r]\ar@{-}[dl]&*=0{}\ar@{-}[dl]&*=0{}&*=0{}\ar@{-}[r]\ar@{-}[dl]&*=0{}\ar@{-}[r]\ar@{-}[dl]&*=0{}\ar@{-}[r]\ar@{-}[dl]&*=0{}\ar@{-}[r]\ar@{-}[dl]&*=0{}\ar@{-}[r]\ar@{-}[dl]&*=0{}\ar@{-}[dl]
\\*=0{}\ar@{-}[r]&*=0{}&*=0{}&*=0{}\ar@{-}[r]\ar@{-}[dl]&*=0{}\ar@{-}[r]\ar@{-}[dl]&*=0{}\ar@{-}[r]\ar@{-}[dl]&*=0{}\ar@{-}[r]&*=0{}\ar@{-}[r]&*=0{}
\\*=0{}&*=0{}&*=0{}\ar@{-}[r]\ar@{-}[dl]&*=0{}\ar@{-}[r]\ar@{-}[dl]&*=0{}\ar@{-}[dl]\ar@{..}[d]
\\*=0{}&*=0{}\ar@{-}[r]&*=0{}\ar@{-}[r]&*=0{}&*=0{}\ar@{-}[r]\ar@{-}[dl]&*=0{}\ar@{-}[r]\ar@{-}[dl]&*=0{}\ar@{-}[r]\ar@{-}[dl]&*=0{}\ar@{-}[r]\ar@{-}[dl]&*=0{}\ar@{-}[dl]
\\*=0{}&*=0{}&*=0{}&*=0{}\ar@{-}[r]&*=0{}\ar@{-}[r]&*=0{}\ar@{-}[r]&*=0{}\ar@{-}[r]&*=0{}
\\*=0{}}
\text{ and }
A \otimes_4 B=\xymatrix @W=1.0pc @H=1.0pc @R=1.0pc @C=1.0pc
{*=0{}&*=0{}&*=0{}&*=0{}&*=0{}&*=0{}&*=0{}&*=0{}\ar@{-}[r]\ar@{-}[dl]&*=0{}\ar@{-}[r]\ar@{-}[dl]&*=0{}\ar@{-}[r]\ar@{-}[dl]&*=0{}\ar@{-}[r]\ar@{-}[dl]&*=0{}\ar@{-}[dl]
\\*=0{}&*=0{}&*=0{}&*=0{}&*=0{}&*=0{}&*=0{}\ar@{-}[r]\ar@{-}[dl]&*=0{}\ar@{-}[r]\ar@{-}[dl]&*=0{}\ar@{-}[r]\ar@{-}[dl]&*=0{}\ar@{-}[r]&*=0{}
\\*=0{}&*=0{}&*=0{}&*=0{}&*=0{}&*=0{}\ar@{-}[r]\ar@{-}[dl]&*=0{}\ar@{-}[r]\ar@{-}[dl]&*=0{}\ar@{-}[dl]\ar@{..}[d]
\\*=0{}&*=0{}&*=0{}&*=0{}&*=0{}\ar@{-}[r]\ar@{-}[dl]&*=0{}\ar@{-}[r]\ar@{-}[dl]&*=0{}\ar@{-}[dl]&*=0{}\ar@{-}[r]\ar@{-}[dl]&*=0{}\ar@{-}[r]\ar@{-}[dl]&*=0{}\ar@{-}[r]\ar@{-}[dl]&*=0{}\ar@{-}[r]\ar@{-}[dl]&*=0{}\ar@{-}[r]\ar@{-}[dl]&*=0{}\ar@{-}[r]\ar@{-}[dl]&*=0{}\ar@{-}[dl]
\\*=0{}&*=0{}&*=0{}&*=0{}\ar@{-}[r]\ar@{-}[dl]&*=0{}\ar@{-}[r]\ar@{-}[dl]&*=0{}&*=0{}\ar@{-}[r]\ar@{-}[dl]&*=0{}\ar@{-}[r]\ar@{-}[dl]&*=0{}\ar@{-}[r]\ar@{-}[dl]&*=0{}\ar@{-}[r]\ar@{-}[dl]&*=0{}\ar@{-}[r]&*=0{}\ar@{-}[r]&*=0{}
\\*=0{}&*=0{}&*=0{}\ar@{-}[r]\ar@{-}[dl]&*=0{}\ar@{-}[dl]&*=0{}&*=0{}\ar@{-}[r]\ar@{-}[dl]&*=0{}\ar@{-}[r]\ar@{-}[dl]&*=0{}\ar@{-}[r]\ar@{-}[dl]\ar@{..}[d]&*=0{}
\\*=0{}&*=0{}\ar@{-}[r]\ar@{-}[dl]&*=0{}\ar@{-}[dl]&*=0{}&*=0{}\ar@{-}[r]\ar@{-}[dl]&*=0{}\ar@{-}[r]\ar@{-}[dl]&*=0{}&*=0{}\ar@{-}[r]\ar@{-}[dl]&*=0{}\ar@{-}[r]\ar@{-}[dl]&*=0{}\ar@{-}[r]\ar@{-}[dl]&*=0{}\ar@{-}[dl]
\\*=0{}\ar@{-}[r]&*=0{}&*=0{}&*=0{}\ar@{-}[r]\ar@{-}[dl]&*=0{}\ar@{-}[dl]&*=0{}&*=0{}\ar@{-}[r]\ar@{-}[dl]&*=0{}\ar@{-}[r]\ar@{-}[dl]&*=0{}\ar@{-}[r]\ar@{-}[dl]&*=0{}
\\*=0{}&*=0{}&*=0{}\ar@{-}[r]\ar@{-}[dl]&*=0{}\ar@{-}[dl]&*=0{}&*=0{}\ar@{-}[r]\ar@{-}[dl]&*=0{}\ar@{-}[r]\ar@{-}[dl]&*=0{}\ar@{-}[dl]\ar@{..}[d]
\\*=0{}&*=0{}\ar@{-}[r]&*=0{}&*=0{}&*=0{}\ar@{-}[r]\ar@{-}[dl]&*=0{}\ar@{-}[r]\ar@{-}[dl]&*=0{}\ar@{-}[dl]&*=0{}\ar@{-}[r]\ar@{-}[dl]&*=0{}\ar@{-}[r]\ar@{-}[dl]&*=0{}\ar@{-}[r]\ar@{-}[dl]&*=0{}\ar@{-}[r]\ar@{-}[dl]&*=0{}\ar@{-}[dl]
\\*=0{}&*=0{}&*=0{}&*=0{}\ar@{-}[r]&*=0{}\ar@{-}[r]&*=0{}&*=0{}\ar@{-}[r]&*=0{}\ar@{-}[r]&*=0{}\ar@{-}[r]&*=0{}\ar@{-}[r]&*=0{}
\\*=0{}}
\]

\begin{example}\label{3-castles}
It might be nice to retain the geometric picture of the products
of Young diagrams in terms of vertical and horizontal stacking, and
stacking in other directions as dimension increases. This is not
found in the just illustrated category, which relies on the merging
viewpoint.  The ``diagram stacking'' point of view is restored if
we restrict to 3-d  Young diagrams. We can represent these objects
as infinite matrices with finitely nonzero natural number entries,
and with monotone decreasing columns and rows.  We  require that
${A_n}_k$ be decreasing in $n$ for constant $k$, and decreasing in
$k$ for constant $n.$ We choose the sequence of rows to represent
the sequence of sequences, i.e. each row represents a Young diagram
which we draw as being parallel to the $xy$ plane.  This choice is
important because it determines the total ordering of matrices and
thus the morphisms of the category.  Thus $y$-axis stacking is
horizontal concatenation (disregarding  trailing zeroes) of matrices
followed by sorting the new longer rows (row merging).  $x$-axis
stacking is addition of matrices.  Now we define $z$-axis stacking
as vertical concatenation  of matrices followed by sorting the new
long columns (column merging).

\end{example}
Here is a visual example of the three new products, beginning with
$z$-axis stacking, labeled $\otimes_1$: if
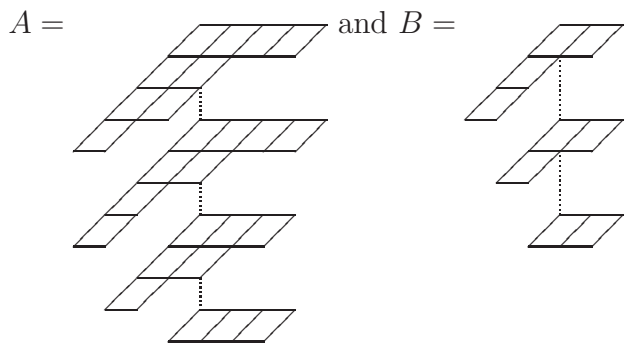
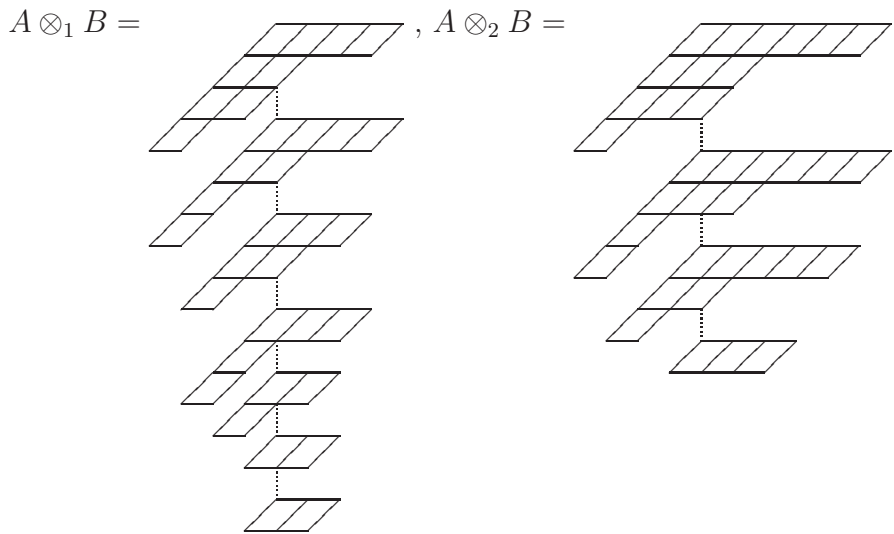
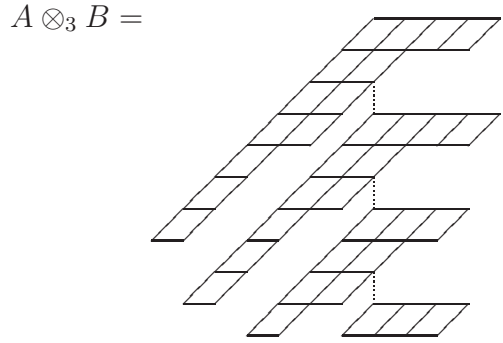
\[
A=\xymatrix @W=1.0pc @H=1.0pc @R=1.0pc @C=1.0pc
{*=0{}&*=0{}&*=0{}&*=0{}&*=0{}\ar@{-}[r]\ar@{-}[dl]&*=0{}\ar@{-}[r]\ar@{-}[dl]&*=0{}\ar@{-}[r]\ar@{-}[dl]&*=0{}\ar@{-}[r]\ar@{-}[dl]&*=0{}\ar@{-}[dl]
\\*=0{}&*=0{}&*=0{}&*=0{}\ar@{-}[r]\ar@{-}[dl]&*=0{}\ar@{-}[r]\ar@{-}[dl]&*=0{}\ar@{-}[r]\ar@{-}[dl]&*=0{}\ar@{-}[r]&*=0{}
\\*=0{}&*=0{}&*=0{}\ar@{-}[r]\ar@{-}[dl]&*=0{}\ar@{-}[r]\ar@{-}[dl]&*=0{}\ar@{-}[dl]\ar@{..}[d]
\\*=0{}&*=0{}\ar@{-}[r]\ar@{-}[dl]&*=0{}\ar@{-}[r]\ar@{-}[dl]&*=0{}  &*=0{}\ar@{-}[r]\ar@{-}[dl]&*=0{}\ar@{-}[r]\ar@{-}[dl]&*=0{}\ar@{-}[r]\ar@{-}[dl]&*=0{}\ar@{-}[r]\ar@{-}[dl]&*=0{}\ar@{-}[dl]
\\*=0{}\ar@{-}[r]&*=0{} &*=0{}&*=0{}\ar@{-}[r]\ar@{-}[dl]&*=0{}\ar@{-}[r]\ar@{-}[dl]&*=0{}\ar@{-}[r]\ar@{-}[dl]&*=0{}\ar@{-}[r]&*=0{}
\\*=0{}&*=0{}&*=0{}\ar@{-}[r]\ar@{-}[dl]&*=0{}\ar@{-}[r]\ar@{-}[dl]&*=0{}\ar@{..}[d]
\\*=0{}&*=0{}\ar@{-}[r]\ar@{-}[dl]&*=0{}\ar@{-}[dl]&*=0{}&*=0{}\ar@{-}[r]\ar@{-}[dl]&*=0{}\ar@{-}[r]\ar@{-}[dl]&*=0{}\ar@{-}[r]\ar@{-}[dl]&*=0{}\ar@{-}[dl]
\\*=0{}\ar@{-}[r]&*=0{}&*=0{}&*=0{}\ar@{-}[r]\ar@{-}[dl]&*=0{}\ar@{-}[r]\ar@{-}[dl]&*=0{}\ar@{-}[r]\ar@{-}[dl]&*=0{}
\\*=0{}&*=0{}&*=0{}\ar@{-}[r]\ar@{-}[dl]&*=0{}\ar@{-}[r]\ar@{-}[dl]&*=0{}\ar@{..}[d]
\\*=0{}&*=0{}\ar@{-}[r]&*=0{}&*=0{}&*=0{}\ar@{-}[r]\ar@{-}[dl]&*=0{}\ar@{-}[r]\ar@{-}[dl]&*=0{}\ar@{-}[r]\ar@{-}[dl]&*=0{}\ar@{-}[dl]&*=0{}
\\*=0{}&*=0{}&*=0{}&*=0{}\ar@{-}[r]&*=0{}\ar@{-}[r]&*=0{}\ar@{-}[r]&*=0{}
\\*=0{}}
\text{ and }B=\xymatrix @W=1.0pc @H=1.0pc @R=1.0pc @C=1.0pc
{*=0{}&*=0{}&*=0{}&*=0{}\ar@{-}[r]\ar@{-}[dl]&*=0{}\ar@{-}[r]\ar@{-}[dl]&*=0{}\ar@{-}[dl]
\\*=0{}&*=0{}&*=0{}\ar@{-}[r]\ar@{-}[dl]&*=0{}\ar@{-}[r]\ar@{-}[dl]\ar@{..}[dd]&*=0{}
\\*=0{}&*=0{}\ar@{-}[r]\ar@{-}[dl]&*=0{}\ar@{-}[dl]
\\*=0{}\ar@{-}[r]&*=0{}&*=0{}&*=0{}\ar@{-}[r]\ar@{-}[dl]&*=0{}\ar@{-}[r]\ar@{-}[dl]&*=0{}\ar@{-}[dl]&*=0{}
\\*=0{}&*=0{}&*=0{}\ar@{-}[r]\ar@{-}[dl]&*=0{}\ar@{-}[r]\ar@{-}[dl]\ar@{..}[dd]&*=0{}
\\*=0{}&*=0{}\ar@{-}[r]&*=0{}
\\
*=0{}&*=0{}&*=0{}&*=0{}\ar@{-}[r]\ar@{-}[dl]&*=0{}\ar@{-}[r]\ar@{-}[dl]&*=0{}\ar@{-}[dl]
\\*=0{}&*=0{}&*=0{}\ar@{-}[r]&*=0{}\ar@{-}[r]&*=0{}
\\*=0{}}
\]
then we let
\[
A \otimes_1 B=\xymatrix @W=1.0pc @H=1.0pc @R=1.0pc @C=1.0pc
{*=0{}&*=0{}&*=0{}&*=0{}&*=0{}\ar@{-}[r]\ar@{-}[dl]&*=0{}\ar@{-}[r]\ar@{-}[dl]&*=0{}\ar@{-}[r]\ar@{-}[dl]&*=0{}\ar@{-}[r]\ar@{-}[dl]&*=0{}\ar@{-}[dl]
\\*=0{}&*=0{}&*=0{}&*=0{}\ar@{-}[r]\ar@{-}[dl]&*=0{}\ar@{-}[r]\ar@{-}[dl]&*=0{}\ar@{-}[r]\ar@{-}[dl]&*=0{}\ar@{-}[r]&*=0{}
\\*=0{}&*=0{}&*=0{}\ar@{-}[r]\ar@{-}[dl]&*=0{}\ar@{-}[r]\ar@{-}[dl]&*=0{}\ar@{-}[dl]\ar@{..}[d]
\\*=0{}&*=0{}\ar@{-}[r]\ar@{-}[dl]&*=0{}\ar@{-}[r]\ar@{-}[dl]&*=0{}&*=0{}\ar@{-}[r]\ar@{-}[dl]&*=0{}\ar@{-}[r]\ar@{-}[dl]&*=0{}\ar@{-}[r]\ar@{-}[dl]&*=0{}\ar@{-}[r]\ar@{-}[dl]&*=0{}\ar@{-}[dl]
\\*=0{}\ar@{-}[r]&*=0{}&*=0{}&*=0{}\ar@{-}[r]\ar@{-}[dl]&*=0{}\ar@{-}[r]\ar@{-}[dl]&*=0{}\ar@{-}[r]\ar@{-}[dl]&*=0{}\ar@{-}[r]&*=0{}
\\*=0{}&*=0{}&*=0{}\ar@{-}[r]\ar@{-}[dl]&*=0{}\ar@{-}[r]\ar@{-}[dl]&*=0{}\ar@{..}[d]
\\*=0{}&*=0{}\ar@{-}[r]\ar@{-}[dl]&*=0{}\ar@{-}[dl]&*=0{}&*=0{}\ar@{-}[r]\ar@{-}[dl]&*=0{}\ar@{-}[r]\ar@{-}[dl]&*=0{}\ar@{-}[r]\ar@{-}[dl]&*=0{}\ar@{-}[dl]
\\*=0{}\ar@{-}[r]&*=0{}&*=0{}&*=0{}\ar@{-}[r]\ar@{-}[dl]&*=0{}\ar@{-}[r]\ar@{-}[dl]&*=0{}\ar@{-}[r]\ar@{-}[dl]&*=0{}
\\*=0{}&*=0{}&*=0{}\ar@{-}[r]\ar@{-}[dl]&*=0{}\ar@{-}[r]\ar@{-}[dl]&*=0{}\ar@{..}[d]
\\*=0{}&*=0{}\ar@{-}[r]&*=0{}&*=0{}&*=0{}\ar@{-}[r]\ar@{-}[dl]&*=0{}\ar@{-}[r]\ar@{-}[dl]&*=0{}\ar@{-}[r]\ar@{-}[dl]&*=0{}\ar@{-}[dl]&*=0{}
\\*=0{}&*=0{}&*=0{}&*=0{}\ar@{-}[r]\ar@{-}[dl]&*=0{}\ar@{-}[r]\ar@{-}[dl]\ar@{..}[d]&*=0{}\ar@{-}[r]&*=0{}&*=0{}
\\*=0{}&*=0{}&*=0{}\ar@{-}[r]\ar@{-}[dl]&*=0{}\ar@{-}[dl]&*=0{}\ar@{-}[r]\ar@{-}[dl]&*=0{}\ar@{-}[r]\ar@{-}[dl]&*=0{}\ar@{-}[dl]&*=0{}
\\*=0{}&*=0{}\ar@{-}[r]&*=0{}&*=0{}\ar@{-}[r]\ar@{-}[dl]&*=0{}\ar@{-}[r]\ar@{-}[dl]\ar@{..}[d]&*=0{}
\\*=0{}&*=0{}&*=0{}\ar@{-}[r]&*=0{}&*=0{}\ar@{-}[r]\ar@{-}[dl]&*=0{}\ar@{-}[r]\ar@{-}[dl]&*=0{}\ar@{-}[dl]
\\*=0{}&*=0{}&*=0{}&*=0{}\ar@{-}[r]&*=0{}\ar@{-}[r]\ar@{..}[d]&*=0{}&*=0{}
\\
*=0{}&*=0{}&*=0{}&*=0{}&*=0{}\ar@{-}[r]\ar@{-}[dl]&*=0{}\ar@{-}[r]\ar@{-}[dl]&*=0{}\ar@{-}[dl]
\\*=0{}&*=0{}&*=0{}&*=0{}\ar@{-}[r]&*=0{}\ar@{-}[r]&*=0{}
\\*=0{}}
\text{ , }
A \otimes_2 B=\xymatrix @W=1.0pc @H=1.0pc @R=1.0pc @C=1.0pc
{*=0{}&*=0{}&*=0{}&*=0{}&*=0{}\ar@{-}[r]\ar@{-}[dl]&*=0{}\ar@{-}[r]\ar@{-}[dl]&*=0{}\ar@{-}[r]\ar@{-}[dl]&*=0{}\ar@{-}[r]\ar@{-}[dl] &*=0{}\ar@{-}[r]\ar@{-}[dl] &*=0{}\ar@{-}[r]\ar@{-}[dl] &*=0{}\ar@{-}[dl]
\\*=0{}&*=0{}&*=0{}&*=0{}\ar@{-}[r]\ar@{-}[dl]&*=0{}\ar@{-}[r]\ar@{-}[dl]&*=0{}\ar@{-}[r]\ar@{-}[dl] &*=0{}\ar@{-}[r]\ar@{-}[dl] &*=0{}\ar@{-}[r] &*=0{}\ar@{-}[r]&*=0{}
\\*=0{}&*=0{}&*=0{}\ar@{-}[r]\ar@{-}[dl]&*=0{}\ar@{-}[r]\ar@{-}[dl]&*=0{}\ar@{-}[r]\ar@{-}[dl]&*=0{}\ar@{-}[dl]
\\*=0{}&*=0{}\ar@{-}[r]\ar@{-}[dl]&*=0{}\ar@{-}[r]\ar@{-}[dl]&*=0{}\ar@{-}[r]&*=0{}\ar@{..}[d]
\\*=0{}\ar@{-}[r]&*=0{}&*=0{}&*=0{}&*=0{}\ar@{-}[r]\ar@{-}[dl]&*=0{}\ar@{-}[r]\ar@{-}[dl]&*=0{}\ar@{-}[r]\ar@{-}[dl]&*=0{}\ar@{-}[r]\ar@{-}[dl]&*=0{}\ar@{-}[r]\ar@{-}[dl]&*=0{}\ar@{-}[r]\ar@{-}[dl]&*=0{}\ar@{-}[dl]
\\*=0{}&*=0{}&*=0{}&*=0{}\ar@{-}[r]\ar@{-}[dl]&*=0{}\ar@{-}[r]\ar@{-}[dl]&*=0{}\ar@{-}[r]\ar@{-}[dl]&*=0{}\ar@{-}[r]\ar@{-}[dl]&*=0{}\ar@{-}[r]&*=0{}\ar@{-}[r]&*=0{}
\\*=0{}&*=0{}&*=0{}\ar@{-}[r]\ar@{-}[dl]&*=0{}\ar@{-}[r]\ar@{-}[dl]&*=0{}\ar@{-}[r]\ar@{..}[d]&*=0{}
\\*=0{}&*=0{}\ar@{-}[r]\ar@{-}[dl]&*=0{}\ar@{-}[dl]&*=0{}&*=0{}\ar@{-}[r]\ar@{-}[dl]&*=0{}\ar@{-}[r]\ar@{-}[dl]&*=0{}\ar@{-}[r]\ar@{-}[dl]&*=0{}\ar@{-}[r]\ar@{-}[dl]&*=0{}\ar@{-}[r]\ar@{-}[dl]&*=0{}\ar@{-}[dl]
\\*=0{}\ar@{-}[r]&*=0{}&*=0{}&*=0{}\ar@{-}[r]\ar@{-}[dl]&*=0{}\ar@{-}[r]\ar@{-}[dl]&*=0{}\ar@{-}[r]\ar@{-}[dl]&*=0{}\ar@{-}[r]&*=0{}\ar@{-}[r]&*=0{}
\\*=0{}&*=0{}&*=0{}\ar@{-}[r]\ar@{-}[dl]&*=0{}\ar@{-}[r]\ar@{-}[dl]&*=0{}\ar@{..}[d]
\\*=0{}&*=0{}\ar@{-}[r]&*=0{}&*=0{}&*=0{}\ar@{-}[r]\ar@{-}[dl]&*=0{}\ar@{-}[r]\ar@{-}[dl]&*=0{}\ar@{-}[r]\ar@{-}[dl]&*=0{}\ar@{-}[dl]&*=0{}
\\*=0{}&*=0{}&*=0{}&*=0{}\ar@{-}[r]&*=0{}\ar@{-}[r]&*=0{}\ar@{-}[r]&*=0{}
\\*=0{}}
\]
and
\[
A \otimes_3 B=\xymatrix @W=1.0pc @H=1.0pc @R=1.0pc @C=1.0pc
{*=0{}&*=0{}&*=0{}&*=0{}&*=0{}&*=0{}&*=0{}&*=0{}\ar@{-}[r]\ar@{-}[dl]&*=0{}\ar@{-}[r]\ar@{-}[dl]&*=0{}\ar@{-}[r]\ar@{-}[dl]&*=0{}\ar@{-}[r]\ar@{-}[dl]&*=0{}\ar@{-}[dl]
\\*=0{}&*=0{}&*=0{}&*=0{}&*=0{}&*=0{}&*=0{}\ar@{-}[r]\ar@{-}[dl]&*=0{}\ar@{-}[r]\ar@{-}[dl]&*=0{}\ar@{-}[r]\ar@{-}[dl]&*=0{}\ar@{-}[r]&*=0{}&*=0{}&*=0{}
\\*=0{}&*=0{}&*=0{}&*=0{}&*=0{}&*=0{}\ar@{-}[r]\ar@{-}[dl]&*=0{}\ar@{-}[r]\ar@{-}[dl]&*=0{}\ar@{-}[dl]\ar@{..}[d]&*=0{}&*=0{}&*=0{}&*=0{}
\\*=0{}&*=0{}&*=0{}&*=0{}&*=0{}\ar@{-}[r]\ar@{-}[dl]&*=0{}\ar@{-}[r]\ar@{-}[dl]&*=0{}\ar@{-}[dl]&*=0{}\ar@{-}[r]\ar@{-}[dl]&*=0{}\ar@{-}[r]\ar@{-}[dl]&*=0{}\ar@{-}[r]\ar@{-}[dl]&*=0{}\ar@{-}[r]\ar@{-}[dl]&*=0{}\ar@{-}[dl]&*=0{}&*=0{}&*=0{}&*=0{}
\\*=0{}&*=0{}&*=0{}&*=0{}\ar@{-}[r]\ar@{-}[dl]&*=0{}\ar@{-}[r]\ar@{-}[dl]&*=0{}&*=0{}\ar@{-}[r]\ar@{-}[dl]&*=0{}\ar@{-}[r]\ar@{-}[dl]&*=0{}\ar@{-}[r]\ar@{-}[dl]&*=0{}\ar@{-}[r]&*=0{}&*=0{}&*=0{}&*=0{}&*=0{}
\\*=0{}&*=0{}&*=0{}\ar@{-}[r]\ar@{-}[dl]&*=0{}\ar@{-}[dl]&*=0{}&*=0{}\ar@{-}[r]\ar@{-}[dl]&*=0{}\ar@{-}[r]\ar@{-}[dl]&*=0{}\ar@{-}[dl]\ar@{..}[d]&*=0{}&*=0{}&*=0{}&*=0{}
\\*=0{}&*=0{}\ar@{-}[r]\ar@{-}[dl]&*=0{}\ar@{-}[dl]&*=0{}&*=0{}\ar@{-}[r]\ar@{-}[dl]&*=0{}\ar@{-}[r]\ar@{-}[dl]&*=0{}&*=0{}\ar@{-}[r]\ar@{-}[dl]&*=0{}\ar@{-}[r]\ar@{-}[dl]&*=0{}\ar@{-}[r]\ar@{-}[dl]&*=0{}\ar@{-}[dl]&*=0{}&*=0{}&*=0{}&*=0{}
\\*=0{}\ar@{-}[r]&*=0{}&*=0{}&*=0{}\ar@{-}[r]\ar@{-}[dl]&*=0{}\ar@{-}[dl]&*=0{}&*=0{}\ar@{-}[r]\ar@{-}[dl]&*=0{}\ar@{-}[r]\ar@{-}[dl]&*=0{}\ar@{-}[r]\ar@{-}[dl]&*=0{}&*=0{}&*=0{}&*=0{}&*=0{}
\\*=0{}&*=0{}&*=0{}\ar@{-}[r]\ar@{-}[dl]&*=0{}\ar@{-}[dl]&*=0{}&*=0{}\ar@{-}[r]\ar@{-}[dl]&*=0{}\ar@{-}[r]\ar@{-}[dl]&*=0{}\ar@{-}[dl]\ar@{..}[d]&*=0{}&*=0{}&*=0{}&*=0{}
\\*=0{}&*=0{}\ar@{-}[r]&*=0{}&*=0{}&*=0{}\ar@{-}[r]\ar@{-}[dl]&*=0{}\ar@{-}[r]\ar@{-}[dl]&*=0{}&*=0{}\ar@{-}[r]\ar@{-}[dl]&*=0{}\ar@{-}[r]\ar@{-}[dl]&*=0{}\ar@{-}[r]\ar@{-}[dl]&*=0{}\ar@{-}[dl]&*=0{}&*=0{}&*=0{}&*=0{}&*=0{}
\\*=0{}&*=0{}&*=0{}&*=0{}\ar@{-}[r]&*=0{}&*=0{}&*=0{}\ar@{-}[r]&*=0{}\ar@{-}[r]&*=0{}\ar@{-}[r]&*=0{}&*=0{}&*=0{}&*=0{}&*=0{}&*=0{}
\\*=0{}&*=0{}&*=0{}&*=0{}&*=0{}&*=0{}&*=0{}&*=0{}&*=0{}}
\]
Note that in this restricted setting of decreasing matrices the
lexicographic merging of sequences (rows) of two matrices does not
preserve the total decreasing property (decreasing in rows and
columns).

These three products just shown preserve the total sum of the entries
in both matrices, and do interact via interchanges to form the
structure of a 3-fold category. Renumbered, they are: $\otimes_1$
($z$-axis stacking) is the vertical concatenation of matrices
followed by sorting the new longer columns, $\otimes_2$ ($y$-axis
stacking) is horizontal concatenation of matrices followed by sorting
the new longer rows and $\otimes_3$  ($x$-axis stacking) is the
addition of matrices.  For comparison, here is the same example of
the products as just given above shown by matrices. Only the non-zero
entries of the matrices are shown.
\[
A = \left[\begin{array}{cccc}4&3&1&1 \\ 4&2&1&1 \\ 3&2&1 \\ 1&1&1\end{array} \right] B= \left[\begin{array}{cc}3&1 \\ 2&1 \\ 1&1 \end{array} \right]
\]
\[
A \otimes_1 B = \left[\begin{array}{cccc} 4&3&1&1 \\ 4&2&1&1 \\ 3&2&1 \\ 3&1&1 \\ 2&1 \\ 1&1 \\ 1&1\end{array} \right]
A\otimes_2 B = \left[\begin{array}{cccccc} 4&3&3&1&1&1 \\ 4&2&2&1&1&1 \\ 3&2&1&1&1 \\ 1&1&1 \end{array} \right]
\text{ and } A\otimes_3 B = \left[\begin{array}{cccc} 7&4&1&1 \\ 6&3&1&1 \\ 4&3&1 \\ 1&1&1\end{array} \right]
\]

\begin{theorem}\label{3dY}
The category of 3-d  Young diagrams with lexicographic ordering and
the products just described possesses the structure  of a 3-fold
monoidal category.
\end{theorem}
\begin{proof}
We already have existence of $\eta^{23}$ by the previous argument
about pointwise application of two interchanging products. To show
existence of $\eta^{13}:(A\otimes_3B)\otimes_1(C\otimes_3D)\to
(A\otimes_1C)\otimes_3(B\otimes_1D)$ we need to check that sorting
each of the columns of two pairs of vertically concatenated matrices
before pointwise adding gives a larger lexicographic result with
respect to rows than adding first and then sorting columns. This
follows from Lemma~\ref{sort}, applied to each pair of sequences
which are the $n^{th}$ columns in the two new matrices formed by
vertically concatenating  $A$ and $C$ and respectively $B$ and $D$,
padded with zeroes so that adding the new matrices results in adding
$A$ and $B$ and respectively $C$ and $D$. From  the lemma then we
have that $(A\otimes_1C)\otimes_3(B\otimes_1D)$ gives a result whose
$n^{th}$ column is lexicographically greater than or equal to the
$n^{th}$ column of $(A\otimes_3B)\otimes_1(C\otimes_3D).$ This
implies that either the pairs of respective columns are each equal
sequences or that there is some least row $i$ and column $j$ such
that all the pairs of columns are identical  in rows less than $i$
and that the two rows $i$ are identical in columns less than $j$,
but that the $i,j$ position in $(A\otimes_3B)\otimes_1(C\otimes_3D)$
is less than the corresponding position in
$(A\otimes_1C)\otimes_3(B\otimes_1D).$  Thus the existence of the
required inequality is shown.

The existence of $\eta^{12}$ is due to the fact
that we are ordering the matrices by giving precedence to the rows.
The two four-fold products can be seen as two alternate operations
on a single large matrix $M$. This matrix is constructed by arranging
$A,B,C,D$ with added zeroes so that $(A\otimes_1C)\otimes_2(B\otimes_1D)$
is the result of first sorting each column vertically, greater
values at the top, and then each row horizontally, greater values
to the left, while $(A\otimes_2B)\otimes_1(C\otimes_2D)$ is achieved
by sorting horizontally first and then vertically.  Recall that in
the ordering rows are given precedence over columns.  Here is an
illustration of the inequality, showing the process of constructing
the large matrix.
\[
A = \left[\begin{array}{ccc} 3&3&2 \\ 1&1 \end{array} \right] B = \left[\begin{array}{c} 9 \\9 \\9 \end{array} \right]
C = \left[\begin{array}{c} 2 \\ 1 \end{array} \right] D= \left[\begin{array}{c} 5 \end{array} \right]
\]
\[
M = \left[\begin{array}{cccc} 3&3&2&9 \\ 1&1&0&9 \\ 0&0&0&9 \\ 2&0&0&5 \\ 1&0&0&0 \end{array} \right]
\]
\[
(A\otimes_2B)\otimes_1(C\otimes_2D) = \left[\begin{array}{cccc} 9&3&3&2 \\ 9&2&1 \\ 9&1 \\ 5 \\ 1 \end{array} \right]
<  \left[\begin{array}{cccc} 9&3&3&2 \\ 9&2&1 \\ 9&1 \\ 5&1 \end{array} \right]= (A\otimes_1C)\otimes_2(B\otimes_1D)
\]
The proof that this inequality always holds requires the following two lemmas.

\begin{lemma}\label{minmax}
For two sequences of $n$ elements each, the first given by $a_1
\dots a_n$ and the second by $b_1 \dots b_n$, then considering pairs
of elements $a_{\sigma(i)}$ and $b_{\tau(i)}$ for permutations
$\sigma, \tau \in S_n$, we have the following inequality:
\[
\max(\min(a_{\sigma(1)},b_{\tau(1)}),\dots,\min(a_{\sigma(n)},b_{\tau(n)})) \le \min(\max(a_1,\dots,a_n), \max(b_1, \dots, b_n)).
\]
\end{lemma}
This is true since for $i=1\dots n$ we have $a_i \le \max(a_1,\dots,a_n)$
and $b_i \le \max(b_1, \dots, b_n).$ Therefore
$\min(a_{\sigma(i)}, b_{\tau(i)})
    \le \min(\max(a_1,\dots,a_n), \max(b_1,  \dots, b_n))$
and the inequality follows.

\begin{lemma}\label{matrixsort}
For a given finite matrix $M$ with $n$ rows, we claim that first
ordering each row (greater to lesser) and then sorting each resulting
column gives a final result that is lexicographically less than or
equal to the final result of sorting each column of $M$ and then
each row.
\end{lemma}

This is seen by a chain of inequalities that each correspond to a
single step in a parallel bubble sorting of the rows of $M.$ Consider
the final result of sorting each column vertically and then each
row. We gradually evolve this into the reverse procedure by performing
a series of steps, each of which begins by comparing two adjacent
columns in the current stage of the evolution. The step consists
of  switches that insure each horizontal pair in the columns is in
order, i.e. switching the positions of the two elements in each row
only if the one in the left column is smaller than the one in the
right.  We call this a parallel switch, or just a switch.  The
result of taking the switched matrix and vertically sorting its
columns and then horizontally sorting its rows will be shown to be
lexicographically less than or equal to the result of vertically
sorting columns and then horizontally sorting rows before the
parallel switch. The entire series of steps together constitute
sorting each row of $M.$ Since after vertically sorting a matrix
which began with sorted rows the new rows still remain sorted, then
at the end of the evolution we are indeed doing the reverse procedure;
that is sorting horizontally first and then vertically.

For a single step in the parallel bubble sort, we claim that after
the parallel switch and then vertical sorting of the two adjacent
columns the pairs in each resulting row will be either all identical
to those in the result of vertically sorting the unswitched columns,
or there will be a first row $k$ in which the pair in the switched
version of the columns consists of one element equal to one element
of the corresponding pair in the unswitched version and one element
less than the other element in the unswitched version.

Since no other columns are changed at this step, then this will
imply that after vertically sorting the other columns and then all
the rows in both matrices, the two resulting matrices will be
identical  or just identical up to the  $k^{th}$ row, where the
switched matrix will be lexicographically less than the unswitched.

The claim for two columns follows from repeated application of
Lemma~\ref{minmax}. Let the two columns be $a_1 \dots a_n$ and
$b_1 \dots b_n$ After the parallel switching, the left column holds the
max of each pair and the right the min.  Vertical sorting moves the
max of each column to the top row, and leaves all the new rows (of
two elements each) still sorted  left to right. Located in the left
position of the new top row is
\[ \max(\max(a_1,b_1), \dots ,\max(a_n,b_n))
    = \max(\max(a_1, \dots ,a_n),\max(b_1, \dots ,b_n))\]
the latter of which is the in the top row of the vertically sorted
unswitched columns. The right position in the top row of the switched
columns is \[\max(\min(a_1,b_1),\dots,\min(a_n,b_n)),\] which is
less than or equal to
\[\min(\max(a_1,\dots,a_n),\max(b_1,\dots,b_n)),\]
the other element in the top row of the vertically sorted unswitched
columns, by the preceeding lemma (with trivial permutations). If
less than, then we are done. If equal then we note  that the remaining
rows $2\dots n$ contain the same collection of elements $a_i$ and
$b_i$ in both the switched and unswitched columns, i.e.  we may
assume that in vertically sorting either version we moved $a_j$ and
$b_l$ to the top row.  Note that since the rows in the switched
version are sorted, $\max(a_l,b_l) \ge \min(a_j,b_j)$ and
$\max(a_j,b_j) \ge \min(a_l,b_l).$ Thus the $\max(a_l,b_j)$ will
always be in the left column and $\min(a_l,b_j)$ in the right.

Then the second row of the vertically sorted switched pair of columns is
\[\max(\max(a_1,b_1),\dots,\widehat{\max(a_j,b_j)},\dots,
    \widehat{\max(a_l,b_l)},\dots,\max(a_n,b_n),\max(a_l,b_j))\]
in the first position and
\[\max(\min(a_1,b_1),\dots,\widehat{\min(a_j,b_j)},\dots,
    \widehat{\min(a_l,b_l)},\dots,\min(a_n,b_n),\min(a_l,b_j))\]
in the second position, where the hats indicate missing elements.
Whereas the second row of the vertically sorted unswitched columns
is made up of
\[\max(\max(a_1,\dots,\hat{a_j},\dots,a_n))\text{ and }
    \max(\max(b_1,\dots,\hat{b_l},\dots,b_n)).\]
Thus the left position in the second row of the switched version
is the same value as one of the elements in the second row of the
unswitched vertically sorted columns.  By Lemma~\ref{minmax} with
the evident permutations, the right position in the second row is
less than or equal to the other element in the second row of the
unswitched vertically sorted columns.  If less than, then we are
done, if equal then the process continues. If the $1^{st}$ through
$(n-1)^{st}$ rows of the switched and unswitched columns contain
the same values after vertical sorting, then so do the $n\th$ rows.

This completes the proof of the lemma. By applying it to the large
matrix $M$ constructed of the four matrices $A,B,C,D$ as described
above, we have the proof of the theorem.
\end{proof}

Now we define the general $n$-fold monoidal category of $n$-dimensional
Young diagrams.  The proof of the theorem for three dimensions plays
an important role in the general theorem,  since each interchanger
involves two products.  Once we have decided to represent Young
diagrams of higher dimension by arrays of natural numbers which
decrease in each index, it is clear that each interchanger will
either involve directly two of the indices of the array or one index
as well as pointwise addition.

\begin{definition} The category of $n$-dimensional Young diagrams
consists of
\begin{enumerate}
\item Objects $A_{i_1 i_2\dots i_{n-1}}$, finitely nonzero
$n$-dimensional arrays of nonnegative integers which are monotone
decreasing in each index, and
\item Morphisms the order relations in the lexicographic ordering
with precedence given to lesser indices.
\end{enumerate}
\end{definition}

There are $n$ ways to take a product of two $n$-dimensional Young
diagrams, which we visualize as arrays of natural numbers in $n-1$
dimensions. The products correspond to merging, i.e. concatenating
and then sorting, in each of the $n-1$ possible directions, as well
as pointwise addition as $\otimes_n$. The order of products is the
reverse of the order of the indices. That is, for $k=1\dots n-1,$
$\otimes_k$ is merging in the direction of the index $i_{n-k}.$

\begin{theorem}
The category of $n$-dimensional Young diagrams with the above
products constitutes an $n$-fold monoidal category.
\end{theorem}

\begin{proof}
We must show the existence of the interchangers $\eta^{jk}$ as
inequalities for $1\le j<k\le n.$ First we demonstrate the existence
of the required inequality when $k < n.$ For $A,B,C,D$ $n$-dimensional
Young diagrams seen as $(n-1)$-dimensional arrays, we let
$M_{i_1 i_2\dots i_{n-1}}$ be a large array made by concatenating $A$
and $B$ in the direction of the index $i_k$, concatenating $C$ and $D$
in the direction of the index $i_k$, and then concatenating those
two results in the direction of the index $i_j$. Zeros are added
(see above for the two dimensional array example) so that the
products
$(A\otimes_k B)\otimes_j(C\otimes_k D)$ and
$(A\otimes_j C)\otimes_k(B\otimes_j D)$
can then both be described as sorting $M_{i_1 i_2\dots i_{n-1}}$
in two directions; first $i_k$ then $i_j$ or vice versa respectively.
That the inequality holds is seen as we compare the results position
by position in the lexicographic order, i.e. reading lower indices
first.  The first differing value we come upon, say in location
$(i_1 i_2\dots i_{n-1}),$ then will necessarily be the first
difference in  the  sub-array of two dimensions in the directions
$i_j$ and $i_k$ determined by the location $(i_1 i_2\dots i_{n-1}).$
Thus by the proof of Lemma~\ref{matrixsort}, the value in
$(A\otimes_k B)\otimes_j(C\otimes_k D)$
is less than the corresponding value in
$(A\otimes_j C)\otimes_k(B\otimes_j D).$

Secondly we check the cases that have $k=n.$ We can see the four-fold
products as operations on two arrays, one made by concatenating $A$
and $C$ in the $i_j$ direction, and another made by concatenating
$B$ and $D$ in the $i_j$ direction, padded with zeroes so that
adding the two pointwise results in pointwise addition of $A$ with
$B$, and of $C$ with $D.$ Then $(A\otimes_k B)\otimes_j(C\otimes_k
D)$ is adding first and then sorting in the $i_j$ direction, while
$(A\otimes_j C)\otimes_k(B\otimes_j D)$ is the reverse process. To
see that the correct inequality holds we again compare the results
position by position in lexicographic order. The first differeing
value is also the first difference between the two coresponding 2
dimensional sub-arrays which are in the directions $i_j$ and
$i_{n-1}.$ These sub-arrays are the results of sorting and then
pointwise addition and vice versa respectively, and so by the proof
for existence of $\eta^{13}$ in Theorem~\ref{3dY} the desired result
is shown.
\end{proof}

\section{Examples of \texorpdfstring{$n$}{N}-fold operads}

The categories from Section~\ref{section:categoryexamples} give us
a domain in which we can exhibit some concrete examples of operads.
To have an operad with an element ${\cal C}(0)$ we will need to
``compactify'' by adjoining an object that is both initial and
terminal to the example categories based on ordered monoids and
sequences. This object we will denote by $\emptyset$ and the unique
maps to $\emptyset$ will be called zero maps. Composition and tensor
product with a zero map both yield a zero map.  We define all
products involving the object $\emptyset$ as an operand to be equal
to $\emptyset.$ $\emptyset$ will be designated the least object in
the total order, except for the purpose of defining the product
using max, as just stated. Thus the unique maps from $\emptyset$
will be the $\le$ relations, and therefore the only way that two
legs of a diagram will not commute is if one of them is a zero map
and the other is not.  In all the examples the composition is
associative since it is based upon ordering, so all we need check
for is the existence of that composition. Note that each of the
following examples satisfy the hypothesis of Theorem~\ref{foo} since
taking the max is a coproduct.  Thus all coproducts are certainly
included among the objects, and the max is the first of the tensor
products in the iterated monoidal structure. Also, max distributes
over each $\otimes_i$ since each product preserves the ordering.

\begin{example}
Of course ${\cal C}(j) = \emptyset$ and ${\cal C}(j) = 0$ for all
$j$ are trivially operads, where $0$ is the monoidal unit.  First
we look at the simplest interesting examples: 2-fold operads in an
ordered monoid such as $\nat$, where $\otimes_1$ is max and $\otimes_2$
is $+$. We always set ${\cal C}(0) = \emptyset$ but often only list the later terms.
Note that the zero map cannot play the role of
operad composition, since it will fail associativity. Therefore a
2-fold operad in $\nat$ is a sequence ${\cal C}(j)$ of
natural numbers which has the property that for any $j_1 \dots j_k$,
$\max({\cal C}(k), \sum{{\cal C}(j_i)}) \le {\cal C}(\sum{j_i})$
and for which ${\cal C}(1) = 0.$ This translates into saying that
for any two whole numbers $x, y$ we have that ${\cal C}(x+y) \ge
{\cal C}(x)+{\cal C}(y)$ and that ${\cal C}(1) = 0.$
 The latter
condition both satisfies the unit axioms and makes it redundant to
also insist that the sequence be monotone increasing. Perhaps the first example
that comes to mind is the Fibonacci numbers. Minimal
examples are formed by choosing a starting term or terms and then
determining each later $n^{th}$ term.  For a starting finite sequence
$0,a_2 , \dots ,a_l$ which obeys the  the axioms of a 2-fold operad
so far, the  operad ${\cal C}_{0,a_2, \dots ,a_l}$ is found by
defining terms ${\cal C}_{a_1 , \dots ,a_l}(n)$ for $n>l$ to be the
maximum of all the values of
$\max({\cal C}(k), \sum_{i=1}^k{{\cal C}(j_i)})$
where the sum of the $j_i$ is $n.$ Some basic examples are the
following sequences.

\begin{alignat*}{2}
{\cal C}_{0,1} &= (0,1,1,2,2,3,3,\dots),
    &\qquad {\cal C}_{0,0,1}&=(0,0,1,1,1,2,2,2,3,3,3,\dots) \\
{\cal C}_{0,2} &= (0,2,2,4,4,6,6,\dots),
    &\qquad {\cal C}_{0,0,2} &= (0,0,2,2,2,4,4,4,6,6,6,\dots)
\end{alignat*}
and
\[
{\cal C}_{0,1,2,4,8} = (0,1,2,4,8,8,9,10,12,16,16,17,18,20,24,\dots).
\]
It is clear that the growth of these sequences oscillates
around linear growth in a predictable way.
\end{example}

\begin{theorem}
If ``arbitrary'' starting terms $0,a_2, \dots ,a_k \in \nat$ are given
(themselves of course obeying the axioms of a 2-fold operad), then
the $n^{th}$ term of the 2-fold operad ${\cal C}_{0,a_2, \dots ,a_k}$ in $\nat$ is given by
\[ a_n  = a_q + pa_k \text{ where } n = pk + q,
    \text{ for } p \in \nat, 0\le q <k.\]
\end{theorem}

\begin{proof}
We need to show that $\max(a_l, \sum_{i=1}^l{a_{j_i}})$ is always
less than or equal to $a_q + pa_k$, where $n = pk + q$, for
$p \in \nat$, $0\le q <k$. We need only consider the cases in which
$l < n.$ Since $a_l$ is always included as one of the  $a_{j_i},$ we
need to show only that $\sum_{i=1}^l{a_{j_i}}$ is always less than
or equal to $a_q + pa_k$ where the sum of the $j_i$ is $n.$ This
follows by strong induction. The base case holds by definition. Let
$j_i = p_ik+q_i$ for $p_i \in \nat$, and $0\le q_i <k.$ Then
$\sum{q_i} = n - k\sum{p_i}  = pk+q -k\sum{p_i} < n .$ Thus we have,
by the growth property of operads and by induction:
\begin{align*}
\sum_{i=1}^l{a_{j_i}} &=\sum{a_{q_i}} + a_k\sum{p_i}    \\
&\le a_{(k(p-\sum{p_i})+q)} + a_k\sum{p_i}      \\
&=a_q + (p-\sum{p_i})a_k + a_k\sum{p_i}         \\
&=a_q+ pa_k.\qedhere
\end{align*}
\end{proof}

\begin{example}\label{bee}
Consider the 3-fold monoidal category  $\seq(\nat,+)$ of lexicographically
ordered finitely nonzero sequences of the natural numbers (here we
use $\nat$ considered as an example of an ordered monoid), with
products $\otimes_1$ the lexicographic max , $\otimes_2$ the
concatenation and $\otimes_3$ the pointwise addition.  An example
of a 2-fold operad in  $\seq(\nat,+)$ that is not a 3-fold operad
is the following:

Let ${\cal B}(0) = \emptyset$ and
let ${\cal B}(j)_i = 1$ for $i<j~,~0$ otherwise. We can picture
these as follows:
\[
{\cal B}(1) = \xymatrix@W=1.2pc @H=1.2pc @R=1.2pc @C=1.2pc {*=0{}\ar@{-}[r] &*=0{}\ar@{-}[r] &*=0{}}~,~~
{\cal B}(2) = \xymatrix@W=.5pc @H=.5pc @R=0pc @C=0pc @*[F-]{~}~,~~
{\cal B}(3) = \xymatrix@W=.5pc @H=.5pc @R=0pc @C=0pc @*[F-]{~&~}~,~~
{\cal B}(4) = \xymatrix@W=.5pc @H=.5pc @R=0pc @C=0pc @*[F-]{~&~&~}~,~~
{\cal B}(5) = \xymatrix@W=.5pc @H=.5pc @R=0pc @C=0pc @*[F-]{~&~&~&~}~,~~ \dots
\]
This is a 2-fold operad, with respect to the lexicographic max and
concatenation.  For instance the instance of composition
$\gamma^{12}:{\cal B}(3)\otimes_1({\cal B}(2)\otimes_2{\cal B}(1)
    \otimes_2{\cal B}(3)) \to {\cal B}(6)$
appears as the relation:
\[
\xymatrix@W=1.2pc @H=1.2pc @R=1.2pc @C=1.2pc {*=0{}\ar@{-}[r]\ar@{-}[d] &*=0{}\ar@{-}[r]\ar@{-}[d] &*=0{}\ar@{-}[r]\ar@{-}[d] &*=0{}\ar@{-}[d]\\
*=0{}\ar@{-}[r] &*=0{}\ar@{-}[r] &*=0{}\ar@{-}[r] &*=0{}\\
}
\raisebox{-0.5em}{ $<$ } \xymatrix@W=1.2pc @H=1.2pc @R=1.2pc @C=1.2pc {*=0{}\ar@{-}[r]\ar@{-}[d] &*=0{}\ar@{-}[r]\ar@{-}[d] &*=0{}\ar@{-}[r]\ar@{-}[d] &*=0{}\ar@{-}[r]\ar@{-}[d] &*=0{}\ar@{-}[r]\ar@{-}[d] &*=0{}\ar@{-}[d]\\
*=0{}\ar@{-}[r] &*=0{}\ar@{-}[r] &*=0{}\ar@{-}[r] &*=0{}\ar@{-}[r] &*=0{}\ar@{-}[r] &*=0{}}
\]
However, the relation
\[
\xymatrix@W=1.2pc @H=1.2pc @R=1.2pc @C=1.2pc {*=0{}\ar@{-}[r]\ar@{-}[d] &*=0{}\ar@{-}[r]\ar@{-}[d] &*=0{}\ar@{-}[r]\ar@{-}[d] &*=0{}\ar@{-}[r]\ar@{-}[d] &*=0{}\ar@{-}[d]\\
*=0{}\ar@{-}[r] &*=0{}\ar@{-}[r] &*=0{}\ar@{-}[r]\ar@{-}[d] &*=0{}\ar@{-}[r]\ar@{-}[d] &*=0{}\\
*=0{} &*=0{} &*=0{}\ar@{-}[r] &*=0{}}
\raisebox{-0.5em}{ $>$ } \xymatrix@W=1.2pc @H=1.2pc @R=1.2pc @C=1.2pc {*=0{}\ar@{-}[r]\ar@{-}[d] &*=0{}\ar@{-}[r]\ar@{-}[d] &*=0{}\ar@{-}[r]\ar@{-}[d] &*=0{}\ar@{-}[r]\ar@{-}[d] &*=0{}\ar@{-}[r]\ar@{-}[d] &*=0{}\ar@{-}[d]\\
*=0{}\ar@{-}[r] &*=0{}\ar@{-}[r] &*=0{}\ar@{-}[r] &*=0{}\ar@{-}[r] &*=0{}\ar@{-}[r] &*=0{}}
\]
shows that
$\gamma^{23}:{\cal B}(3)\otimes_2 ({\cal B}(1)\otimes_3{\cal B}(3)
    \otimes_3{\cal B}(2)) \to {\cal B}(6)$
does not exist, so that ${\cal B}$ is not a 3-fold operad.
\end{example}

\begin{example}
Next we give an example of a 3-fold operad in $\seq(\nat,+)$.
Let ${\cal C}(0) = \emptyset$ and
let ${\cal C}(j) = (j-1,0 \dots).$ We can picture these as follows:
\[
{\cal C}(1) = \xymatrix@W=1.2pc @H=1.2pc @R=1.2pc @C=1.2pc {*=0{}\ar@{-}[r] &*=0{}\ar@{-}[r] &*=0{}}~,~~
{\cal C}(2) = \xymatrix@W=.5pc @H=.5pc @R=0pc @C=0pc @*[F-]{~}~,~~
{\cal C}(3) = \xymatrix@W=.5pc @H=.5pc @R=0pc @C=0pc @*[F-]{~\\~}~,~~
{\cal C}(4) = \xymatrix@W=.5pc @H=.5pc @R=0pc @C=0pc @*[F-]{~\\~\\~}~,~~
{\cal C}(5) = \xymatrix@W=.5pc @H=.5pc @R=0pc @C=0pc @*[F-]{~\\~\\~\\~}~,~~\dots
\]

First we note that the operad ${\cal C}$ just given is a 3-fold
operad since we have that the
$\gamma^{23}:{\cal C}(k)\otimes_2({\cal C}(j_i)\otimes_3\dots
    \otimes_3{\cal C}(j_k)) \to {\cal C}(j)$
exists.  For instance
$\gamma^{23}:{\cal C}(3)\otimes_2 ({\cal C}(1)\otimes_3
    {\cal C}(3)\otimes_3{\cal C}(2)) \to {\cal C}(6)$
appears as the relation
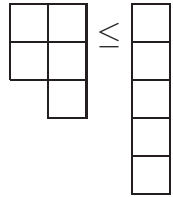
\[
\xymatrix@W=1.2pc @H=1.2pc @R=1.2pc @C=1.2pc {*=0{}\ar@{-}[r]\ar@{-}[d] &*=0{}\ar@{-}[r]\ar@{-}[d] &*=0{}\ar@{-}[d]\\
*=0{}\ar@{-}[r]\ar@{-}[d] &*=0{}\ar@{-}[r]\ar@{-}[d] &*=0{}\ar@{-}[d]\\
*=0{}\ar@{-}[r] &*=0{}\ar@{-}[r]\ar@{-}[d] &*=0{}\ar@{-}[d]\\
*=0{} &*=0{}\ar@{-}[r] &*=0{}}
\raisebox{-1em}{ $\le$ }
\xymatrix@W=1.2pc @H=1.2pc @R=1.2pc @C=1.2pc {*=0{}\ar@{-}[r]\ar@{-}[d] &*=0{}\ar@{-}[d]\\
*=0{}\ar@{-}[r]\ar@{-}[d] &*=0{}\ar@{-}[d]\\
*=0{}\ar@{-}[r]\ar@{-}[d] &*=0{}\ar@{-}[d]\\
*=0{}\ar@{-}[r]\ar@{-}[d] &*=0{}\ar@{-}[d]\\
*=0{}\ar@{-}[r]\ar@{-}[d] &*=0{}\ar@{-}[d]\\
*=0{}\ar@{-}[r] &*=0{}}
\]

Then we remark that as expected the composition
$\gamma^{12}:{\cal C}(k)\otimes_1({\cal C}(j_i)\otimes_2\dots
    \otimes_2{\cal C}(j_k)) \to {\cal C}(j)$
also exists.  For instance
$\gamma^{12}:{\cal C}(3)\otimes_1 ({\cal C}(1)\otimes_2{\cal C}(2)
    \otimes_2{\cal C}(3)) \to {\cal C}(6)$
appears as the relation
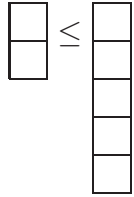
\[
\xymatrix@W=1.2pc @H=1.2pc @R=1.2pc @C=1.2pc {*=0{}\ar@{-}[r]\ar@{-}[d] &*=0{}\ar@{-}[d]\\
*=0{}\ar@{-}[r]\ar@{-}[d] &*=0{}\ar@{-}[d]\\
*=0{}\ar@{-}[r] &*=0{}}
\raisebox{-1em}{ $\le$ }
\xymatrix@W=1.2pc @H=1.2pc @R=1.2pc @C=1.2pc {*=0{}\ar@{-}[r]\ar@{-}[d] &*=0{}\ar@{-}[d]\\
*=0{}\ar@{-}[r]\ar@{-}[d] &*=0{}\ar@{-}[d]\\
*=0{}\ar@{-}[r]\ar@{-}[d] &*=0{}\ar@{-}[d]\\
*=0{}\ar@{-}[r]\ar@{-}[d] &*=0{}\ar@{-}[d]\\
*=0{}\ar@{-}[r]\ar@{-}[d] &*=0{}\ar@{-}[d]\\
*=0{}\ar@{-}[r] &*=0{}}
\]
\end{example}

\begin{example}
Now we consider some products of the previous two described operads
in $\seq(\nat,+)$. We expect ${\cal B} {\otimes'}{\cal C}$ given by
$({\cal B} {\otimes'}{\cal C})(j) = {\cal B}(j) \otimes_3{\cal C}(j)$
to be a 2-fold operad and it is. It appears thus:
\[
\emptyset~,~~ \xymatrix@W=1.2pc @H=1.2pc @R=1.2pc @C=1.2pc {*=0{}\ar@{-}[r] &*=0{}\ar@{-}[r] &*=0{}}~,~~
\xymatrix@W=.5pc @H=.5pc @R=0pc @C=0pc @*[F-]{~\\~}~,~~
\xymatrix@W=.5pc @H=.5pc @R=0pc @C=0pc @*[F-]{~&~\\~\\~}~,~~
\xymatrix@W=.5pc @H=.5pc @R=0pc @C=0pc @*[F-]{~&~&~\\~\\~\\~}~,~~
\xymatrix@W=.5pc @H=.5pc @R=0pc @C=0pc @*[F-]{~&~&~&~\\~\\~\\~\\~}~,~~
\dots
\]

We demonstrate the tightness of the existence of products of operads
by pointing out that $D(j) = {\cal B}(j) \otimes_2{\cal C}(j)$ does
not form an operad. We leave it to the reader to demonstrate this
fact.
\end{example}

Now we pass to the categories of Young diagrams in which the
interesting products are given by horizontal and vertical stacking.
It is important that we do not restrict the morphisms to those
between diagrams of the same total number of blocks in order that
all  the operad compositions exist.

\begin{theorem}\label{suff}
A sequence of Young diagrams ${\cal C}(n),~ n\in\nat$, in the
category $\modseq(\nat,+)$, is a 2-fold operad if $C(0) = \emptyset$
and for $n \ge 1$, $h({\cal C}(n)) = f(n)$ where
$f\colon\nat\to\nat$ is a function such that $f(1) = 0$ and
$f(i+j)> f(i) + f(j)$.
\end{theorem}

\begin{proof}
These conditions are not necessary, but they are sufficient since
the first implies that ${\cal C}(1) =0$ which shows that the unit
conditions are satisfied; and the second implies that the maps
$\gamma$ exist. We see existence of $\gamma^{12}$ since for $j_i > 0$,
$h({\cal C}(k)\otimes_1({\cal C}(j_1)\otimes_2\dots
    \otimes_2{\cal C}(j_k))) =\max(f(k),\max(f(j_i))) \le f(j).$
We have existence of $\gamma^{13}$  and $\gamma^{23}$ since
$\max(f(k), \sum{f(j_i)}) \le f(j)$.
\end{proof}

\begin{example}
Examples of $f$ include  $(x-1)P(x)$ where $P$ is a polynomial with
coefficients in $\nat$. This is easy to show since then $P$ will
be monotone increasing for $x\ge 1$ and thus $(i+j-1)P(i+j) =
(i-1)P(i+j)+jP(i+j) \ge (i-1)P(i)+jP(j) - P(j).$ By this argument
we can also use any $f= (x-1)g(x)$ where $g\colon\nat\to\nat$
is monotone increasing for $x \ge 1$.
\end{example}

For a specific example with a handy picture that also illustrates
again the nontrivial use of the interchange $\eta$ we simply let
$f = x-1.$ Then we have to actually describe the elements of
$\modseq(\nat)$ that make up the operad. One nice choice is the
operad ${\cal C}$ where ${\cal C}(n) = \{n-1,n-1,...,n-1\},$ the
$(n-1)\times (n-1)$ square Young diagram.
\[ {\cal C}(1) = 0,
\ {\cal C}(2) = \xymatrix@W=.6pc @H=.6pc @R=0pc @C=0pc @*[F-]{~},
\ {\cal C}(3) = \xymatrix@W=.6pc @H=.6pc @R=0pc @C=0pc @*[F-]{~&~\\~&~}\ ,
    \ldots
\]
For instance $\gamma^{23}:{\cal C}(3)\otimes_2 ({\cal C}(1)\otimes_3{\cal C}(3)\otimes_3{\cal C}(2)) \to {\cal C}(6)$
appears as the relation
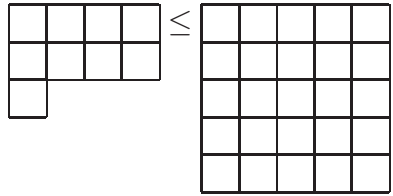
\[
\xymatrix@W=.6pc @H=.6pc @R=0pc @C=0pc @*[F-]{~&~&~&~\\~&~&~&~\\~} \le \xymatrix@W=.6pc @H=.6pc @R=0pc @C=0pc @*[F-]{~&~&~&~&~\\~&~&~&~&~\\~&~&~&~&~\\~&~&~&~&~\\~&~&~&~&~}
\]
An instance of the associativity diagram with upper left position
${\cal C}(2)\otimes_2({\cal C}(3)\otimes_3{\cal C}(2))\otimes_2
({\cal C}(2)\otimes_3{\cal C}(2)\otimes_3{\cal C}(4)\otimes_3{\cal C}(5)\otimes_3{\cal C}(3))$ is as follows:

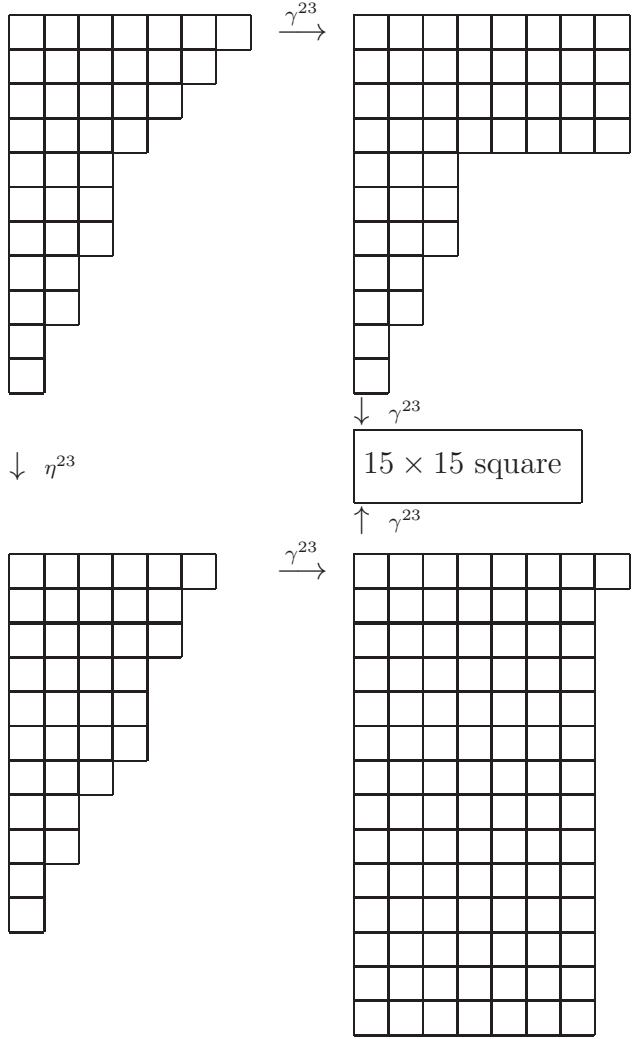
\begin{tabular}{llll}
&$\xymatrix@W=.5pc @H=.5pc @R=0pc @C=0pc @*[F-]{~&~&~&~&~&~&~\\~&~&~&~&~&~\\~&~&~&~&~\\~&~&~&~\\~&~&~\\~&~&~\\~&~&~\\~&~\\~&~\\~\\~}$&$\buildrel\gamma^{23}\over\longrightarrow$&$\xymatrix@W=.5pc @H=.5pc @R=0pc @C=0pc @*[F-]{~&~&~&~&~&~&~&~\\~&~&~&~&~&~&~&~\\~&~&~&~&~&~&~&~\\~&~&~&~&~&~&~&~\\~&~&~\\~&~&~\\~&~&~\\~&~\\~&~\\~\\~}$ \\
&&&$\downarrow~{\scriptstyle\gamma^{23}}$\\
&$\downarrow~{\scriptstyle\eta^{23}}$&&$\xymatrix@W=1.7pc @H=1.7pc @R=0pc @C=0pc @*[F-]{15\times 15 \text{ square }}$\\
&&&$\uparrow~{\scriptstyle\gamma^{23}}$\\
&$\xymatrix@W=.5pc @H=.5pc @R=0pc @C=0pc @*[F-]{~&~&~&~&~&~\\~&~&~&~&~\\~&~&~&~&~\\~&~&~&~\\~&~&~&~\\~&~&~&~\\~&~&~\\~&~\\~&~\\~\\~}$&$\buildrel\gamma^{23}\over\longrightarrow$&$\xymatrix@W=.5pc @H=.5pc @R=0pc @C=0pc @*[F-]{~&~&~&~&~&~&~&~\\~&~&~&~&~&~&~\\~&~&~&~&~&~&~\\~&~&~&~&~&~&~\\~&~&~&~&~&~&~\\~&~&~&~&~&~&~\\~&~&~&~&~&~&~\\~&~&~&~&~&~&~\\~&~&~&~&~&~&~\\~&~&~&~&~&~&~\\~&~&~&~&~&~&~\\~&~&~&~&~&~&~\\~&~&~&~&~&~&~\\~&~&~&~&~&~&~}$\\
&\\
\end{tabular}

\begin{example}
Again we note that the conditions in Theorem~\ref{suff} are not
neccessary ones. In fact, given any Young diagram $B$ we can construct
a unique operad that is minimal in each term with respect to ordering
of the diagrams.
\end{example}
\begin{definition}
The 2-fold operad in the category of Young diagrams generated by a
Young diagram $B$ is  denoted by ${\cal C}_B$ and defined as follows:
${\cal C}_B(1) =0$ and ${\cal C}_B(2) = B.$ Each sucessive term is
defined to be the lexicographic maximum of all the products of prior
terms which compose to the term in question; for $n>2$ and over
$\sum{j_i} = n$:
\[ {\cal C}_B(n) = \max\{{\cal C}_B(k)\otimes_2({\cal C}_B(j_1)
    \otimes_3\dots\otimes_3{\cal C}_B(j_k))\}. \]
\end{definition}

\begin{theorem}
If a Young diagram $B$ consists of total number of blocks $q$, then the term ${\cal C}_B(n)$ of the operad generated by $B$ consists of $q(n-1)$ blocks.
\end{theorem}

\begin{proof}
The proof is by strong induction. The number of blocks is given for
${\cal C}_B(1)$ and ${\cal C}_B(2)$. Since the definition is in
terms of a maximum over composable products, if the number of blocks
in each piece of any such a product is assumed by induction to be
respectively $q(k-1)$, and  $q(j_1 -1) \dots q(j_k -1)$, then the
total number of blocks in each product (and thus the maximum) is
$q(n-1)$ since $\sum{j_i} = n$.
\end{proof}

Here are the first few terms of the operad thus generated by
$B = 
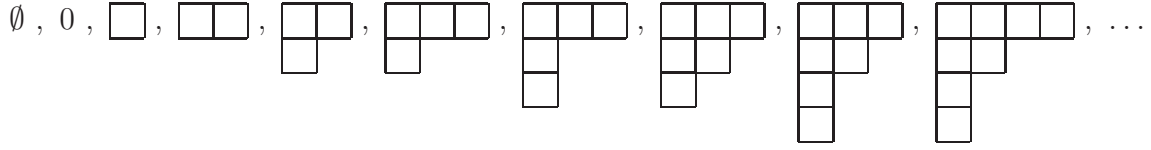
\xymatrix@W=.5pc @H=.5pc @R=0pc @C=0pc @*[F-]{~}~~$.
\[
\emptyset~,~~0~,~~\xymatrix@W=.5pc @H=.5pc @R=0pc @C=0pc @*[F-]{~}~,~~
\xymatrix@W=.5pc @H=.5pc @R=0pc @C=0pc @*[F-]{~&~}~,~~
\xymatrix@W=.5pc @H=.5pc @R=0pc @C=0pc @*[F-]{~&~\\~}~,~~
\xymatrix@W=.5pc @H=.5pc @R=0pc @C=0pc @*[F-]{~&~&~\\~}~,~~
\xymatrix@W=.5pc @H=.5pc @R=0pc @C=0pc @*[F-]{~&~&~\\~\\~}~,~~
\xymatrix@W=.5pc @H=.5pc @R=0pc @C=0pc @*[F-]{~&~&~\\~&~\\~}~,~~
\xymatrix@W=.5pc @H=.5pc @R=0pc @C=0pc @*[F-]{~&~&~\\~&~\\~\\~}~,~~
\xymatrix@W=.5pc @H=.5pc @R=0pc @C=0pc @*[F-]{~&~&~&~\\~&~\\~\\~}~,~~\dots
\]
Note that  height of any given column grows linearly, but that the
length of any row grows logarithmically.
\begin{theorem}
The minimal operad ${\cal C}_{\xymatrix@W=.35pc @H=.35pc @R=.35pc @C=.35pc {*=0{}\ar@{-}[r]\ar@{-}[d] &*=0{}\ar@{-}[d]\\*=0{}\ar@{-}[r]&*=0{}}}$ of Young diagrams which begins with ${\cal C}_{\xymatrix@W=.35pc @H=.35pc @R=.35pc @C=.35pc {*=0{}\ar@{-}[r]\ar@{-}[d] &*=0{}\ar@{-}[d]\\*=0{}\ar@{-}[r]&*=0{}}}(1) =0$
and ${\cal C}_{\xymatrix@W=.35pc @H=.35pc @R=.35pc @C=.35pc {*=0{}\ar@{-}[r]\ar@{-}[d] &*=0{}\ar@{-}[d]\\*=0{}\ar@{-}[r]&*=0{}}}(2) = \xymatrix@W=.5pc @H=.5pc @R=0pc @C=0pc @*[F-]{~}~~$,
has terms ${\cal C}_{\xymatrix@W=.35pc @H=.35pc @R=.35pc @C=.35pc {*=0{}\ar@{-}[r]\ar@{-}[d] &*=0{}\ar@{-}[d]\\*=0{}\ar@{-}[r]&*=0{}}}(n)$ that are built of $n-1$ blocks each, and whose monotone decreasing sequence representation
is given by the formula
\[{\cal C}_{\xymatrix@W=.35pc @H=.35pc @R=.35pc @C=.35pc {*=0{}\ar@{-}[r]\ar@{-}[d] &*=0{}\ar@{-}[d]\\*=0{}\ar@{-}[r]&*=0{}}}(n)_k = \operatorname{Round}\left(n/2^k\right); k=1,2,\dots\]
where rounding is done to the nearest integer and .5 is rounded to zero.
\end{theorem}

\begin{proof}
The proof of the formula for the column heights is by way of first showing that
each term in
${\cal C}_{\xymatrix@W=.35pc @H=.35pc @R=.35pc @C=.35pc {*=0{}\ar@{-}[r]\ar@{-}[d] &*=0{}\ar@{-}[d]\\*=0{}\ar@{-}[r]&*=0{}}}$
can be built canonically as follows:
\[
{\cal C}_{\xymatrix@W=.35pc @H=.35pc @R=.35pc @C=.35pc {*=0{}\ar@{-}[r]\ar@{-}[d] &*=0{}\ar@{-}[d]\\*=0{}\ar@{-}[r]&*=0{}}}(n)
= \begin{array}{cc}  & \lceil\frac{n}{2}\rceil \\
{\cal C}_{\xymatrix@W=.35pc @H=.35pc @R=.35pc @C=.35pc {*=0{}\ar@{-}[r]\ar@{-}[d] &*=0{}\ar@{-}[d]\\*=0{}\ar@{-}[r]&*=0{}}}(\lceil\frac{n}{2}\rceil)\otimes_2(& \overbrace{\underbrace{{\cal C}_{\xymatrix@W=.35pc @H=.35pc @R=.35pc @C=.35pc {*=0{}\ar@{-}[r]\ar@{-}[d] &*=0{}\ar@{-}[d]\\*=0{}\ar@{-}[r]&*=0{}}}(2)\otimes_3 \dots \otimes_3{\cal C}_{\xymatrix@W=.35pc @H=.35pc @R=.35pc @C=.35pc {*=0{}\ar@{-}[r]\ar@{-}[d] &*=0{}\ar@{-}[d]\\*=0{}\ar@{-}[r]&*=0{}}}(2)} \otimes_3{\cal C}_{\xymatrix@W=.35pc @H=.35pc @R=.35pc @C=.35pc {*=0{}\ar@{-}[r]\ar@{-}[d] &*=0{}\ar@{-}[d]\\*=0{}\ar@{-}[r]&*=0{}}}(1)}~~) \\
&\lfloor\frac{n}{2}\rfloor \end{array}
\]
We must demonstrate that the
maximum of all ${\cal C}_{\xymatrix@W=.35pc @H=.35pc @R=.35pc @C=.35pc {*=0{}\ar@{-}[r]\ar@{-}[d] &*=0{}\ar@{-}[d]\\*=0{}\ar@{-}[r]&*=0{}}}(k)\otimes_2({\cal C}_{\xymatrix@W=.35pc @H=.35pc @R=.35pc @C=.35pc {*=0{}\ar@{-}[r]\ar@{-}[d] &*=0{}\ar@{-}[d]\\*=0{}\ar@{-}[r]&*=0{}}}(j_1)\otimes_3\dots\otimes_3{\cal C}_{\xymatrix@W=.35pc @H=.35pc @R=.35pc @C=.35pc {*=0{}\ar@{-}[r]\ar@{-}[d] &*=0{}\ar@{-}[d]\\*=0{}\ar@{-}[r]&*=0{}}}(j_k))$
where $\sum{j_i} = n$
is precisely given by the above canonical construction. We make the assumption (of strong induction)
that this holds for terms before the $n^{th}$ term, and check for the inequality
${\cal C}_{\xymatrix@W=.35pc @H=.35pc @R=.35pc @C=.35pc {*=0{}\ar@{-}[r]\ar@{-}[d] &*=0{}\ar@{-}[d]\\*=0{}\ar@{-}[r]&*=0{}}}(k)\otimes_2({\cal C}_{\xymatrix@W=.35pc @H=.35pc @R=.35pc @C=.35pc {*=0{}\ar@{-}[r]\ar@{-}[d] &*=0{}\ar@{-}[d]\\*=0{}\ar@{-}[r]&*=0{}}}(j_1)\otimes_3\dots\otimes_3{\cal C}_{\xymatrix@W=.35pc @H=.35pc @R=.35pc @C=.35pc {*=0{}\ar@{-}[r]\ar@{-}[d] &*=0{}\ar@{-}[d]\\*=0{}\ar@{-}[r]&*=0{}}}(j_k))$
less than or equal to the canonical construction.
The case in which there are only 0 or 1 odd integers among the $j_k$'s is
directly observed using the strong induction. If there are two or more
odd integers among the $j_k$'s  and the first column of the diagram they help determine
is greater than or equal to the first column of
${\cal C}_{\xymatrix@W=.35pc @H=.35pc @R=.35pc @C=.35pc {*=0{}\ar@{-}[r]\ar@{-}[d] &*=0{}\ar@{-}[d]\\*=0{}\ar@{-}[r]&*=0{}}}(k)$
then the inequality holds by induction on the size of the first column.
If there are two or more
odd integers among the $j_k$'s  and the first column of the diagram they help determine
is less than the first column of
${\cal C}_{\xymatrix@W=.35pc @H=.35pc @R=.35pc @C=.35pc {*=0{}\ar@{-}[r]\ar@{-}[d] &*=0{}\ar@{-}[d]\\*=0{}\ar@{-}[r]&*=0{}}}(k)$
then we check the sub-cases $n$ odd and $n$ even. For $n$ even the result is seen directly, and for
$n$ odd we again rely on induction.
\end{proof}

For comparison to the previous example of the operad with square terms,  the
instance of the associativity diagram with upper left position
${\cal C}_{\xymatrix@W=.35pc @H=.35pc @R=.35pc @C=.35pc {*=0{}\ar@{-}[r]\ar@{-}[d] &*=0{}\ar@{-}[d]\\*=0{}\ar@{-}[r]&*=0{}}}(2)\otimes_2({\cal C}_{\xymatrix@W=.35pc @H=.35pc @R=.35pc @C=.35pc {*=0{}\ar@{-}[r]\ar@{-}[d] &*=0{}\ar@{-}[d]\\*=0{}\ar@{-}[r]&*=0{}}}(3)\otimes_3{\cal C}_{\xymatrix@W=.35pc @H=.35pc @R=.35pc @C=.35pc {*=0{}\ar@{-}[r]\ar@{-}[d] &*=0{}\ar@{-}[d]\\*=0{}\ar@{-}[r]&*=0{}}}(2))\otimes_2
({\cal C}_{\xymatrix@W=.35pc @H=.35pc @R=.35pc @C=.35pc {*=0{}\ar@{-}[r]\ar@{-}[d] &*=0{}\ar@{-}[d]\\*=0{}\ar@{-}[r]&*=0{}}}(2)\otimes_3{\cal C}_{\xymatrix@W=.35pc @H=.35pc @R=.35pc @C=.35pc {*=0{}\ar@{-}[r]\ar@{-}[d] &*=0{}\ar@{-}[d]\\*=0{}\ar@{-}[r]&*=0{}}}(2)\otimes_3{\cal C}_{\xymatrix@W=.35pc @H=.35pc @R=.35pc @C=.35pc {*=0{}\ar@{-}[r]\ar@{-}[d] &*=0{}\ar@{-}[d]\\*=0{}\ar@{-}[r]&*=0{}}}(4)\otimes_3{\cal C}_{\xymatrix@W=.35pc @H=.35pc @R=.35pc @C=.35pc {*=0{}\ar@{-}[r]\ar@{-}[d] &*=0{}\ar@{-}[d]\\*=0{}\ar@{-}[r]&*=0{}}}(5)\otimes_3{\cal C}_{\xymatrix@W=.35pc @H=.35pc @R=.35pc @C=.35pc {*=0{}\ar@{-}[r]\ar@{-}[d] &*=0{}\ar@{-}[d]\\*=0{}\ar@{-}[r]&*=0{}}}(3))$ is as follows:

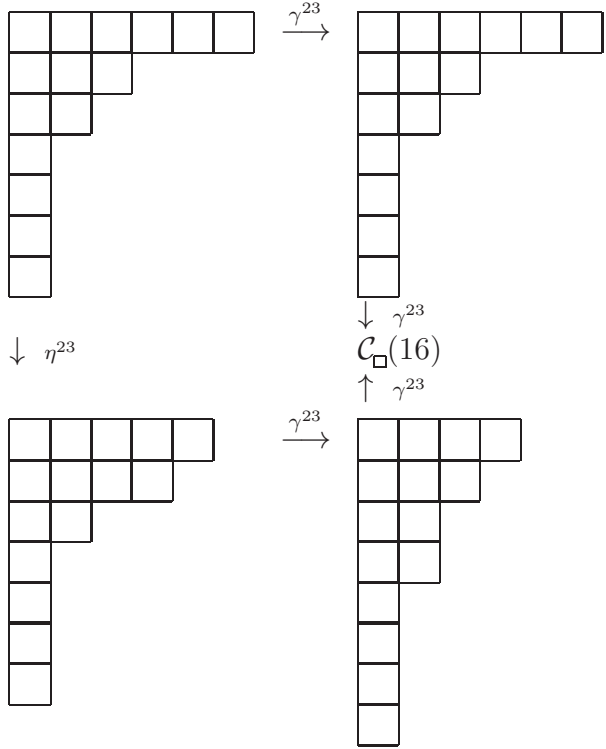
\begin{tabular}{llll}
&$\xymatrix@W=.7pc @H=.7pc @R=0pc @C=0pc @*[F-]{~&~&~&~&~&~\\~&~&~\\~&~\\~\\~\\~\\~}$&$\buildrel\gamma^{23}\over\longrightarrow$&$\xymatrix@W=.7pc @H=.7pc @R=0pc @C=0pc @*[F-]{~&~&~&~&~&~\\~&~&~\\~&~\\~\\~\\~\\~}$ \\
&&&$\downarrow~{\scriptstyle\gamma^{23}}$\\
&$\downarrow~{\scriptstyle\eta^{23}}$&&${\cal C}_{\xymatrix@W=.35pc @H=.35pc @R=.35pc @C=.35pc {*=0{}\ar@{-}[r]\ar@{-}[d] &*=0{}\ar@{-}[d]\\*=0{}\ar@{-}[r]&*=0{}}}(16)$\\
&&&$\uparrow~{\scriptstyle\gamma^{23}}$\\
&$\xymatrix@W=.7pc @H=.7pc @R=0pc @C=0pc @*[F-]{~&~&~&~&~\\~&~&~&~\\~&~\\~\\~\\~\\~}$&$\buildrel\gamma^{23}\over\longrightarrow$&$\xymatrix@W=.7pc @H=.7pc @R=0pc @C=0pc @*[F-]{~&~&~&~\\~&~&~\\~&~\\~&~\\~\\~\\~\\~}$\\
&\\
\end{tabular}

There may be interesting applications of the type of growth modeled
by operads in iterated monoidal categories.  Since the growth is
in multiple dimensions it suggests applications in the science of
allometric measurements, broadly used to refer to any $n$ characteristics
that grow in tandem. These measurements are often used in biological
sciences to try to predict values of one characteristic from others,
such as tree height from trunk diameter or crown diameter, or
skeletal mass from total body mass or dimensions, or even genomic
diversity from various geographical features.  Allometric comparisons
are often used in geology and chemistry, for instance when predicting
the growth of speleothems or crystals. There are also potential
applications to networks, where the growth of diameter or linking
distance of a network is related logarithmically to the growth in
number of nodes. In computational geometry, the number of vertices
of the convex hull of $n$ uniformly scattered points in a polygon
grows as the log of $n$.

This sort of minimal growth in the terms of the operad could be
perturbed, for example by replacing the term
$\xymatrix@W=.5pc @H=.5pc @R=0pc @C=0pc @*[F-]{~&~}\ $
in the above with the alternate term
\raisebox{0.5em}
{$\xymatrix@W=.5pc @H=.5pc @R=0pc @C=0pc @*[F-]{~\\~}\ $},
which would affect the later terms in turn. An interesting avenue
for further investigation would be the comparison of such perturbations
to determine the relative effects of a given perturbation's size
and position of occurence in the sequence.


We conclude with a description of the concepts of $n$-fold operad
algebra and of the tensor products of operad algebras.

\begin{definition}\label{opalg}
Let ${\cal C}$ be  an $n$-fold operad in ${\cal V}$. A
${\cal C}$-algebra is an object $A\in{\cal V}$ and maps
\[
\theta^{pq}:{\cal C}(j)\otimes_p(\otimes_q^j A)\to A
\]
for $n\ge q>p \ge 1$, $j\ge 0$.

\begin{enumerate}
\item Associativity: The following diagram is required to commute
for all $n\ge q>p \ge 1$, $k\ge 1$, $j_s\ge 0$ , where
$j = \sum\limits_{s=1}^k j_s$.

\[
\xymatrix{{\cal C}(k) \otimes_p ({\cal C}(j_1) \otimes_q \dots \otimes_q {\cal C}(j_k))\otimes_p(\otimes_q^j A)
\ar[rr]^>>>>>>>>>>>>{\gamma^{pq} \otimes_p \text{id}}
\ar[dd]_{\text{id} \otimes_p \eta^{pq}}
&& {\cal C}(j)\otimes_p(\otimes_q^j A)
\ar[d]_{\theta^{pq}}\\
&&A\\
{\cal C}(k) \otimes_p (({\cal C}(j_1)\otimes_p (\otimes_q^{j_1} A))\otimes_q \dots \otimes_q ({\cal C}(j_k)\otimes_p (\otimes_q^{j_k} A)))
\ar[rr]^>>>>>>>>>>>>{\text{id} \otimes_p (\otimes_q^k \theta^{pq})}
&&{\cal C}(k)\otimes_p(\otimes_q^k A)
\ar[u]^{\theta^{pq}}
}
\]

\item Units: The following diagram is required to commute for all
$n\ge q>p \ge 1$.
\[
\xymatrix{
I\otimes_p A
\ar[d]_{{\cal J}\otimes_p 1}
\ar@{=}[r]^{}
&A\\
{\cal C}(1)\otimes_p A
\ar[ur]^{\theta^{pq}}
}
\]
\end{enumerate}
\end{definition}

\begin{example}
Of course the initial object is always an algebra for every operad,
and every object is an algebra for the initial operad. For a slightly
less trivial example we turn to the height preordered category of
Remark~\ref{hpre}. Define the operad ${\cal B}(j)$ as in
Example~\ref{bee}. Then any nonzero sequence $A$ is an algebra for
this operad.
\end{example}

\begin{remark}
Let ${\cal C}$ and  ${\cal D}$ be $m$-fold operads in an $n$-fold
monoidal category.  If $A$ is an algebra of ${\cal C}$ and $B$ is
an algebra of ${\cal D}$ then $A\otimes_{i+m} B$ is an algebra for
${\cal C}\otimes'_i{\cal D}.$

That the product of $n$-fold operad algebras is itself an $n$-fold
operad algebra is easy to verify once we note that the new $\theta$
is in terms of the two old ones:
\[
\theta^{pq}_{A\otimes_{i+m} B} =
(\theta^{pq}_{A}\otimes_{i+m}\theta^{pq}_{B})\circ \eta^{p(i+m)}
    \circ (1 \otimes_p \eta^{q(i+m)})
\]
Maps of operad algebras are straightforward to describe--they are
required to preserve  structure; that is to commute with  $\theta.$
\end{remark}


\end{document}